\documentclass[preprint, fleqn,  times, 3p, 10pt]{elsarticle}

\usepackage{docmute}

\usepackage{hyperref}
\usepackage{geometry}  
\usepackage{amsbsy}  
\usepackage{mathtools}  
\usepackage{isomath}  
\usepackage{xfrac}  
\usepackage{xspace}  
\usepackage{float}  
\usepackage{graphicx}  
\usepackage{subcaption}  
\usepackage{natbib}  
\usepackage[capitalise]{cleveref}  
\usepackage{appendix}  
\usepackage{enumitem}  
\usepackage{multicol}  
\usepackage[noend]{algpseudocode}  
\usepackage{algorithm}  
\usepackage{acro}  
\usepackage{xr}  
\usepackage{bm} 
\usepackage{tabularx}  
\usepackage{etoolbox}  
\usepackage{booktabs}
\usepackage{multirow}
\usepackage{chngcntr}
\usepackage{apptools}
\usepackage{amsfonts}
\usepackage{mathrsfs}
\usepackage{amsmath} 
\usepackage{amsthm}
\usepackage{xcolor}
\usepackage{comment}
\AtAppendix{\counterwithin{lemma}{section}}

\hypersetup{colorlinks=true, allcolors=blue}

\newtheorem{theorem}{Theorem}
\newtheorem{lemma}{Lemma}
\newtheorem{proposition}{Proposition}
\newtheorem*{definition}{Definition}
 
\newtheorem{remark}{Remark}

\newcommand{\trueparam}{\theta^\dagger}
\newcommand{\truebeta}{\beta^\dagger}
\newcommand{\truesigma}{\sigma^\dagger}

\newcommand{\trueE}{\mathbb{E}_{\trueparam}}
\newcommand{\truedist}{\nu_{\trueparam}}

\newcommand{\mrm}{\mathrm} 
\newcommand{\mbf}{\mathbf}  
\newcommand{\limit}{$n \to \infty$, $\Delta_n \to 0$ and $n \Delta_n \to \infty$} 
\newcommand*{\filtration}[1]{\mathcal{F}_{t_{#1}}} 
\newcommand*{\sample}[1]{X_{{#1}}}

\newcommand{\probconv}{\xrightarrow{\mathbb{P}_{\trueparam}}}
\newcommand{\distconv}{\xrightarrow{\mathcal{L}_{\trueparam}}}

\everymath {\displaystyle}

\newcommand{\ruby}[2]{
\leavevmode
\setbox0=\hbox{#1}
\setbox1=\hbox{\tiny #2}
\ifdim\wd0>\wd1 \dimen0=\wd0 \else \dimen0=\wd1 \fi
\hbox{
\kanjiskip=0pt plus 2fil
\xkanjiskip=0pt plus 2fil
\vbox{
\hbox to \dimen0{
\small \hfil#2\hfil}
\nointerlineskip
\hbox to \dimen0{\mathstrut\hfil#1\hfil}}}}
\allowdisplaybreaks[1]

\newtheorem{condition}{Condition}

\journal{Stochastic Processes and their Applications}

\begin{document}

\begin{frontmatter}



\title{Parameter Inference for Degenerate Diffusion Processes}

\author[1]{Yuga Iguchi}
\corref{cor1} 
\ead{yuga.iguchi.21@ucl.ac.uk}
\author[1]{Alexandros Beskos} 
\author[2]{Matthew M. Graham} 

\cortext[cor1]{Corresponding author}
\affiliation[1]{organization={Department of Statistical Science, University College London}, addressline={1-19 Torrington Place},
postcode={WC1E 6BT}, city={London}, country={United Kingdom}}


\affiliation[2]{organization={Advanced Research Computing Centre, University College London}, addressline={Bidborough House},
postcode={WC1H 9BT}, city={London}, country={United Kingdom}}

\begin{abstract}
We study parametric inference for ergodic diffusion processes with a degenerate diffusion matrix. Existing research focuses on a particular class of hypo-elliptic Stochastic Differential Equations (SDEs), with components split into `rough'/`smooth' and noise from rough components propagating directly onto smooth ones, but some critical model classes arising in applications have yet to be explored. We aim to cover this gap, thus analyse the \emph{highly degenerate} class of SDEs, where components split into further sub-groups. Such models include e.g.~the notable case of generalised Langevin equations. We propose a tailored time-discretisation scheme and provide asymptotic results supporting our scheme in the context of high-frequency, full observations. The proposed discretisation scheme is applicable in much more general data regimes and is shown to overcome biases via simulation studies also in the practical case when only a smooth component is observed. Joint consideration of our study for highly degenerate SDEs and existing research provides a general `recipe' for the development of time-discretisation schemes to be used within statistical methods for general classes of hypo-elliptic SDEs.
\end{abstract}


\end{frontmatter} 

\section{Introduction} \label{sec:intro}
This work addresses the statistical calibration of a wide class of hypo-elliptic diffusions.  
Stochastic Differential Equations (SDEs) are widely used as an effective tool to describe dynamics of the time evolution of phenomena of interest across a multitude of disciplines. Consider SDE models of the following general form: 
\begin{align} \label{eq:sde_1}
   d X_t = V_0 (X_t, \theta) dt + \sum_{j = 1}^d V_j (X_t, \theta) d B_{j,t}, \qquad  X_0 = x \in \mathbb{R}^N,
\end{align}
with $V_j (\cdot, \theta) : \mathbb{R}^N \to \mathbb{R}^N$, $0 \le j \le d$, for parameter $\theta$, driven by the $d$-dimensional standard Brownian motion $B = (B_{1, t}, \ldots, B_{d, t})$, $t \geq 0$, defined upon the filtered probability space $(\Omega, \mathcal{F}, \{\mathcal{F}_t\}_{t \geq 0}, \mathbb{P})$, with $d, N\ge 1$. Several theoretical results about parameter inference for SDEs have been established under positive definiteness conditions on the diffusion matrix $a = V V^\top\in \mathbb{R}^{N\times N}$, with $V = [V_1, \ldots, V_d]$. In such a case, the solution of~(\ref{eq:sde_1}) is referred to as an \emph{elliptic} diffusion. However, many important applications give rise to diffusion processes that 
allow matrix $a$ to be degenerate. We give below examples for such classes of SDEs. 
Under the weak H\"ormander's condition, discussed later in this work, the process defined via the SDE~(\ref{eq:sde_1}) with degenerate diffusion matrix $a$ permits a density with respect to (w.r.t.)~the Lebesgue measure for its transition dynamics, and is referred to as a \emph{hypo-elliptic} diffusion. 
\subsection{Classes of Diffusion Models}  \label{sec:intro_model}
We can summarise the SDE models we consider in this work via two classes of hypo-elliptic diffusions. 
The first hypo-elliptic class is determined via 
the following degenerate SDE: 
\begin{align}
\begin{aligned} \label{eq:hypo-I}  
d X_t & = 
 \begin{bmatrix}
 d X_{S,t} \\[0.2cm] 
 d X_{R,t} 
 \end{bmatrix}
 = 
 \begin{bmatrix}
  V_{S, 0} ( X_t,  \beta_S)  \\[0.2cm] 
  V_{R, 0} ( X_t,  \beta_R) 
 \end{bmatrix} dt
 + 
 \sum_{j = 1}^d 
 \begin{bmatrix}
 \mathbf{0}_{N_S} \\[0.2cm]
 V_{R, j} (X_t, \sigma)
 \end{bmatrix} d B_{j, t},  \qquad 
 X_0 = x = \bigl[ x_{S}^\top,  x_{R}^\top \bigr]^\top \in \mathbb{R}^N. 
\end{aligned} \tag{Hypo-I}  
\end{align} 
Here, the involved SDE functionals are specified as follows: 
\begin{gather*}
 V_{S,0} : \mathbb{R}^N  \times \Theta_{\beta_S} \to \mathbb{R}^{N_S}, \qquad   
V_{R,0} : \mathbb{R}^N  \times \Theta_{\beta_R} \to \mathbb{R}^{N_R}, \qquad 
 V_{R,j} : \mathbb{R}^N  \times \Theta_{\sigma} \to \mathbb{R}^{N_R},  \quad  1 \le j \le d, 
\end{gather*} 
for positive integers $N_S, N_R$ such that $N= N_S + N_R$, and unknown parameter vector
\begin{align*} 
\theta = (\beta_{S}, \beta_{R}, \sigma) \in \Theta 
= \Theta_{\beta_S} \times \Theta_{\beta_R} \times \Theta_{\sigma}
 \subseteq \mathbb{R}^{N_{\beta_S}} \times \mathbb{R}^{N_{\beta_R}} \times\mathbb{R}^{N_{\sigma}},
\end{align*}
for positive integers $N_{\beta_S}, \, N_{\beta_R}, \, N_\sigma$, and  a compact set $\Theta$. For class (\ref{eq:hypo-I}), we will later on introduce a condition upon the vector-valued functionals $\{V_{S,0},V_{R,0},\ldots, V_{R,d}\}$ that is sufficient for the law of $X_t$ to admit a Lebesgue density, and which is related to the weak H\"ormander's condition. In brief, the condition stipulates that $X_{R,t}$ is indeed a rough component, and that \emph{all} coordinates of the drift function $V_{S, 0} (X_t, \beta_S)$ properly relate with the rough component $X_{R ,t}$, so that randomness from $X_{R, t}$ is propagated onto all 
coordinates of vector $X_{S, t}$.

The development of a theoretical and algorithmic framework for parametric inference over class (\ref{eq:hypo-I}) has been the topic of several recent works, see e.g.~\cite{poke:09, dit:19, glot:21, iguchi:22}.
%
%
However, we stress that several important hypo-elliptic SDEs used in practice do not belong in class (\ref{eq:hypo-I}) {-- in the 
sense that \emph{not all} coordinates of the drift function for the smooth components involve the rough component $X_{R,t}$ -- thus} are not covered by recent investigations.   
We specify a key class of practically useful but under-explored hypo-elliptic SDEs via the following equation: 
\begin{align}
\begin{aligned} \label{eq:hypo-II}
d X_t &  
= 
\begin{bmatrix}
 d X_{S_1,t} \\[0.2cm] 
 d X_{S_2,t} \\[0.2cm]  
 d X_{R,t} 
 \end{bmatrix}
 =  
\begin{bmatrix}
    V_{S_1, 0} ( X_{S_1, t}, X_{S_2, t}, \beta_{S_1}) \\[0.2cm]
    V_{S_2, 0} ( X_{t} , \beta_{S_2}) \\[0.2cm] 
    V_{R, 0} ( X_{t}, \beta_{R}) 
\end{bmatrix} dt 
+ \sum_{j = 1}^d 
\begin{bmatrix}
   \mathbf{0}_{N_{S_1}} \\[0.2cm]
   \mathbf{0}_{N_{S_2}} \\[0.2cm] 
   V_{R, j} (X_t, \sigma)
\end{bmatrix} d B_{j, t}; \\[0.2cm] 
X_0 & =  x
=
\bigl[ x_{S_1}^\top, x_{S_2}^\top, x_{R}^\top  \bigr]^\top  
\in \mathbb{R}^N.
\end{aligned}  \tag{Hypo-II}  
\end{align}
Thus, the drift functions of the smooth components are now specified as: 
\begin{align*} 
   V_{S_1, 0} : \mathbb{R}^{N_{S_1}+N_{S_2}} \times \Theta_{{\beta_{S_1}}} 
   \to \mathbb{R}^{N_{S_1}}, \qquad  
   V_{S_2, 0} : \mathbb{R}^{N} \times \Theta_{{\beta_{S_2}}} 
   \to \mathbb{R}^{N_{S_2}},
\end{align*}
and the parameter vector writes as:
\begin{align*}
\theta = (\beta_{S_1}, \beta_{S_2}, \beta_R, \sigma) \in \Theta =
\Theta_{\beta_{S_1}} \times \Theta_{\beta_{S_2}} \times \Theta_{\beta_{R}} \times \Theta_\sigma \subseteq 
\mathbb{R}^{N_{\beta_{S_1}}} \times
\mathbb{R}^{N_{\beta_{S_2}}} \times 
\mathbb{R}^{N_{\beta_R}} \times
\mathbb{R}^{N_{\sigma}},
\end{align*}
for positive integers $N_{S_1}$, $N_{S_2}$, $N_{\beta_{S_1}}$, $N_{\beta_{S_2}}$, and $N=N_{S_1}+N_{S_2}+N_R$, where $\Theta$ is again a compact set. 
{We stress again that the drift function $V_{S_1, 0} (X_{S_1, t}, X_{S_2, t}, \beta_{S_1})$ does not depend on the rough component $X_{R,t}$, 
so noise from $X_{R, t}$ is not directly propagated onto $X_{S_1, t}$. Thus, (\ref{eq:hypo-II}) is treated as a different model class than (\ref{eq:hypo-I}).} We refer to (\ref{eq:hypo-II}) as the \emph{highly degenerate} class of  SDEs.


For {model class} (\ref{eq:hypo-II}) we will {later on} set-up a restriction {over} its constituent vector-valued functionals $\{V_{S_{1,0}},V_{S_{2,0}},V_{R,0},\ldots,V_{R,d}\}$, 
related to  H\"ormander's condition, {which will differ from the corresponding condition assumed for (\ref{eq:hypo-I})}. Roughly, 
such a requirement will guarantee that 
noise from the rough component $X_{R,t}$ indeed propagates onto all coordinates of $X_{S_2, t}$, first, then moving onto $X_{S_1, t}$. Thus,  $X_{S_1, t}$ is `smoother' than $X_{S_2, t}$.
Importantly, a consequence of such a behaviour is that class (\ref{eq:hypo-II}) is \emph{not} included within (\ref{eq:hypo-I}), instead the two classes, (\ref{eq:hypo-I}) and (\ref{eq:hypo-II}), are intrinsically distinct and must be treated separately in terms of theoretical and algorithmic considerations. 
%
\subsection{A Motivating Class of Models} \label{sec:intro_ex}
The \emph{non-Markovian Langevin equation} (or \emph{generalised Langevin equation} (GLE)) is used in a wide range of applications due to its effectiveness in describing complex stochastic systems with memory effects (thus, of non-Markovian structure).
Examples include dynamics observed in protein folding \citep{ay:21}, cancer cells \citep{mi:20}, flocks of birds \citep{fe:20}, molecules \citep{ne:15} and coarse-grained systems \citep{kal:15, li:17}. For simplicity, we consider here a one-dimensional particle with unit mass, and denote its position and momentum, respectively,  by $(q, p)$. Then, a GLE describes the particle dynamics as follows:
\begin{align} \tag{GLE} \label{eq:gle}
\begin{aligned} 
  \dot{q}_t & = \, {p}_t;  \\[0.1cm]
  \dot{p}_t & =  - U'(q_t) 
  - \int_0^t K (t-s) p_s ds 
  + \eta_t, 
 \end{aligned} 
\end{align}
where $U : \mathbb{R} \to \mathbb{R}$ is an appropriate potential function, $K : [0, \infty) \to \mathbb{R}$ is the \emph{memory kernel}, and  $ \eta_t $ is a zero-mean stationary Gaussian noise with auto-correlation specified via a 
fluctuation-dissipation relation in equilibrium, i.e.,   
$\mathbb{E} [ \eta_t  \eta_s ] =  K (t-s)$, $s, t > 0$, given a unit temperature.   
Due to the presence of $K(\cdot)$, particle dynamics will depend on the full state history, with such a property being quite desirable in applications, see e.g.~the references given above. However, the cost of generating the dynamics of model (\ref{eq:gle}) can be overly expensive. Thus, a
standard approach followed in practice is to introduce a parametrisation for the memory kernel $K(\cdot)$ and represent the non-Markovian system (\ref{eq:gle}) as a Markovian one on an extended  space, the latter system being referred to as \emph{Quasi-Markovian Generalised Langevin Equation (QGLE)}. 
Such parametrisation is extremely rich, thus being able to accurately capture the behaviour of systems with general 
true kernel
$K(\cdot)$.  
In particular, a common parametrisation of the memory kernel is the following: 
\begin{align*} 
K (t) = \alpha \delta (t) - \langle e^{- t A} \lambda, \lambda \rangle, 
\qquad  \alpha > 0, \quad \lambda \in \mathbb{R}^m, \quad A \in \mathbb{R}^{m \times m}, \quad m\ge 1,
\end{align*}
with $\delta=\delta(\cdot)$ the Dirac function. 
In this case, the original system in (\ref{eq:gle}) can be equivalently re-written as the following Markovian one: 
\begin{align} 
\tag{QGLE-I} \label{eq:qgle-I} 
\begin{aligned}
\begin{bmatrix}
d q_t  \\[0.1cm] 
d p_t  \\[0.1cm]
d s_t  
\end{bmatrix}
& = 
\begin{bmatrix} 
p_t  \\[0.1cm]
- {U'} (q_t) - \alpha p_t - \langle \lambda, s_t \rangle \\[0.1cm]
- p_t \lambda  - A s_t 
\end{bmatrix} dt 
+ \sum_{j = 1}^{m + 1} 
\begin{bmatrix}
{0} \\[0.1cm]
\sigma_j  
\end{bmatrix} d B_{j, t},  \qquad   s_0 \sim  \mathscr{N} \, ( \mathbf{0}_{m}, \,  I_m), 
\end{aligned} 
\end{align}
with $s_t \in \mathbb{R}^m$ an auxiliary component and $\sigma_j\in \mathbb{R}^{m+1}_+$, $1\le j\le m$. Another typical choice for the memory kernel is the following:
\begin{align*} 
K (t) = \langle e^{- t A} \lambda, \lambda \rangle, 
\end{align*}
in which case the equivalent QGLE writes as: 
\begin{align} 
\tag{QGLE-II} \label{eq:qgle-II}  
\begin{aligned} 
\begin{bmatrix}
d q_t  \\[0.1cm] 
d p_t  \\[0.1cm]
d s_t  
\end{bmatrix}
& = 
\begin{bmatrix}
p_t \\[0.1cm]
-  {U'}(q_t) + \langle \lambda,  s_t \rangle \\[0.1cm] 
- p_t \lambda - A s_t 
\end{bmatrix} dt 
+ 
\sum_{j= 1}^m  
\begin{bmatrix}
 0 \\[0.1cm]
 0 \\[0.1cm] 
 \sigma_j
\end{bmatrix} d B_{j, t}, \qquad  s_0 \sim  \mathscr{N} \, ( \mathbf{0}_{m}, \,  I_m), 
\end{aligned}
\end{align}
with $\sigma_j \in \mathbb{R}_+^m$. Class (\ref{eq:qgle-I}) is investigated, e.g., in \cite{ce:10}. Then, class (\ref{eq:qgle-II}) is popular, e.g., in thermodynamics modelling, see \cite{pav:14, lei:15}. 
Class (\ref{eq:qgle-I}) belongs in  (\ref{eq:hypo-I}), with the rough component comprised of $p_t$, $s_t$.  
For class (\ref{eq:qgle-II}), the rough component consists  only of $s_t$, with~$q_t$ 
depending on the smooth component $p_t$ and not on $s_t$. Thus, (\ref{eq:qgle-II}) lies within class (\ref{eq:hypo-II}). Recently, parametric inference for GLEs within the QGLE setting, under discrete-time observations of the smooth component $q_t$, has been of interest for applications, see e.g.~\cite{fe:20, vr:22}.    
\subsection{Related Works and Objectives}
In this paper we investigate parameter estimation for the two classes of degenerate diffusion processes,  (\ref{eq:hypo-I}) and (\ref{eq:hypo-II}), given discrete-time observations obtained at instances 
$0 \le t_0 < t_1 < \cdots < t_n$, $n \in \mathbb{N}$, with equidistant observation intervals $\Delta_n := t_{i} - t_{i-1}$, $1 \le i \le n$. In particular, we consider the following scenarios for 
the observations: 
\begin{itemize}
\item[1.] \emph{Complete observation regime}, i.e., with all $N$ coordinates of 
$X_{t}$ being observed. 
\item[2.] \emph{Partial observation regime}, i.e., with a strict subset of coordinates being observed. In agreement with applications, in this setting only the upper-most smooth component is assumed to be observed. 
\end{itemize}
Within class  (\ref{eq:hypo-I}), and for the complete observation regime, \cite{dit:19} and \cite{glot:20, glot:21} develop judicious discrete-time (conditionally) Gaussian approximations for the transition distribution. 
Such a proxy provides contrast estimators proven to be asymptotically normal in a high-frequency observation setting, i.e., \limit, with a requirement that the step-size scales as $\Delta_n = o (n^{-1/2})$. {It should be remarked that \cite{dit:19} considered hypo-elliptic SDEs in (\ref{eq:hypo-I}) with particular restrictions, e.g., $N_S = 1$, and proposed contrast estimators separately for $\beta_S \in \Theta_{\beta_S}$ and $(\beta_R, \sigma) \in \Theta_{\beta_R} \times \Theta_{\sigma}$, while \cite{glot:20, glot:21} provided a joint contrast estimator for all parameter $\theta = (\beta_S, \beta_R, \sigma)$ under a more general model structure}. 
In \cite{iguchi:22}, the step-size condition $\Delta_n = o (n^{-1/2})$ to obtain asymptotic normality is weakened to $\Delta_n = o (n^{-1/3})$. In the partial observation regime, with the upper-most smooth component being observed, the missing components must be carefully imputed given the available observations. For SDEs in class (\ref{eq:hypo-I}), it is often the case that the dynamics of the smooth component is determined as 
\begin{align} \label{eq:smooth_FD}
d X_{S, t} = X_{R ,t} dt.
\end{align}
Such a remark also applies for class (\ref{eq:hypo-II}), 
with the role of $X_{S, t}$ taken up by the upper-most smooth component, and the one of $X_{R ,t}$ by the second smooth component.
We keep the discussion within class (\ref{eq:hypo-I}), as this is the context typically looked at in earlier literature.  
For the described setting, it is tempting and, indeed, widely used in practice, to recover the hidden rough component via 
finite-differences, using the observations $\{X_{S, t_{i+1}} \}_i$, i.e.~via 
$ X_{R, t_{i}} = (X_{S, t_{i+1}} - X_{S, t_i}) / \Delta_n$, 
if the step-size $\Delta_n$ is small enough. However, \cite{poke:09, sam:12} show that, in the context of bivariate models within class (\ref{eq:hypo-I}), such an approach delivers (asymptotically) biased estimates of the diffusion parameter $\sigma$. {To side-step the bias, \cite{sam:12} further developed a corrected contrast function built upon a conditionally Gaussian approximation for the rough component $X_{R, t}$.} 
\cite{poke:09, dit:19} argue against applying finite-differences and, instead, consider appropriate It\^o-Taylor schemes leading to non-degenerate conditionally Gaussian approximations for the SDE transition density. Such proxies are then embedded  within MCMC Gibbs samplers or Monte-Carlo 
Expectation-Maximisation (MC-EM) methods to impute the missing components conditionally on observations. {Note that the approach does not require (\ref{eq:smooth_FD}) for the smooth component.
\cite{poke:09} then illustrate empirically that the scheme omitting drift terms of size $\mathcal{O} (\Delta_n^2)$ from the It\^o-Taylor expansion of the smooth component leads to a biased estimation for the drift parameter, $\beta_R$, of the rough component. Subsequent analytical works \citep{dit:19, glot:20, iguchi:22} illustrated that the bias is resolved by adding the drift terms of size $\mathcal{O} (\Delta_n^2)$ in the smooth component within class (\ref{eq:hypo-I})}, as such terms are needed to counterbalance the noise terms of size $\mathcal{O} (\Delta_n^{3/2})$ arising in such an expansion.

{The above discussion suggests that one of the possible recipes for accurate estimation of both hidden components and parameters would be} the development of a conditionally Gaussian approximation for the full coordinates (as such Gaussianity allows for access to computationally effective inference methodologies) obtained via careful inclusion of higher-order terms from the relevant It\^o-Taylor expansion. 
Such an insight for the design of a `correct' discretisation scheme for the purposes of statistical inference has, arguably, not been clearly spelled out in the literature. 
{Furthermore, the benefit of developing such a (locally) Gaussian discretisation is that it produces a closed-form approximate likelihood which is widely applicable to `likelihood-based' inference of degenerate diffusions under both a frequentist/Bayesian framework. We hereby mention that other numerical schemes based on operator splitting (not on It\^o-Taylor expansion) can be utilised in `likelihood-free' inference approaches where only the simulation of SDE paths is required. Recently, exploiting the measure-preserving property of splitting schemes, several works (\cite{buc:20, dit:23}) have presented empirical results showcasing that embedding such schemes in likelihood-free inference can lead to effective parameter estimation for a partially observed hypo-elliptic diffusion belonging in (\ref{eq:hypo-I}). \cite{pil:24} developed an asymptotically unbiased contrast estimator based on splitting schemes for elliptic SDEs having a super-linear drift function and an additive noise, i.e., constant diffusion coefficients.}  
\\ 

Our work aims to provide a comprehensive study of statistical calibration for a wide class of degenerate diffusion models. For this purpose, we review previous works for class  (\ref{eq:hypo-I}), and, then, we establish new analytical results for class (\ref{eq:hypo-II}). 
The main contributions of our work can be  summarised as follows: 
\begin{itemize}
\item[(i)] 
For the highly degenerate class (\ref{eq:hypo-II}), we construct a conditionally Gaussian time-discretisation scheme. The corresponding transition density is well-defined (i.e.~non-degenerate) under a suitable assumption on functionals $\{V_{S_{1,0}},V_{S_{2,0}},V_{R,0},\ldots,V_{R,d}\}$ motivated both by modelling considerations and by adherence to the {weak} H\"ormander's condition. 
We refer to the new proxy as the `locally Gaussian scheme' in agreement with the name assigned by \cite{glot:21} to a conditionally Gaussian scheme developed for class (\ref{eq:hypo-I}).
\item[(ii)] 
For class (\ref{eq:hypo-II}), we define a joint contrast estimator based on the transition density of the locally Gaussian scheme. Then, we show that the estimator is asymptotically normal in the complete observation regime, under a high-frequency observation setting, i.e., for $n \to \infty$, $\Delta_n \to 0$, $n \Delta_n \to \infty$, with the additional condition that the step-size must scale as $ \Delta_n = o (n^{-1/2})$. 
\item[(iii)] Under the partial observation regime often encountered in practical applications, we show via analytical consideration of some case studies that use of a finite-difference method for estimation of hidden components leads to asymptotically biased estimation of the diffusion parameter $\sigma$ for class (\ref{eq:hypo-II}). Thus, we put forward the developed locally Gaussian scheme for (\ref{eq:hypo-II}) as an effective tool to impute hidden components and estimate parameters. 
\item[(iv)] 
By reviewing the methodology already produced in the literature for class (\ref{eq:hypo-I}) and examining the new one produced in this work for (\ref{eq:hypo-II}), we can provide a complete guideline for the development of a discretisation scheme for general degenerate diffusion processes so that the corresponding contrast function does not introduce bias in parameter estimation procedures.
\end{itemize}
\begin{table}
\caption{{Parametric inference for hypo-elliptic SDEs from high-frequency, complete observations}} 
\label{table:complete}
\centering 
\begin{tabular}{c c c c} 
\toprule 
\\[-10pt] 
Work  & Model & Joint estimation & Step-size condition 
\\
\midrule  
\cite{dit:19} 
& \begin{tabular} {c}
   (\ref{eq:hypo-I}) with $N_S = 1$ \\
   and diagonal diffusion matrix
\end{tabular} 
& $\times$ 
& $\Delta_n = o (n^{-1/2})$
\\[0.3cm]
\cite{glot:20, glot:21} 
& (\ref{eq:hypo-I}) 
& $\circ$ 
& $\Delta_n = o (n^{-1/2})$ \\[0.2cm]
\cite{iguchi:22} 
& (\ref{eq:hypo-I}) 
& $\circ$ 
& $\Delta_n = o (n^{-1/3})$ \\[0.2cm] 
This paper 
& (\ref{eq:hypo-II}) 
& $\circ$ 
& $\Delta_n = o (n^{-1/2})$ 
\\ 
\bottomrule   
\end{tabular}
\end{table}
\begin{table}
\caption{{Parametric inference for hypo-elliptic SDEs from partial observations}} 
\label{table:partial}
\centering 
\begin{tabular}{c c c c} 
\toprule 
\\[-10pt] 
Work  
& Model 
& Approach
& Estimation
\\ 
\midrule  
\cite{poke:09}
& 
\begin{tabular}{c} 
(\ref{eq:hypo-I}) satisfying (\ref{eq:smooth_FD}) \\
 with $N_S = N_R = 1$
\end{tabular}
& 
\begin{tabular}{c}
    Finite-differences \\ 
    Gaussian approximation
\end{tabular}
& 
\begin{tabular}{c}
   Bias in $\sigma$ \\ 
   Bias in $\beta_R$ 
\end{tabular}
\\[0.5cm] 
\cite{sam:12} 
& 
\begin{tabular}{c} 
(\ref{eq:hypo-I}) satisfying (\ref{eq:smooth_FD}) \\
 with $N_S = N_R = 1$
\end{tabular}
& \begin{tabular}{c}
    Finite-differences with \\
    a corrected contrast function
\end{tabular}
& Unbiased
\\[0.3cm]
\cite{dit:19} 
& \begin{tabular} {c}
   (\ref{eq:hypo-I}) with $N_S = 1$ \\
   and diagonal diffusion matrix
\end{tabular} 
& Gaussian approximation
& Unbiased
\\[0.3cm]
This paper 
& (\ref{eq:hypo-II}) 
& Gaussian approximation 
& Unbiased 
\\ 
\bottomrule   
\end{tabular}
\end{table} 
A comparison of our work with early works is briefly summarised in Table \ref{table:complete} and \ref{table:partial}.
The rest of the paper is organised as follows. Section \ref{sec:hormander} specifies the class of hypo-elliptic SDEs of relevance for this work, with reference to H\"ormander's condition. Section \ref{sec:scheme}
revisits the correct (in terms of its statistical properties) discretisation scheme for class \mbox{(\ref{eq:hypo-I})} and introduces the one for (\ref{eq:hypo-II}). 
Section \ref{sec:main} provides  our core analytical results of  asymptotic consistency and normality for the statistical estimates obtained via the new scheme, in a complete observation setting. All proofs are collected in an Appendix. We present case studies showcasing the emergence of bias when
standard alternative schemes are called upon or when finite-differences are used to impute unobserved components (a common practice in applications).
Under the correct schemes shown here for classes (\ref{eq:hypo-I}) and (\ref{eq:hypo-II}), we set up a simple Kalman filter for fitting a non-linear sub-class of models commonly arising in applications (we term these \emph{conditional Gaussian non-linear systems}) in the practical partial observation setting.
Section \ref{sec:numerical_studies} presents numerical studies, for the partial observation regime, both for simple models and ones relevant to real applications,   within class (\ref{eq:hypo-II}). The code used in the numerical studies is available at {https://github.com/YugaIgu/calibration-hypoSDEs}.
We finish with some conclusions in Section \ref{sec:end}. \\[0.2cm]

\noindent \textbf{Notation.} 
For the highly degenerate class (\ref{eq:hypo-II}), to establish a common notation with (\ref{eq:hypo-I}), we use the argument $x_S=(x_{S_1},x_{S_2})$ and also set: 
\begin{gather*}
X_{S, t}
= \bigl[ X_{S_1, t}^\top, X_{S_2, t}^\top \bigr]^\top\in \mathbb{R}^{N_S}, \quad  
N_S = N_{S_1} + N_{S_2};\\[0.2cm]
\beta_S = \bigl[\beta_{S_1}^{\top}, \beta_{S_2}^{\top}\bigr]^{\top}\in 
\Theta_{\beta_S}, \quad \Theta_{\beta_S} = \Theta_{\beta_{S_1}}\times 
\Theta_{\beta_{S_2}}, \quad N_{\beta_{S}} = N_{\beta_{S_1}} + N_{\beta_{S_2}};
\\[0.2cm]
V_{S,0} (x, \beta_S) 
= 
\bigl[ V_{S_1, 0} (x_S, \beta_{S_1})^\top, \; 
V_{S_2, 0} (x, \beta_{S_2})^\top  \bigr]^\top.
\end{gather*}
For $x  \in \mathbb{R}^N$ and $\theta = (\beta_S, \beta_R, \sigma) \in \Theta$, 
we write
\begin{align}
V_0 (x, \theta)  & = 
\left[ V_{S, 0} (x, \beta_S)^\top, \; 
V_{R,0}  (x, \beta_R)^\top \right]^\top;  \label{eq:V0} \\[0.1cm] 
V_j (x, \theta)  & = 
\left[ \mathbf{0}_{N_S}^\top, \;  V_{R,j} (x, \sigma)^\top \right]^\top, \qquad 1 \le j \le d.  \nonumber
\end{align}  
For $\varphi (\cdot , \theta) : \mathbb{R}^N \to \mathbb{R}$, $\theta \in \Theta$, bounded up to 2nd order derivatives, we define the differential operators $\mathcal{L}$ and $\mathcal{L}_j, \; 1 \le  j  \le d$, as: 
\begin{gather*}
\mathcal{L} \varphi (x ,  \theta) 
= \sum_{i=1}^N V_0^i (x, \theta) \frac{\partial \varphi} {\partial x_i}(x, \theta)  
+ \frac{1}{2}  \sum_{i_1, i_2 = 1}^N \sum_{k=1}^d  V_k^{i_1} (x, \theta) V_k^{i_2} (x, \theta)  \frac{\partial^2 \varphi }{\partial x_{i_1} \partial x_{i_2} } (x,\theta);  \\
\mathcal{L}_j \varphi (x , \theta)
= \sum_{i=1}^N V_j^i (x, \theta) \frac{\partial \varphi}{\partial x_i}(x , \theta), \ \ 1 \le j \le d,  
\end{gather*} 
for $(x, \theta) \in \mathbb{R}^N \times \Theta$. 
Application of the above differential operators is extended to vector-valued functions in the apparent way, via separate consideration of each scalar component. 
 We denote the probability law of the process $\{ X_{t} \}_{t \geq 0}$ under a parameter $\theta \in \Theta$ as $\mathbb{P}_{\theta}$, and we write
$\probconv, \ \ \distconv$
for convergence in probability and distribution, respectively, under the true parameter $\trueparam$.  We write the expectation under the probability law $\mathbb{P}_\theta$ as $\mathbb{E}_\theta$ to emphasise the dependence on $\theta \in \Theta$. For $u \in \mathbb{R}^n$, $n \in \mathbb{N}$ and the multi-index $\alpha \in \{1, \ldots, n \}^l$, $l \in \mathbb{N}$, we define $\textstyle{\partial^u_\alpha = \partial^l/\partial u_{\alpha_1} \cdots \partial u_{\alpha_l}}$, i.e.~an operator acting on maps $\mathbb{R}^n \to \mathbb{R}$, and then extended, by separate application on each co-ordinate, on maps $\mathbb{R}^n \to \mathbb{R}^{m}$, $m \in \mathbb{N}$. {We denote by 
$C_p^\infty (\mathbb{R}^{n_1} \times \Theta; \mathbb{R}^{n_2})$, $n_1, n_2 \in \mathbb{N}$, 
the space of functions $f: \mathbb{R}^{n_1} \times \Theta \to \mathbb{R}^{n_2}$ such that $f (x, \theta)$ is infinitely differentiable w.r.t.~$x \in \mathbb{R}^{n_1}$ for all $\theta \in \Theta$, and for any $\alpha \in \{1, \ldots, n_1\}^\ell, \, \ell \ge 0$, $\partial_x^{\alpha} f (x, \theta)$ is of polynomial growth in $x \in \mathbb{R}^{n_1}$ uniformly in $\theta \in \Theta$. For $x, y \in \mathbb{R}^N$, we write $\textstyle \langle x, y \rangle = \sum_{i = 1}^N x_i y_i$ and $\|x\| = \sqrt{\langle x, x \rangle}$.}  

\section{Hypo-Elliptic SDEs} \label{sec:hormander}
%
We fully specify the classes of SDEs of interest in (\ref{eq:hypo-I}) and (\ref{eq:hypo-II}), by providing, in each case, appropriate conditions on the collection of functionals
$\{V_{0}, V_{1},\ldots, V_{d}\}$, motivated by modelling considerations and the existence of a Lebesgue density for the SDE transition dynamics. 
We illustrate later on that the imposed conditions suffice so that the locally Gaussian scheme for $X_t$ in (\ref{eq:hypo-II}) we put forward in this paper is non-degenerate. 


%
\subsection{H\"ormander's Condition} \label{sec:hormander_review} 
We quickly review the definition of H\"ormander's condition.
Consider the class of SDEs with the general form in~(\ref{eq:sde_1}). We define
\begin{align*}
{\tilde{V}_0 (x, \theta) = V_0 (x, \theta) 
- \tfrac{1}{2} \sum_{k = 1}^d \mathcal{L}_k V_k (x, \theta)}, \quad (x, \theta) \in \mathbb{R}^N  \times \Theta.
\end{align*}
From standard properties of It\^o's processes,  $\tilde{V}_0$ is the drift function of (\ref{eq:sde_1}) when written as a Stratonovich-type SDE. The functionals $\{\tilde{V}_0, \ldots, V_d\}$ of the Stratonovich SDE can be corresponded to differential operators, the latter applying on mappings on $\mathbb{R}^{N}\to \mathbb{R}^{N}$ and giving as outcome mappings, again, on the same spaces. In particular, we have: 
\begin{align*}
 \tilde{V}_0 
\mapsto  \sum_{i = 1}^N \tilde{V}_0^i (x) {\partial_{x_i}}, 
 \qquad  
 V_k \mapsto \sum_{i = 1}^N {V}_k^i (x) {\partial_{x_i}}, 
 \quad  1 \le k \le d.  
\end{align*}
Without confusion, we use the same notation both for the SDE 
functionals and the corresponding differential operators.
Parameter $\theta$ is removed from the expressions for simplicity. For two functionals (equivalently, differential operators) as above, 
$\textstyle{ W = \sum_{i = 1}^N W^i (x)\partial_{x_i}}$ and $\textstyle{ Z = \sum_{i = 1}^N Z^i (x)\partial_{x_i}}$, the Lie bracket is defined as
\begin{align*} 
[W, Z] = W\,Z - Z\,W, 
\end{align*} 
that is, for a given $x \in \mathbb{R}^N$, 
\begin{align*} 
[W, Z] (x) = \sum_{i = 1}^N \bigl\{ W^i (x) \partial_{x_i} Z (x)
- Z^i (x) \partial_{x_i} W (x) \bigr\}\in \mathbb{R}^{N}
. 
\end{align*} 
We introduce the collections of functionals
\begin{gather*} 
\mathscr{H}_0 = \big\{V_1, \ldots, V_d \big\},  
\qquad 
\mathscr{H}_k = \Big\{ \big\{\,[\tilde{V}_0, V], [V_r, V]\,\big\}  :  V  \in \mathscr{H}_{k-1}, \, 1  \le r \le d \Big\}, 
\quad 
k \ge 1.
\end{gather*}
Then, H\"ormander's condition is stated as follows: 
\begin{definition}
H\"ormander's condition is said to hold at a point $x \in \mathbb{R}^N$ if there exists $M \ge 0$ such that  
\begin{align} \label{eq:H}
\mrm{span} \big\{ V (x) : V \in \mathscr{H}_j, \ \  0 \le j \le M \big\} = \mathbb{R}^N.  
\end{align}
%
\end{definition}
\noindent H\"ormander's condition implies that for any $t > 0$ and an initial condition $X_0=x\in\mathbb{R}^{N}$, the law of $X_t$ is absolutely continuous w.r.t.~the Lebesgue measure.  Also, if the coefficients of the SDE are infinitely-times differentiable, with partial derivatives of all orders being bounded, then the Lebesgue density is smooth, see, e.g., \cite{nua:06, pav:14}. {When H\"ormander's condition holds with $M=0$, the condition is typically referred to as \emph{strong} H\"ormander's condition and is equivalent to the SDE being an elliptic one. It is easy to see that in our setting of a degenerate diffusion matrix, H\"ormander's condition cannot hold with $M = 0$.  If (\ref{eq:H}) holds with $M \ge 1$, then the term \emph{weak} H\"ormander's condition is often used, see e.g. \cite{piga:18}. Furthermore, if (\ref{eq:H}) holds uniformly in the initial point $x \in \mathbb{R}^N$, we have the \emph{uniform} H\"ormander's condition, which is invoked frequently in the literature, see e.g. \cite{bally:96, cass:09}. 
Precisely, the uniform H\"ormander's condition is stated as follows.  
There exists an integer $M \ge 0 $ such that 
\begin{align} \label{eq:inf_H}
\inf_{x \in \mathbb{R}^N}
\ 
\inf_{\substack{ \xi \in \mathbb{R}^N \mrm{s.t.} \\ \| \xi \| = 1}}
\, 
\sum_{j = 0}^M
\sum_{V \in \mathscr{H}_j}
\Bigl\langle 
V (x), \, \xi 
\Bigr\rangle^2 > 0. 
\end{align} 
To establish a statistical theory for models (\ref{eq:hypo-I}) and (\ref{eq:hypo-II}), we work below with assumptions associated with the uniform H\"ormander's condition.}
\subsection{Diffusion Classes (\ref{eq:hypo-I}) and (\ref{eq:hypo-II})}
We  now set up separate conditions for the  
SDEs in classes (\ref{eq:hypo-I}) and (\ref{eq:hypo-II}). These will make use of the drift function, 
$\tilde{V}_0=\tilde{V}_0(x,\theta)$, of the Stratonovich version of the SDEs.  
%
%
%
%
%
%
%
%
For $1 \le j \leq k \le N$, we define the projection operator $\mathrm{proj}_{j,k} : \mathbb{R}^N \to \mathbb{R}^{k - j + 1}$ as 
\begin{align*}
x = \bigl[x_1,\ldots, x_N\bigr]^{\top} \mapsto \mathrm{proj}_{j, k} (x)
= \bigl[ x_{j}, \ldots, x_k \bigr]^\top. 
\end{align*}
We also introduce: 
\begin{gather*}
\widetilde{\mathscr{H}}_0 = \mathscr{H}_0, 
\qquad 
\widetilde{\mathscr{H}}_1 = \Bigl\{ [\widetilde{V}_0, V] : \, V \in \widetilde{\mathscr{H}}_0 \,  \Bigr\} 
\qquad 
\widetilde{\mathscr{H}}_2 = \Bigl\{ 
\bigl[\widetilde{V}_0, [\widetilde{V}_0, V]\bigr] : \, V \in \widetilde{\mathscr{H}}_0 \, \Bigr\}.  
\end{gather*}
For (\ref{eq:hypo-I}) and (\ref{eq:hypo-II}), we assign the following conditions to fully specify the structure of the corresponding degenerate system of SDEs. 
{
\begin{condition}[Classes of SDEs]
\label{assump:hypo}
\begin{enumerate}
\item[I.] For class (\ref{eq:hypo-I}), it holds that: 
\begin{align} 
& \inf_{(x, \sigma) \in \mathbb{R}^N \times \Theta_\sigma}
\ 
\inf_{\substack{ \xi \in \mathbb{R}^{N_R} \\  \mrm{s.t.} \, 
\| \xi \| = 1}} 
\sum_{V \in \widetilde{\mathscr{H}}_0} 
\Bigl\langle 
\mathrm{proj}_{N_{S}+1, N} 
\bigl\{ V (x, \theta) \bigr\}, \ \xi 
\Bigr\rangle^2
= 
\inf_{(x, \sigma) \in \mathbb{R}^N \times \Theta_\sigma}
\ 
\inf_{\substack{ \xi \in \mathbb{R}^{N_R} \\  \mrm{s.t.} \, 
\| \xi \| = 1}} 
\sum_{k = 1}^d 
\Bigl\langle V_{R, k} (x, \sigma), \ \xi 
\Bigr\rangle^2 > 0; 
\nonumber
\\[0.2cm]  
& 
\inf_{(x, \theta) \in \mathbb{R}^N \times \Theta}
\ 
\inf_{\substack{ \xi \in \mathbb{R}^{N} \\ \mrm{s.t.} \, 
\| \xi \| = 1}}
\sum_{j = 0, 1} 
\sum_{V \in \widetilde{\mathscr{H}}_j} 
\Bigl\langle 
V (x, \theta), \ \xi 
\Bigr\rangle^2 > 0.
\label{eq:span-I} 
\end{align} 
\item[II.] In the case of class (\ref{eq:hypo-II}), it holds that: 
\begin{align}
&
\inf_{(x, \sigma) \in \mathbb{R}^N \times \Theta_\sigma}
\ 
\inf_{\substack{ \xi \in \mathbb{R}^{N_R} \\  \mrm{s.t.} \, 
\| \xi \| = 1}} 
\sum_{V \in \widetilde{\mathscr{H}}_0} 
\Bigl\langle 
\mathrm{proj}_{N_{S}+1, N} \bigl\{ V (x, \theta) \bigr\}, \ \xi 
\Bigr\rangle^2
= 
\inf_{(x, \sigma) \in \mathbb{R}^N \times \Theta_\sigma}
\ 
\inf_{\substack{ \xi \in \mathbb{R}^{N_R} \\  \mrm{s.t.} \, 
\| \xi \| = 1}} 
\sum_{k = 1}^d 
\Bigl\langle V_{R, k} (x, \sigma), \ \xi 
\Bigr\rangle^2 > 0; 
\nonumber
\\[0.2cm]  
& \inf_{(x, \theta) \in \mathbb{R}^N \times \Theta}
\  \inf_{\substack{ \xi \in \mathbb{R}^{N_{S_2} + N_R} \\ 
\mrm{s.t.} \,  \| \xi \| = 1}} 
\sum_{j = 0,1}
\sum_{V \in \widetilde{\mathscr{H}}_j}
\Bigl\langle 
\mathrm{proj}_{N_{S_{1}}+1,N} \bigl\{ V (x, \theta)
\bigr\}, 
\ \xi 
\Bigr\rangle^2  > 0; 
\nonumber
\\[0.2cm]
& \inf_{(x, \theta) \in \mathbb{R}^N \times \Theta}
\  \inf_{\substack{ \xi \in \mathbb{R}^{N} \\  \mrm{s.t.} \, \| \xi \| = 1}}
\sum_{j = 0,1,2}
\sum_{V \in \widetilde{\mathscr{H}}_j}
\Bigl\langle 
V (x, \theta), \ \xi 
\Bigr\rangle^2  
> 0. \nonumber
\end{align} 
\end{enumerate}
\end{condition}
}
\noindent Note that the uniform  H\"ormander's condition holds for (\ref{eq:hypo-I}) and (\ref{eq:hypo-II}) under \ref{assump:hypo}-I and \ref{assump:hypo}-II (for $M=1$ and $M=2$) respectively.
%
%
\begin{remark} 
Conditions \ref{assump:hypo}-I and \ref{assump:hypo}-II 
 separate classes (\ref{eq:hypo-I}) and (\ref{eq:hypo-II}).
In particular, the top equation in both \ref{assump:hypo}-I, \ref{assump:hypo}-II implies that the diffusion matrix of the rough component is of full rank, thus $X_{R,t}$ acquires the roughness of an elliptic SDE.
The second equation in \ref{assump:hypo}-I, \ref{assump:hypo}-II ensures that all coordinates of component $X_{S,t}$ and $X_{S_2,t}$, respectively, possess the same smoothness as integrals of elliptic SDEs, {that is, contain a Gaussian noise with the size of $\mathcal{O} (\Delta_n^{3/2})$}. 
Finally, the third equation in \ref{assump:hypo}-II implies that component  $X_{S_1,t}$ has the smoothness of second integrals of elliptic SDEs, {i.e., contain a Gaussian noise of size $\mathcal{O} (\Delta_n^{5/2})$}. Note, e.g., that (\ref{eq:span-I}) will not hold for (\ref{eq:hypo-II}), since for the highly degenerate case, due to the drift function of the upper-most component $X_{S_1, t}$ not involving $X_{R, t}$, 
we have 
\begin{align*} 
\mathrm{proj}_{1, N_{S_1}} [\tilde{V}_{0}, V_{k}]= \mathbf{0}_{N_{S_1}},\quad  1 \le k \le d.
\end{align*}
To check this latter equation, notice that for both (\ref{eq:hypo-I}) and (\ref{eq:hypo-II}) we have: i) $V_0$ and $\tilde{V}_{0}$ coincide on the smooth coordinates; ii) $\tilde{V}_0 V_k$ is zero on the smooth coordinates. Then, for (\ref{eq:hypo-II}) we additionally have that $V_k \tilde{V}_0$ is zero on the upper-most $N_{S_1}$ coordinates due to the particular choice of $V_0=V_0(x_{S_1},x_{S_2})$.   
\end{remark} 
\begin{remark}
We introduced conditions \ref{assump:hypo}-I \& II upon functionals $\{V_{0}, V_{1},\ldots, V_{d}\}$, so that 
the two classes of SDEs, (\ref{eq:hypo-I}) and (\ref{eq:hypo-II}), possess sufficient structure to allow for their intended use for the modelling objectives in mind.
It turns out that the exact same conditions \ref{assump:hypo}-I \& II play a key role so that the locally Gaussian approximation for $X_{t_{i+1}} | X_{t_i}= x$ written down later in Section \ref{sec:scheme} are 
well-defined with a positive definite covariance matrix for all $(x, \theta) \in \mathbb{R}^N \times \Theta$.  
\end{remark}
{
\subsubsection{An Example for \ref{assump:hypo}-I} \label{sec:ex_hypo1}
We provide an example for \ref{assump:hypo}-I via the following bivariate underdamped Langevin equation: 
\begin{align} 
\begin{aligned} \label{eq:ULE} 
 d q_t & = p_t d t; \\[0.2cm] 
 d p_t & = \bigl( -  U' (q_t) -  \gamma  p_t \bigr) dt
 + \sigma d B_{1, t}, 
\end{aligned}
\end{align} 
where $U : \mathbb{R} \to \mathbb{R}$, $\gamma > 0$ and $\sigma > 0$. We write $x = (q, p) \in \mathbb{R}^2$ and $\theta = (\gamma, \sigma) \in \Theta$, where $\Theta$ is an appropriate compact space. We have: 
\begin{align*}
\widetilde{V}_0 = p \partial_q + (- U'(q) - \gamma p) \partial_p, \qquad 
V_1 = \sigma \partial_p, 
\qquad 
[\widetilde{V}_0, V_1] 
= - \sigma \partial_q + \gamma \sigma \partial_p, 
\end{align*}
and then 
\begin{align*}
V_{1} (x, \sigma) = 
\begin{bmatrix}
0 \\
V_{R, 1} (x, \sigma)  
\end{bmatrix}
= 
\begin{bmatrix}
0 \\
\sigma  
\end{bmatrix}, 
\qquad 
[\widetilde{V}_0, V_{1}] (x, \theta) = 
\begin{bmatrix}
- \sigma  \\
\gamma \sigma  
\end{bmatrix}. 
\end{align*} 
Thus, SDE (\ref{eq:ULE}) satisfies condition \ref{assump:hypo}-I and lies within the framework of class (\ref{eq:hypo-I}). 
}
\subsubsection{An Example for \ref{assump:hypo}-II} \label{sec:ex_hypo2}
We provide an example for \ref{assump:hypo}-II via the following three-dimensional hypo-elliptic diffusion motivated from model class (\ref{eq:qgle-II}):
\begin{align} 
\begin{aligned} \label{eq:qgle} 
 d q_t & = p_t d t; \\[0.2cm] 
 d p_t & = \bigl( - U'(q_t) + \lambda  s_t \bigr) dt;  \\[0.1cm] 
 d s_t & = \bigl( - \lambda p_t  - \alpha s_t \bigr) dt +  \sigma d B_{1, t}, 
\end{aligned}
\end{align} 
where $U:\mathbb{R} \to \mathbb{R}$, $\alpha  > 0$, $\sigma > 0$ and $\lambda \in \mathbb{R} \setminus \{0\}$. 
{Notice that the drift function of the positional component $q_t$ is independent of the rough component $s_t$.} 
In this case, for $ x=(p,q,s) \in \mathbb{R}^3$, $\theta = (\lambda, \alpha, \sigma) \in \Theta$ with some appropriate compact space $\Theta$, we have:
\begin{gather*}
\tilde{V}_{0} = p\,\partial_{q}
+ \bigl( - U'(q) + \lambda s \bigr) \partial_{p} 
+ \bigl( -\lambda p - \alpha s \bigr) \partial_{s}, 
\qquad  V_{1} = \sigma \partial_{s}, \\[0.2cm]
[\tilde{V}_{0}, V_{1}]  = 
- \lambda \sigma \partial_p + \alpha \sigma \partial_{s}, \qquad  
\bigl[\tilde{V}_{0}, [\tilde{V}_{0}, V_{1}] \bigr] = \lambda \sigma \partial_{q}  
- \lambda \alpha \sigma \partial_{p} 
+ \sigma (- \lambda^2  + \alpha^2 ) \partial_{s}. 
\end{gather*}
{
We obtain: 
\begin{align*}
 V_1 (x, \sigma) = 
 \begin{bmatrix}
0 \\ 
0 \\ 
V_R (x, \sigma)
\end{bmatrix} 
= \begin{bmatrix}
0 \\ 
0 \\ 
\sigma
\end{bmatrix}, 
\qquad 
[\widetilde{V}_0, V_1] (x, \theta) 
= \begin{bmatrix}
0 \\ 
-\lambda \sigma  \\ 
\alpha \sigma
\end{bmatrix},
\qquad 
\bigl[
 \widetilde{V}_0, [\widetilde{V}_0, V_1]
\bigr] (x, \theta) 
= \begin{bmatrix}
\lambda \sigma  \\ 
-\lambda \alpha \sigma \\ 
\sigma (- \lambda^2 + \alpha^2) 
\end{bmatrix}.  
\end{align*}
} 
%
Thus, SDE (\ref{eq:qgle}) satisfies condition \ref{assump:hypo}-II and lies within the framework of class (\ref{eq:hypo-II}). 
\section{Time-Discretisation of Hypo-Elliptic SDEs} \label{sec:scheme}
We discuss time-discretisation schemes for hypo-elliptic diffusions  within classes (\ref{eq:hypo-I}) and (\ref{eq:hypo-II}), with the focus being on the performance of the schemes for the purposes of parametric inference. We set up the context by reviewing  schemes proposed in literature for class (\ref{eq:hypo-I}). 
We then propose a new scheme for the highly degenerate diffusion class (\ref{eq:hypo-II}) that will later be proven to possess desirable statistical properties. 
Hereafter, to distinguish among the two classes of SDEs, we use the notation $\textstyle{ X_t^{(\mrm{I})}}$ and $\textstyle{X_t^{(\mrm{II})}}$ for  processes in (\ref{eq:hypo-I}) and (\ref{eq:hypo-II}), respectively.     
\subsection{Time-Discretisation of (\ref{eq:hypo-I}) 
-- Brief Review}
We review relevant schemes for the hypo-elliptic class (\ref{eq:hypo-I}) used in the literature. 
First, the classical \emph{Euler-Maruyama} scheme is defined as follows, for $0 \le i \le n$, 
\begin{align*} 
\begin{aligned}
 X_{S, i+1}^{\mrm{EM}, (\mrm{I})}   
 & = X_{S, i}^{\mrm{EM}, (\mrm{I})} 
+ V_{S, 0} \bigl( X_{i}^{\mrm{EM}, (\mrm{I})}, \beta_S \bigr) \Delta_n ; \\ 
 X_{R, i+1}^{\mrm{EM}, (\mrm{I})} 
 & = X_{R, i}^{\mrm{EM}, (\mrm{I})} 
 + V_{R, 0} \bigl( X_{i}^{\mrm{EM}, (\mrm{I})}, \beta_R \bigr) \Delta_n 
 + \sum_{j = 1}^d V_{R, j} \bigl( X_{i}^{\mrm{EM}, (\mrm{I})}, \sigma \bigr) 
  \times \bigl( B_{j, i+1} - B_{j, i} \bigr), 
\end{aligned} 
\end{align*}
with subscript $i$ in $X_{i}^{\mrm{EM}, (\mrm{I})}$ and $B_{j, i}$ indicating the time instance $t_i=i\Delta_n$. The approximation of the smooth component does not involve noise, thus the Euler-Maruyama scheme is degenerate. \cite{poke:09} studied some example bivariate SDEs with drift function $V_{S, 0}(x, \beta_S \bigr) = x_{R}$, and showed that use of the Euler-Maruyama scheme in the high-frequency partial observation regime, where only the smooth component $X_{S,t}$ is observed, induces bias in parameter estimates. Note that in this setting the unobserved component $X_{R,t}$ is estimated via a 
finite-difference approach, i.e., 
\begin{align*} 
   X_{R, i}^{\, (\mrm{I})} \approx 
    \hat{X}_{R, i}^{\, (\mrm{I})} 
    = \frac{X_{S, i+1}^{\, (\mrm{I})} - X_{S, i}^{\, (\mrm{I})}}{\Delta_n},
\end{align*}
and such an imputation is a main cause for the presence of bias in the estimation of  parameter $\sigma$. To sidestep the above issue, \cite{poke:09} proposed the following conditionally Gaussian scheme, for $0 \le i \le n$: 
\begin{align} \label{eq:IS}
\begin{aligned} 
 \widetilde{X}_{S, i+1}^{\, (\mrm{I})} & = 
 \widetilde{X}_{S, i}^{\, (\mrm{I})} 
 + V_{S, 0} \bigl( \widetilde{X}_{i}^{\, (\mrm{I})} ,  \beta_S \bigr) \Delta_n  
+ \sum_{j = 1}^d \mathcal{L}_j V_{S, 0} (\widetilde{X}_{i}^{\, (\mrm{I})}, \theta) \int_{t_i}^{t_{i+1}} \int_{t_i}^{u} d B_{j,v} du;   \\ 
 \widetilde{X}_{R, i+1}^{\, (\mrm{I})} 
 & = \widetilde{X}_{R, i}^{\, (\mrm{I})} 
 + V_{R, 0} \bigl( \widetilde{X}_{i}^{\, (\mrm{I})} , \beta_R \bigr) \Delta_n 
 + \sum_{j = 1}^d V_{R, j} \bigl( \widetilde{X}_{i}^{\, (\mrm{I})} , \sigma \bigr) \times \bigl( B_{j, i+1} - B_{j, i} \bigr). 
\end{aligned} 
\end{align}
%
Note that now the smooth component $\widetilde{X}_{S,i+1}^{\, (\mrm{I})}$ involves Gaussian noise after application of an It\^o-Taylor expansion for $\textstyle V_{S, 0} (X_t, \beta_S)$. Under condition \ref{assump:hypo}-I, $(\widetilde{X}^{\, (\mrm{I})}_{S, i+1}, \widetilde{X}_{R, i+1}^{\, (\mrm{I})})$ is conditionally Gaussian with an invertible covariance matrix. Then, \cite{poke:09} utilised the well-posed likelihood of scheme (\ref{eq:IS}) to estimate  both the hidden paths of the rough components and the parameters via a Bayesian approach, namely Gibbs sampling. Under 
a high-frequency observation setting, they empirically showed that the estimate of parameter $\sigma$ is asymptotically unbiased, 
but the estimator of the drift parameter $\beta_R$ based on scheme (\ref{eq:IS}) suffers from bias even in the complete observation regime. 

\cite{glot:20} introduced the `local Gaussian' scheme, 
where, for $0 \le i \le n$, 
%
\begin{align}  \label{eq:lg_hypo_I}
\begin{aligned} 
\bar{X}_{S, i+1}^{\, (\mrm{I})} & = \bar{X}_{S, i}^{\, (\mrm{I})}  +
V_{S, 0} \bigl( \bar{X}_{i}^{\, (\mrm{I})},  \beta_S \bigr) \Delta_n + \tfrac{\Delta_n^2}{2} \mathcal{L} V_{S, 0} 
\bigl(\bar{X}_{i}^{\, (\mrm{I})}, \theta \bigr) 
+ \sum_{j = 1}^d \mathcal{L}_j V_{S, 0} (\bar{X}_{i}^{\, (\mrm{I})}, \theta) \int_{t_i}^{t_{i+1}}  \int_{t_i}^{u}  d B_{j,v} du;  \\ 
\bar{X}_{R, i+1}^{\, (\mrm{I})} & = \bar{X}_{R, i}^{\, (\mrm{I})} + V_{R, 0} \bigl( \bar{X}_{i}^{\, (\mrm{I})}, \beta_R \bigr) \Delta_n + \sum_{j = 1}^d V_{R, j} 
\bigl( \bar{X}_{i}^{\, (\mrm{I})}, \sigma \bigr) \times \bigl( B_{j, i+1} - B_{j, i} \bigr).  
\end{aligned} \tag{LG-I}
\end{align} 
Compared to (\ref{eq:IS}), scheme (\ref{eq:lg_hypo_I}) includes term $\textstyle{ \Delta_n^2 (\mathcal{L} V_{S, 0}) (\bar{X}^{\, (\mrm{I})}_i, \theta)}/2$ in the smooth component.
\cite{glot:20} illustrate the significance of this term for the purposes of parameter inference, 
by proving asymptotic consistency and normality for the contrast estimator derived from the likelihood of the discretisation scheme (\ref{eq:lg_hypo_I}), in the high-frequency, complete observation regime, namely, $n \to \infty$, $\Delta_n \to 0$, $n \Delta_n \to \infty$, under the step-size condition $\Delta_n = o (n^{-1/2})$. 
%
%
\begin{remark}
\cite{dit:19} applied a strong 1.5 order scheme (\citep{kloe:92}) to construct contrast estimators for class (\ref{eq:hypo-I}) with $N_S = 1$, under strong conditions for the diffusion matrix so that the scheme becomes conditionally Gaussian. Then, they provided two separate contrast functions for estimating $\beta_S$ and $(\beta_R, \sigma)$ from the  approximate Gaussian density for $X_S$ and $X_R$, respectively, rather than the joint density. As noted in Remark 4.6 in \cite{glot:20}, the separate contrast functions result in a larger asymptotic variance for the estimation for  $\beta_S$ compared with the single contrast estimator defined via the joint density of rough and smooth components. 
\end{remark}
\subsection{Time-Discretisation of (\ref{eq:hypo-II})}  \label{sec:lg-hypo-II} 
We propose a time-discretisation scheme for the second hypo-elliptic class (\ref{eq:hypo-II}),
with desirable properties for the purposes of parameter inference. 
The brief review of schemes for (\ref{eq:hypo-I}) in the previous section suggests that the discretisation scheme for (\ref{eq:hypo-II}) should satisfy the following two key criteria: 
\begin{itemize}
    \item[I.] 
    The scheme should be conditionally non-degenerate, i.e., the law of $X_{t_{i+1}}$ given $X_{t_i}$ should admit a Lebesgue transition density for the full coordinates. This will allow to impute unobserved paths conditionally on observations without making use of bias-inducing finite-difference approximations. 
    \item[II.]  The scheme should involve deterministic terms 
    obtained from careful truncation of the stochastic Taylor expansion for the drift of the smooth component,  $V_{S,0} (X_{t}^{\, (\mrm{II})}, \beta_{S})$, so that the contrast estimator corresponding to the scheme is asymptotically unbiased under the high-frequency, complete observation regime.   
\end{itemize}
As for Criterion I, we will explain later in Section \ref{sec:case_2} that, indeed, use of a degenerate discretisation scheme or of 
finite-differences to estimate hidden components induces a bias in the estimation of parameters. Based upon the above key criteria, we propose the following discretisation scheme for (\ref{eq:hypo-II}): 
%
%
\begin{align}
\begin{aligned} \label{eq:lg_hypo_II} 
\bar{X}_{S_1, i+1}^{\, \mathrm{(II)}} 
    & = 
    \mu_{S_1} \bigl( \Delta_n, \bar{X}_{i}^{\, \mathrm{(II)}}, \theta \bigr)  
    + \sum_{j = 1}^d \mathcal{L}_{j} \mathcal{L} V_{S_1, 0} 
    \bigl(\bar{X}_{i}^{\, \mathrm{(II)}} , \theta \bigr) \int_{t_i}^{t_{i+1}} \int_{t_i}^u \int_{t_i}^v  
    dB_{j, w} dv du;  \\[0.2cm] 
    \bar{X}_{S_2, i+1}^{\, \mathrm{(II)}} 
    & = 
    \mu_{S_2} \bigl( \Delta_n,  \bar{X}_{i}^{\, \mathrm{(II)}} , \theta \bigr) 
    + \sum_{j = 1}^d \mathcal{L}_{j} V_{S_2, 0} (\bar{X}_{i}^{\, \mathrm{(II)}} , \theta) \int_{t_i}^{t_{i+1}}  \int_{t_i}^u d B_{j, v} du;  \\[0.2cm]  
    \bar{X}_{R, i+1}^{\, \mathrm{(II)}} 
    & = 
    \mu_{R} \bigl( \Delta_n, \bar{X}_{i}^{\, \mathrm{(II)}} , \theta \bigr)  
    + \sum_{j = 1}^d V_j \bigl( \bar{X}_{i}^{\, \mathrm{(II)}} , \sigma \bigr) \times \bigl(B_{j,i+1} - B_{j, i} \bigr),   
\end{aligned} \tag{LG-II}
\end{align} 
where we have set, for $(\Delta, x, \theta) \in (0, \infty) \times \mathbb{R}^N \times \Theta$, 
\begin{align*}  
\begin{bmatrix}
    {\mu}_{S_1}  (\Delta, x, \theta)  \\[0.3cm] 
    {\mu}_{S_2}  (\Delta, x, \theta)  \\[0.3cm]  
    {\mu}_{R}  (\Delta, x, \theta)  
\end{bmatrix}
=
\begin{bmatrix}
    x_{S_1} + V_{S_1,0} (x_S, \beta_{S_1}) \Delta 
    + \mathcal{L} V_{S_1, 0} (x, \theta) \tfrac{\Delta^2}{2}
    + \mathcal{L}^2 V_{S_1, 0} (x, \theta) \tfrac{\Delta^3}{6}   
    \\[0.3cm]
    x_{S_2} + V_{S_2,0} (x, \beta_{S_2}) \Delta 
    + \mathcal{L} V_{S_2, 0} (x, \theta) \tfrac{\Delta^2}{2}  \\[0.3cm] 
    x_R + V_{R, 0} (x, \beta_R) \Delta
\end{bmatrix}. 
\end{align*}   
%
Notice that the scheme involves $3d$ Gaussian random variables: 
\begin{align*} 
B_{j, i+1} - B_{j, i}, \quad  
\int_{t_i}^{t_{i+1}} \int_{t_i}^u d B_{j,  v} du, \quad 
\int_{t_i}^{t_{i+1}} \int_{t_i}^u \int_{t_i}^v  d B_{j, w} dv du, \qquad  1 \le j \le d. 
\end{align*}
The latter of the above integrals appears due to the  
application of a third-order stochastic Taylor expansion on $V_{S_1,0} (x_S, \beta_{S_1})$, in the smoothest component $\bar{X}_{S_1, i+1}^{\, (\mrm{II})}$.
%
%
As we will show in Section \ref{sec:main}, the log-likelihood based on the local Gaussian scheme (\ref{eq:lg_hypo_II}) produces a contrast estimator that is asymptotically unbiased in the high-frequency, complete observation regime. In order for the deduced contrast function to provide desirable asymptotic properties, it is required to include terms up to $\mathcal{O} (\Delta_n^3)$ in the definition of $\mu_{S_1}$, otherwise estimation of the parameter $\beta_R$ in the model (\ref{eq:hypo-II}) can be asymptotically biased as \cite{poke:09} observed for some bivariate hypo-elliptic diffusions in the framework of (\ref{eq:hypo-I}). 
\\ 

\noindent We denote by $\Sigma  (\Delta, x, \theta)$ the covariance matrix for one-step implementation of scheme $\bar{X}^{(\mrm{II})}$, with step-size $\Delta >0$, current state $x \in \mathbb{R}^N$ and parameter $\theta \in \Theta$. The covariance matrix is given as: 
\begin{align}  \label{eq:covariance_lg}
\Sigma (\Delta, x, \theta) 
=
\begin{bmatrix}
\Sigma_{S_1 S_1} (\Delta, x, \theta) & 
\Sigma_{S_1 S_2}  (\Delta, x, \theta) & 
\Sigma_{S_1 R}  (\Delta, x, \theta) \\[0.2cm]
\Sigma_{S_2 S_1}  (\Delta, x, \theta) & 
\Sigma_{S_2 S_2} (\Delta, x, \theta) & 
\Sigma_{S_2 R } (\Delta, x, \theta) \\[0.2cm]
\Sigma_{R S_1} (\Delta, x, \theta) & 
\Sigma_{R S_2}  (\Delta, x, \theta) &
\Sigma_{R R}  (\Delta, x, \theta) 
\end{bmatrix},    
\end{align}
where each block matrix is specified as: for $x=(x_{S_1}, x_{S_2}, x_R) \in \mathbb{R}^N$, $\theta = (\beta_{S_1}, \beta_{S_2}, \beta_R, \sigma) \in \Theta $, 
\begin{align*}
\begin{aligned} 
   & \Sigma_{S_1 S_1}  (\Delta, x, \theta)
   \equiv \tfrac{\Delta^5}{20} \, a_{S_1} (x, \theta), \quad  
   \Sigma_{S_1 S_2}  (\Delta, x, \theta)
   \equiv \tfrac{\Delta^4}{8} \, 
   \partial_{x_{S_2}}^\top  V_{S_1, 0} (x_{S}, \beta_{S_1}) \, 
   a_{S_2} (x, \theta);  \\[0.3cm]
   & \Sigma_{S_1 R}  (\Delta, x, \theta)
   \equiv \tfrac{\Delta^3}{6}
    \partial_{x_{S_2}}^\top  V_{S_1, 0} (x_{S}, \beta_{S_1}) 
   \partial_{x_R}^\top V_{S_2, 0} (x, \beta_{S_2}) 
   a_R (x, \sigma), \quad  
   \Sigma_{S_2 S_1}  (\Delta, x, \theta) 
   \equiv  \Sigma_{S_1 S_2}  (\Delta, x, \theta)^\top; \\[0.3cm]
   & \Sigma_{S_2 S_2}  (\Delta, x, \theta) 
   \equiv \tfrac{\Delta^3}{3} a_{S_2} (x, \theta), 
   \quad \Sigma_{S_2 R} (\Delta, x, \theta) 
   \equiv \tfrac{\Delta^2}{2}  
    \partial_{x_R}^\top V_{S_2, 0} (x, \beta_{S_2}) 
   \, a_R (x, \sigma), \quad 
   \Sigma_{R S_1} (\Delta, x, \theta) \equiv  \Sigma_{S_1 R} (\Delta, x, \theta)^\top;  \\[0.3cm] 
   & \Sigma_{R S_2} (\Delta, x, \theta) 
   \equiv  \Sigma_{S_2 R} (\Delta, x, \theta)^\top,  \quad 
   \Sigma_{RR} (\Delta, x, \theta) \equiv \Delta \, 
   a_R (x, \sigma).  
\end{aligned}
\end{align*}
In the above, we have set
\begin{align*}
    a_R (x, \sigma) & = \sum_{k = 1}^d V_{R, k} (x, \sigma) V_{R, k} (x, \sigma)^\top, \quad 
    a_{S_2} (x,\theta) = \partial_{x_R}^\top V_{S_2, 0} (x, \beta_{S_2}) \, a_R (x, \sigma) \,  \bigl( \partial_{x_R}^\top V_{S_2, 0} (x, \beta_{S_2}) \bigr)^\top;   \\[0.2cm] 
    a_{S_1} (x, \theta) & = 
    \partial_{x_{S_2}}^\top  V_{S_1, 0} (x_{S}, \beta_{S_1}) \, a_{S_2} (x, \theta) \, 
    \bigl( \partial_{x_{S_2}}^\top  V_{S_1, 0} (x_{S}, \beta_{S_1}) \bigr)^\top. 
\end{align*}
{
\begin{proposition} \label{prop:positive_definite}
Under condition \ref{assump:hypo}-II, it holds that: 
\begin{align} \label{eq:inf_dets}
\inf_{(x, \sigma) \in \mathbb{R}^N \times \Theta_\sigma} 
| a_{R} (x, \sigma) | > 0, 
\qquad 
\inf_{(x, \theta) \in \mathbb{R}^N \times \Theta} 
| a_{S_2} (x, \theta) | > 0, 
\qquad 
\inf_{(x, \theta) \in \mathbb{R}^N \times \Theta} 
| a_{S_1} (x, \theta) | > 0.   
\end{align} 
Then, it also holds that: for any $\Delta > 0$, 
\begin{align} \label{eq:det_cov}
\inf_{(x, \theta) \in \mathbb{R}^N \times \Theta} 
| \Sigma (\Delta, x , \theta) | > 0. 
\end{align} 
%
%
\end{proposition}
\noindent The proof is given in Appendix \ref{appendix:positive_definite}. Due to Proposition \ref{prop:positive_definite}, the covariance $\Sigma (\Delta, x, \theta)$ is invertible uniformly in $\Delta, x$ and $\theta$, thus the approximate log-likelihood based on the local Gaussian discretisation scheme (\ref{eq:lg_hypo_II}) is well-defined for the highly degenerate class (\ref{eq:hypo-II}). We note that, in brief,  result (\ref{eq:det_cov}) follows from (\ref{eq:inf_dets}) under condition \ref{assump:hypo}-II. 
} 
\begin{remark} \label{rem:matrix}
{
For elliptic diffusions, \cite{kess:97, uchi:12} assumed 
$\textstyle \inf_{(x, \sigma) \in \mathbb{R}^N \times \Theta_\sigma} | a_{R} (x, \sigma) | > 0$, and then established asymptotic properties of developed contrast estimators by introducing additional technical conditions. 
For the degenerate diffusion class (\ref{eq:hypo-II}), we work with condition \ref{assump:hypo}-II that is sufficient for (\ref{eq:inf_dets}). We make use of (\ref{eq:inf_dets}) in the proof of main results, i.e.~when studying the asymptotic properties of the contrast estimator proposed later on.}
\end{remark}
\section{Parameter Inference for Class (\ref{eq:hypo-II})}
\label{sec:main}
We explore analytically parameter inference procedures for hypo-elliptic diffusions in class (\ref{eq:hypo-II}). We prove in Section \ref{sec:complete} that a contrast estimator constructed from the conditionally Gaussian discretisation scheme (\ref{eq:lg_hypo_II}) is asymptotically unbiased under the high-frequency, complete observation regime. We illustrate the precise impact of the drift terms involved in 
scheme (\ref{eq:lg_hypo_II}) on the asymptotic results. In Section \ref{sec:partial}, we consider the partial observation regime. As observed for the case of class (\ref{eq:hypo-I}) in the literature, we show via analytical case studies that use of finite-differences for the estimation of hidden paths leads to biased parameter estimates within (\ref{eq:hypo-II}). Also, we explain that the local Gaussian scheme can be put into  effective use within  computational approaches for filtering hidden components and estimating parameters. 
\subsection{Complete Observation Regime} 
\label{sec:complete}
\subsubsection{Contrast Estimator} 
Based on the proposed scheme (\ref{eq:lg_hypo_II}) and the corresponding tractable transition density, we construct a contrast estimator for the hypo-elliptic class (\ref{eq:hypo-II}). We write the transition density of the local Gaussian scheme (\ref{eq:lg_hypo_II}), for given $\Delta > 0$, current position $x \in \mathbb{R}^N$ and parameter $\theta \in \Theta$ as:
\begin{align*} 
y  \mapsto \bar{p}_\Delta (x, y ; \theta) 
= \frac{1}{\sqrt{ (2 \pi)^N \Delta^{5 N_{S_1} + 3 N_{S_2} + N_R}
| \Sigma  (x, \theta) | }} 
\exp \Big( -\tfrac{1}{2} m (\Delta, x, y, \theta)^\top
    \Sigma^{-1}  (x, \theta) 
   m (\Delta, x, y, \theta) \Big),  
\end{align*}
where we have set, for $y = (y_{S_1}, y_{S_2}, y_R) \in \mathbb{R}^{N_{S_1}} \times \mathbb{R}^{N_{S_2}} \times \mathbb{R}^{N_{R}}$, 
\begin{align*} 
   m (\Delta, x, y, \theta) 
  =
  \begin{bmatrix}
  \tfrac{1}{\sqrt{\Delta^{5}}} \bigl(y_{S_1} - {\mu}_{S_1} (\Delta, x, \theta) \bigr) \\[0.3cm]  
  \tfrac{1}{{\sqrt{\Delta^{3}}}} 
  \bigl(y_{S_2} - {\mu}_{S_2} (\Delta, x, \theta) \bigr) \\[0.3cm] 
  \tfrac{1}{\sqrt{\Delta}}  \bigl(y_{R} - {\mu}_{R} (\Delta, x, \theta) \bigr) 
  \end{bmatrix}, \qquad  
  \Sigma (x, \theta) = \Sigma (1, x, \theta).  
\end{align*}
Note that from Proposition \ref{prop:positive_definite}, under Condition \ref{assump:hypo}-II the covariance matrix $\Sigma (x, \theta)$ is invertible for any $(x, \theta) \in \mathbb{R}^N \times \Theta$. We denote by $X^{\, (\mrm{II})}_i$ the (complete) observation of a diffusion within class (\ref{eq:hypo-II}) at time~$t_i$, $0\le i\le n$. Then, after removing some constant terms from $-2\,\textstyle{\sum_{i = 1}^n \log \bar{p}_{\Delta_n} ( X_{i-1}^{(\mrm{II})}, X_{i}^{(\mrm{II})}  ; \theta)}$, we define the following contrast function:  
\begin{align} 
\begin{aligned} \label{eq:contrast}  
\ell_{n}( \theta ) 
& =  \sum_{i = 1}^n  m (\Delta_n, 
\sample{i-1}^{(\mrm{II})}, \sample{i}^{(\mrm{II})} , \theta)^\top \, 
\Sigma^{-1} (\sample{i-1}^{(\mrm{II})}, \theta) \, 
m (\Delta_n, \sample{i-1}^{(\mrm{II})},  
\sample{i}^{(\mrm{II})}, \theta) 
+ \sum_{i = 1}^n \log \bigl| \Sigma 
(\sample{i-1}^{\, (\mrm{II})}, \theta) \bigr|.   
\end{aligned}
\end{align}
Thus, the contrast estimator for the hypo-elliptic class
(\ref{eq:hypo-II}) is defined as: 
\begin{align} \label{eq:MLE}
  \hat{\theta}_{n} 
  = \bigl( 
  \hat{\beta}_{S_1, n}, \, 
  \hat{\beta}_{S_2, n}, \,  
  \hat{\beta}_{R, n}, \, 
  \hat{\sigma}_{n}  \bigr) 
  = \mathrm{argmin}_{\theta \in \Theta} \; 
  \ell_{n} \bigl( \theta \bigr).   
\end{align}
\subsubsection{Asymptotic Results} \label{sec:condition_main}
Before we state our main results, we introduce some conditions for class (\ref{eq:hypo-II}). 
\begin{enumerate}[resume]
\item \label{assump:coeff}
{ 
$V_j \in C_p^{\infty} (\mathbb{R}^N \times \Theta; \mathbb{R}^N)$, $0 \le j \le d$.
} 
\item  \label{assump:bd_deriv}
For any $1 \le i \le N, \, 0 \le j \le d$ and $\alpha \in \{1, \ldots, N \}^l, \, l \ge 0$ the following function
\begin{align*} 
\theta \mapsto \partial^x_\alpha V_j^i (x , \theta)  
\end{align*}
is three times differentiable for all $x \in \mathbb{R}^N$. Furthermore, derivatives of the above map up to the third order are of polynomial growth in $x \in \mathbb{R}^N$ uniformly in $\theta \in \Theta$.  
%
%
\item \label{assump:moments}
The diffusion process $\{X_t\}_{t \geq 0}$ defined via  (\ref{eq:hypo-II}) is ergodic under $\theta = \trueparam$, with invariant distribution $\truedist$ on $\mathbb{R}^N$. Furthermore, all moments of $\truedist$ are finite. 
\item \label{assump:finite_moment}
It holds that for all $p \ge 1$, 
$\textstyle \sup_{t > 0} \mathbb{E}_{\trueparam} [|X_t|^p] < \infty$. 
\item \label{assump:ident}
If it holds that
\begin{align*}
V_{S, 0} (x , \beta_S) = V_{S, 0} (x , \truebeta_S), \ \ 
V_{R, 0} (x , \beta_R) = V_{R, 0} (x , \truebeta_R), \ \  
V_R (x , \sigma) = V_R (x , \truesigma), 
\end{align*}
for $x$ in set of probability 1 under  $\nu_{\trueparam}$, then 
$\beta_S = \truebeta_S$,  $\beta_R = \truebeta_R$,  $\sigma = \truesigma$.
\end{enumerate}
%
We write the true value of the parameter for a model in (\ref{eq:hypo-II}) as 
$ 
\theta^{\dagger} = (
\beta_{S_1}^{\dagger},  
\beta_{S_2}^{\dagger},     
\beta_{R}^{\dagger},  
\sigma^{\dagger} ) \in \Theta
$. 
The latter, $\trueparam$, is assumed to lie in the interior of $\Theta$. 
Recall the definition of function $V_0 : \mathbb{R}^N \times \Theta \to \mathbb{R}^N$ in (\ref{eq:V0}). Then, the contrast estimator defined in (\ref{eq:MLE}) has the following asymptotic properties in the high-frequency observation setting.
\begin{theorem}[Consistency] \label{thm:consistency}
Assume that conditions \ref{assump:hypo}-II, (\ref{assump:coeff})--(\ref{assump:ident}) hold. If \limit, then
\begin{align*}
\hat{\theta}_{n} \probconv \trueparam. 
\end{align*}
\end{theorem}
\begin{theorem}[Asymptotic Normality]  \label{thm:clt} 
Assume that conditions \ref{assump:hypo}-II,  (\ref{assump:coeff})--(\ref{assump:ident}) hold. If \limit \, with $\Delta_n =  o ( n^{-1/2})$, then
\begin{align*}    
  \begin{bmatrix}
    \sqrt{\tfrac{n}{\Delta_n^3}} \bigl( \hat{\beta}_{S_1, n} - \truebeta_{S_1} \bigr) \\[0.3cm]   
    \sqrt{\tfrac{n}{\Delta_n}} 
    \bigl( \hat{\beta}_{S_2, n} - \truebeta_{S_2} \bigr) \\[0.3cm]
    \sqrt{n \Delta_n} \bigl( \hat{\beta}_{R, n} - \truebeta_{R} \bigr) \\[0.3cm]
    \sqrt{n} \bigl( \hat{\sigma}_{n} - \truesigma \bigr) 
  \end{bmatrix}
  \distconv \mathscr{N} \bigl( 
   \mathbf{0}_{N_\theta}, 
    \Gamma^{-1} ( \trueparam )  
   \bigr),  
\end{align*}
where the asymptotic precision matrix $\Gamma (\trueparam)$ is given as:
\begin{align} \label{eq:precision_matrix} 
\Gamma ( \trueparam )  
= \mrm{diag} \Bigl(
\Gamma^{\beta_{S_1}} (\trueparam), \,  
\Gamma^{\beta_{S_2}} (\trueparam), \,
\Gamma^{\beta_{R}} (\trueparam), \,
\Gamma^{\sigma} (\trueparam)
\Bigr), 
\end{align} 
with the involved block matrices 
$\Gamma^{\beta_{S_1}} (\trueparam)  \in \mathbb{R}^{N_{\beta_{S_1}} \times N_{\beta_{S_1}}}$, 
$\Gamma^{\beta_{S_2}} (\trueparam)  \in \mathbb{R}^{N_{\beta_{S_2}} \times N_{\beta_{S_2}}}$,
$\Gamma^{\beta_{R}} (\trueparam)  \in \mathbb{R}^{N_{\beta_{R}} \times N_{\beta_{R}}}$,
$\Gamma^{\sigma} (\trueparam)  \in \mathbb{R}^{N_{\sigma} \times N_{\sigma}}$ 
specified as: 
\begin{align*}
&  \Gamma_{ij}^{\beta_{S_1}} ( \trueparam )
=  720  \int  
\partial_i^{\beta_{S_1}} V_{S_1, 0} (x_S , \truebeta_{S_1})^\top 
\, a^{-1}_{S_1} (x, \trueparam)  \, 
\partial_{{j}}^{\beta_{S_1}} 
V_{S_1, 0} (x_S , \truebeta_{S_1})  \; 
 \nu_{\trueparam} (dx); \\[0.1cm]
 &  \Gamma_{ij}^{\beta_{S_2}} ( \trueparam )
= 12 \int 
\partial^{\beta_{S_2}}_{i} 
 V_{S_2, 0} (x , \truebeta_{S_2} )^\top 
\, a_{S_2}^{-1} (x, \trueparam )  \, 
\partial^{\beta_{S_2}}_{j} V_{S_2, 0} (x , \truebeta_{S_2} ) \; 
 \nu_{\trueparam} (dx); \\[0.1cm] 
&  \Gamma_{ij}^{\beta_R} ( \trueparam )
= \int 
\partial_{{i}}^{\beta_R}
 V_{R, 0} (x , \truebeta_R )^\top 
\, a_R^{-1}  (x ,  \truesigma)  
\, 
\partial_{{j}}^{\beta_R} V_{R, 0} (x , \truebeta_R) \; 
 \nu_{\trueparam} (dx); \\[0.1cm]
& \Gamma_{ij}^\sigma (\trueparam)  
= \tfrac{1}{2} \int \mathrm{tr} 
\bigl( 
\partial_{i}^\sigma 
\Sigma (x ,  \trueparam)
\Sigma^{-1}  (x ,  \trueparam)
\partial_{j}^\sigma
\Sigma (x ,  \trueparam)
\Sigma^{-1} (x ,  \trueparam) \bigr)   
\;  \nu_{\trueparam}(dx). 
\end{align*} 
\end{theorem} 
\noindent The proofs of Theorems \ref{thm:consistency} \& \ref{thm:clt} are given in Appendix \ref{sec:pf}. 
\begin{remark} \label{rem:pf_main}
\cite{glot:20}
prove consistency and asymptotic normality of a contrast estimator constructed via a locally Gaussian scheme for  class~(\ref{eq:hypo-I}). 
Our proofs follow a different approach from the one in the above work. 
Indicatively, a condition on the step-size of $\Delta_n = o (n^{-1/2})$ is not required for our proof of consistency, whereas it is needed in \cite{glot:20}. 
We make use of preliminary convergence rates for the estimators $\hat{\beta}_{S_1,n}$, $\hat{\beta}_{S_2,n}$ and some key identities in the involved matrix calculations to control terms arising in the expansion of the logarithm of the contrast function -- this method then provides the consistency for $\hat{\beta}_{R,n}$.
{A detailed strategy to show the consistency can be found in Appendix \ref{sec:pf_consistency}. Also, the strategy used in our proofs can be applied to show asymptotic properties of a contrast estimator for a much wider class of highly degenerate diffusions such that condition \ref{assump:hypo}-I or II does not hold, and higher order Lie-brackets, e.g., $[\widetilde{V}_0, [\widetilde{V}_0,[\widetilde{V}_0, V_j ]]]$, are required for the uniform H\"ormander's condition. For instance, let $m \ge 3$ and consider the following $(m+1)$-dimensional highly degenerate SDE:
\begin{align*}
    d X_{S_1, t} & = V_{S_1, 0} (X_{S_2, t}, \beta_{S_1}) dt; \\
    d X_{S_2, t} & = V_{S_2, 0} (X_{S_3, t}, \beta_{S_2}) dt; \\ 
    & \vdots \\
    d X_{S_{m-1}, t} & = V_{S_{m-1}, 0} (X_{S_m, t}, \beta_{S_{m-1}}) dt; \\ 
    d X_{S_m, t} & = V_{S_m, 0} (X_{R, t}, \beta_{S_m}) dt; \\ 
    d X_{R, t} & = V_{R, 0} (X_{R, t}, \beta_R) dt +  
    V_{R, 1} (X_{R, t}, \sigma) dB_{1,t}, 
\end{align*}
where $V_{S_i, 0}:\mathbb{R} \times \Theta_{\beta_{S_i}} \to \mathbb{R}, \, 1 \le i \le m$, $V_{R, 0} : \mathbb{R} \times \Theta_{\beta_R} \to \mathbb{R}$ 
and $V_{R, 1} : \mathbb{R} \times \Theta_{\sigma} \to \mathbb{R}$ with some compact parameter spaces $\Theta_{\beta_{S_i}}$, $\Theta_{\beta_R}$ and $\Theta_\sigma$. We assume a similar condition to \ref{assump:hypo} so that the uniform H\"ormander's condition (\ref{eq:inf_H}) holds with $M=m$.  We provide a recipe to construct an asymptotically unbiased contrast estimator for such a general class of degenerate diffusions in Section \ref{sec:end}. In brief, given the deduced contrast estimator, the proof of consistency will start by showing the consistency of  $\hat{\beta}_{S_1, n}$ and obtaining a convergence rate. Then, utilising this result, we would next prove the consistency of $\hat{\beta}_{S_2, n}$ and get its convergence rate. We would repeat this procedure until we show the consistency of the estimators for $(\beta_R, \sigma)$.  
}

%
%
\end{remark}
\subsubsection{Case Study --  Bias due to Incorrect Drift Expansion} \label{sec:case_drift}
We have proven that the contrast estimator based on the proposed scheme (\ref{eq:lg_hypo_II}) is asymptotically unbiased under the high-frequency, complete observation regime. The inclusion of an appropriate number terms from the stochastic Taylor expansion of $V_{S_1, 0} (X_{S_1, t}, X_{S_2, t}, \beta_{S_1})$ and $V_{S_2, 0} (X_{t}, \beta_{S_2})$ in scheme (\ref{eq:lg_hypo_II}) is critical for obtaining desirable asymptotic properties.
Omission of such terms will typically give rise to an asymptotic bias. 
In this subsection, we briefly highlight the effect of the `drift correction' via a simple three-dimensional hypo-elliptic model from (\ref{eq:hypo-II}). 
We consider the following SDE: 
\begin{align}
\begin{aligned} \label{eq:toy-3}  
 d q_ t & = p_t dt; \\ 
 d p_t & = s_t dt ; \\ 
 d s_t & = - \beta s_t dt + \sigma dB_t, 
\end{aligned} 
\end{align}
where $\theta = (\beta, \sigma)$ is the parameter vector. We assume that all components of the system are observed, and consider the following discretisation scheme for SDE (\ref{eq:toy-3}): 
\begin{align} \label{eq:toy_app_1}
\begin{aligned}
  \bar{x}_{i+1} = 
\begin{bmatrix}
  \bar{q}_{i+1} \\[0.2cm]
  \bar{p}_{i+1} \\[0.2cm] 
  \bar{s}_{i+1} 
\end{bmatrix}
= 
\begin{bmatrix}
\bar{q}_i + \bar{p}_i \Delta_n + \bar{s}_i \tfrac{\Delta_n^2 }{2} 
\\[0.3cm] 
\bar{p}_i + \bar{s}_i \Delta_n 
- \beta \bar{s}_i  \tfrac{\Delta^2_n}{2} \\[0.3cm] 
\bar{s}_i - \beta \bar{s}_i \Delta_n  
\end{bmatrix}
+ 
\begin{bmatrix}
\sigma \times \int_{t_i}^{t_{i+1}} \int_{t_i}^u \int_{t_i}^v dB_w dv du \\[0.3cm] 
\sigma \times \int_{t_i}^{t_{i+1}} \int_{t_i}^u d B_{v} du \\[0.3cm] 
\sigma \times ( B_{t_{i+1}} - B_{t_i} )
\end{bmatrix}. 
\end{aligned}
\end{align}
Thus, terms of size $\mathcal{O} (\Delta_n^3)$ are not included in the deterministic part of the approximation $\bar{q}_{i+1}$ of the smoothest component.
Based on the conditionally Gaussian scheme (\ref{eq:toy_app_1}), we define a contrast function as $(- 2) \times$ (complete log-likelihood), that is, after some constants are removed: 
\begin{align*}
\ell_{n} (\theta ; x_{0:n})  
= 6n \log \sigma +  \tfrac{1}{\sigma^2} \sum_{i = 1}^n m(\Delta_n, x_{i-1}, x_i, \theta)^\top 
\, 
\Sigma^{-1} \, 
m (\Delta_n, x_{i-1}, x_i, \theta), 
\end{align*}
where we have set $x_i = [q_i, p_i, s_i ]^\top$,  $0 \le i \le n$, and 
\begin{align*}
m (\Delta_n, x_{i-1}, x_i, \theta) 
=  
\begin{bmatrix}
\tfrac{1}{\sqrt{\Delta_n^5}} \bigl(q_{i} -  q_{i-1} - p_{i-1} \Delta_n -  s_{i-1} \tfrac{\Delta_n^2 }{2} \bigr) \\[0.3cm] 
\tfrac{1}{\sqrt{\Delta_n^3}} \bigl(p_{i} - p_{i-1} - s_{i-1} \Delta_n + \beta s_{i-1} \tfrac{\Delta^2_n}{2} \bigr) \\[0.3cm]
\tfrac{1}{\sqrt{\Delta_n}} \bigl( s_{i} - s_{i-1} + \beta s_{i-1} \Delta_n  \bigr) 
\end{bmatrix}
, \quad  
{\Sigma} 
=         
\begin{bmatrix}
 \tfrac{1}{20} & \tfrac{1}{8} & \tfrac{1}{6} \\[0.2cm]
 \tfrac{1}{8} & \tfrac{1}{3} & \tfrac{1}{2} \\[0.3cm]  
 \tfrac{1}{6} & \tfrac{1}{2} & 1
\end{bmatrix}.   
\end{align*}
%
%
%
Solving $\partial_\beta \, \ell_n (\theta; x_{0:n}) = 0$, we obtain the contrast estimator for $\beta$ as $\widetilde{\beta}_n = {\widetilde{g}_n}/{\widetilde{f}_n}$, where we have defined:   
\begin{align*}  
 \widetilde{f}_n  & = 
 \Bigl\{ \tfrac{1}{2} \Sigma^{-1}_{32}  + \Sigma^{-1}_{33}  
+ \tfrac{1}{4} \Sigma^{-1}_{22}  + \tfrac{1}{2} \Sigma^{-1}_{23} \Bigr\} \times \tfrac{1}{n}  \sum_{i = 1}^n s_{i-1}^2; \\ 
 \widetilde{g}_n & = - \tfrac{1}{n \sqrt{\Delta_n}} \sum_{i = 1}^n s_{i-1} \times \Bigl\{ 
 \Sigma^{-1}_{31} \times \tfrac{q_{i} - q_{i-1} - p_{i-1} \Delta_n - s_{i-1} \tfrac{\Delta_n^2}{2} }{\sqrt{\Delta_n^5}} 
 \; + \;  \Sigma^{-1}_{32} \times \tfrac{p_{i} - p_{i-1} - s_{i-1} \Delta_n}{\sqrt{\Delta_n^3}}  + \Sigma^{-1}_{33} \times \tfrac{s_{i} - s_{i-1}}{\sqrt{\Delta_n}}  \\ 
 & \quad + \tfrac{1}{2} \Sigma^{-1}_{21} 
 \times \tfrac{q_{i} - q_{i-1} - p_{i-1} \Delta_n - s_{i-1} \tfrac{\Delta_n^2}{2} }{\sqrt{\Delta_n^5}} 
 + \tfrac{1}{2} \Sigma^{-1}_{22}  
 \times \tfrac{p_{i} - p_{i-1} - s_{i-1} \Delta_n}{\sqrt{\Delta_n^3}} + \tfrac{1}{2} \Sigma^{-1}_{23} 
 \times \tfrac{s_{i} - s_{i-1}}{\sqrt{\Delta_n}} 
 \Bigr\}.  
\end{align*}
From the ergodicity of the process $\{s_t\}$ and Lemma \ref{lemma:ergodic_thm} in the Appendix, we have that as \limit, 
\begin{align*} 
\widetilde{f}_n \probconv c_1 \times 
\int s^2 \, \truedist (ds), 
\end{align*}
for a non-zero constant $c_1 =  {\Sigma}_{32}^{-1}/2 
+ {\Sigma}_{33}^{-1}   
+  {\Sigma}_{22}^{-1}/4
+ {\Sigma}_{23}^{-1}/2$, where $\nu_\trueparam(ds)$ is the invariant distribution of $\{ s_t \}$ under the true parameter $\trueparam$. 
For the numerator $\widetilde{g_n}$, we apply Lemmas \ref{lemma:ergodic_thm} \& \ref{lemma:canonical_conv} in Appendix to obtain that 
\begin{align*} 
\widetilde{g}_n \probconv (c_1 + c_2) \times \truebeta \times 
\int s^2 \, \truedist (ds), 
\end{align*}
for a non-zero constant $c_2 \equiv  {\Sigma}_{31}^{-1}/6  
 +  {\Sigma}_{21}^{-1} / 12$. Hence, it holds that, if \limit, then
\begin{align*}  
\widetilde{\beta}_n 
\probconv \bigl( 1 + \tfrac{c_2}{c_1} \bigr) \times \truebeta. 
\end{align*}
Thus, the drift estimation based on the discretisation scheme (\ref{eq:toy_app_1}) with inappropriate drift expansion is, in general, asymptotically biased. One can check that the above bias is removed upon use of our locally Gaussian scheme (\ref{eq:lg_hypo_II}) instead of (\ref{eq:toy_app_1}). 
\subsection{Partial Observation Regime} 
\label{sec:partial}
\subsubsection{Case Study -- Bias due to Finite-Differences} \label{sec:case_2} 
To motivate the `appropriateness' of the proposed locally Gaussian scheme (\ref{eq:lg_hypo_II}) in the context of parameter inference in a partial observation setting, we illustrate that naive use of finite-differences to impute hidden components (a quite common in applications) induces a bias in the estimation of the SDE parameters. To observe this, consider the model (\ref{eq:toy-3}) again but with the drift parameter $\beta$ fixed to 1: 
\begin{align}
\begin{aligned} \label{eq:toy-2}
   d q_t & = p_t dt;  \\ 
   d p_t & =  s_t dt; \\ 
   d s_t & = - s_t dt + \sigma d B_t,  
\end{aligned}
\end{align}
for $\sigma > 0$.  We then apply the Euler-Maruyama scheme for the first equation of (\ref{eq:toy-2}) and the locally Gaussian scheme (\ref{eq:lg_hypo_I}) for the remaining dynamics, i.e., 
\begin{align}
\begin{aligned} \label{eq:toy_app_2}
\begin{bmatrix}
 \bar{q}_{i+1} \\[0.1cm]
 \bar{p}_{i+1} \\[0.1cm]
 \bar{s}_{i+1}
\end{bmatrix}
=
 \begin{bmatrix}
 \bar{q}_{i} +  \bar{p}_{i} \Delta_n \\[0.1cm]
 \bar{p}_{i} +  \bar{s}_{i} \Delta_n 
- \bar{s}_{i} \tfrac{\Delta_n^2}{2} \\[0.1cm]
 \bar{s}_{i} - \bar{s}_i \Delta_n  
 \end{bmatrix}
 + 
 \begin{bmatrix} 
 0 \\[0.3cm] 
 \sigma \times \int_{t_i}^{t_{i+1}} \int_{t_i}^{u} 
 d B_{v} du \\[0.3cm]
 \sigma \times \bigl( B_{i+1} - B_{i} \bigr)
 \end{bmatrix}. 
\end{aligned}
\end{align} 
Scheme (\ref{eq:toy_app_2}) is  degenerate since the upper-most equation does not involve noise.  
We now consider the estimator based on the likelihood provided by (\ref{eq:toy_app_2}), given the discrete-time observations $\{q_{0:n}, s_{0:n} \}$, and with the hidden paths $p_{0:n}$ imputed via the first equation of (\ref{eq:toy_app_2}) using the observations $q_{0:n}$. 
\begin{remark}
In practice, the rough component $s_t$ is often not observed, so one must impute the missing components $s_{0:n}$ conditionally on the observations $q_{0:n}$ by making use of the transition density (or some approximation of it) for both coordinates $(p,s)$ of (\ref{eq:toy-2}). 
One would reasonably expect that presence of bias will be typical in such a practical scenario, if it is found to be present in the simpler case when $s_t$ is directly observed.
\end{remark} 
%
\noindent The complete likelihood of the discretisation scheme (\ref{eq:toy_app_2}) is given as: 
\begin{align*} 
\begin{aligned}
    & \prod_{i =1}^n \Biggl\{ \tfrac{1}{\sqrt{(2 \pi)^2 \Delta_n^4 \sigma^4 | \Sigma | }} 
    \exp \Bigl(  
    - \frac{1}{2 \sigma^2} \, 
    m (\Delta_n, y_{i-1}, y_{i})^\top \, 
    \Sigma^{-1}  \,
    m (\Delta_n, y_{i-1}, y_{i}) 
    \Bigr) \times  \delta \bigl( q_{i} - q_{i-1} - p_{i-1} \Delta_n \bigr) \Biggr\}, 
\end{aligned} 
\end{align*}
where we have defined $y_i = [p_i, s_i]^\top$, $0 \le i \le n$, and  
\begin{align*}
m (\Delta_n, y_{i-1}, y_i) = 
\begin{bmatrix} 
  \tfrac{1}{\sqrt{\Delta_n^3}} 
  \bigl( p_i - p_{i-1} - s_{i-1} \Delta_n + s_{i-1} \tfrac{\Delta_n^2}{2} \bigr) \\[0.2cm] 
  \tfrac{1}{\sqrt{\Delta_n}} \bigl( s_i - s_{i-1} + s_{i-1} \Delta_n \bigr)
\end{bmatrix}, 
\quad 
\Sigma = 
\begin{bmatrix}
 \tfrac{1}{3}  &  \tfrac{1}{2} \\[0.2cm]  
 \tfrac{1}{2} &  1
\end{bmatrix}. 
\end{align*}  
Integrating out $p_{0:n}$,  we obtain the marginal likelihood $f_n (\sigma \, ; q_{0:n+1}, s_{0:n})$ as: 
\begin{align*}
f_n (\sigma \, ; q_{0: n + 1}, s_{0:n}) 
=  \prod_{i =1}^n \Biggl\{ \tfrac{1}{\sqrt{(2 \pi)^2 \Delta_n^4 \sigma^4 | \Sigma | }} 
\exp \Bigl(  - \frac{1}{2 \sigma^2} \, 
m (\Delta_n, \hat{y}_{i-1}, \hat{y}_{i} )^\top \, 
\Sigma^{-1}  \,
m ( \Delta_n, \hat{y}_{i-1}, \hat{y}_{i} ) 
\Bigr) \Biggr\}, 
\end{align*}
where $\hat{y}_i = [\hat{p_i}, s_i]^\top$, with $\hat{p}_i = \bigl( q_{i+1} - q_i \bigr) / \Delta_n$. 
%
Then, we obtain the following contrast function for $\sigma$, after removing constant terms from $(-2) \times \log f_n (\sigma ;  q_{0:n+1}, s_{0:n})$: 
\begin{align} \label{eq:toy2_cont}
\ell_{n} (\theta ; q_{0: n+1}, s_{0:n}) 
= 4 n \log \sigma 
+ \tfrac{1}{\sigma^2} \sum_{i = 1}^n 
m (\Delta_n, \hat{y}_{i-1}, \hat{y}_i )^\top  \, 
{\Sigma}^{-1} \, 
m (\Delta_n, \hat{y}_{i-1}, \hat{y}_i).  
\end{align}
Solving $\partial_\sigma  \ell_n (\sigma ; q_{0:n+1}, s_{0:n}) = 0 $, we obtain the estimator $\hat{\sigma}_n$, such that:
\begin{align*} 
(\hat{\sigma}_n)^2 = \tfrac{1}{2n} \sum_{i=1}^n 
m (\Delta_n, \hat{y}_{i-1}, \hat{y}_i)^\top 
\, {\Sigma}^{-1} \, 
m (\Delta_n, \hat{y}_{i-1}, \hat{y}_i). 
\end{align*}
%
It holds that, if \limit, then 
\begin{align} \label{eq:conv_toy_2}
 ( \hat{\sigma}_n )^2 \probconv \tfrac{8}{5} (\truesigma)^2. 
\end{align}
We prove convergence (\ref{eq:conv_toy_2}) in Appendix \ref{appendix:pf_case_study}. The proof indicates that the bias arises from the  higher order stochastic Taylor expansion terms of $q_{i+1} $ which are ignored by the estimate $\hat{p}_i$. Thus, use of the finite-difference approximation 
for imputation of component $p_t$, 
induces an asymptotic bias at the estimation of $\sigma$. 
%
\subsubsection{Filtering and Parameter Inference via the Proposed Scheme (\ref{eq:lg_hypo_II})}
\label{sec:kalman_filter}
We put forward the locally Gaussian scheme (\ref{eq:lg_hypo_II}) for imputing hidden components and performing parameter inference under a partial observation regime. The scheme and its transition density on the full set of coordinates can be combined with various computational methods, e.g., Monte-Carlo Expectation Maximisation (MC-EM) and Markov Chain Monte Carlo (MCMC), similarly to earlier works \citep{dit:19, poke:09} that applied some conditionally Gaussian schemes for inference of specific hypo-elliptic models within the class (\ref{eq:hypo-I}). 

We now highlight the use of a relatively straightforward Kalman filter recursion for carrying out statistical inference once the locally Gaussian scheme is adopted,
for a rich sub-class of hypo-elliptic models, referred to here as  \emph{conditionally Gaussian non-linear systems}.
That is, the system is originally specified as a non-linear SDE but can be treated as a linear system given components that correspond to observations. For elliptic diffusions with such a structure, continuous-time filtering and smoothing have been investigated in engineering, see e.g.~Chapter 8 of \cite{ch:23}. Several important hypo-elliptic models used in applications fall within this sub-class, e.g., standard Langevin equations, Quasi-Markovian generalised Langevin Equations (\ref{eq:qgle-I}, \ref{eq:qgle-II}).  Here, our interest lies in the sub-class derived via the general model (\ref{eq:hypo-II}) once the constituent coefficients are specified as:  
\begin{gather*}
V_{S_1, 0} (x_S, \beta_{S_1}) 
= 
C_{\beta_{S_1}} x_{S_1}
+ \hat{C}_{\beta_{S_1}} x_{S_2}, \quad   
V_{S_2, 0} (x, \beta_{S_2}) 
= f_{S_2} (x_{S_1}, \beta_{S_1}) 
+ C_{\beta_{S_2}} x_H;  
\\[0.2cm]
V_{R, 0} (x, \beta_R) 
= f_{R} (x_{S_1}, \beta_{R}) 
+ C_{\beta_R} x_H, \quad 
V_{R, j} (x, \sigma) = f_{R, j} (\sigma), \quad 1 \le j \le d, 
\end{gather*} 
for $x = (x_{S_1}, x_{S_2}, x_R) = (x_{S_1}, x_{H}) \in \mathbb{R}^{N_{S_1}} \times \mathbb{R}^{N_{S_2} + N_R}$ and $\theta = (\beta_{S_1}, \beta_{S_2}, \beta_R, \sigma) \in \Theta$, where: (i) $f_{S_2}$, $f_R$, $f_{R,j}$ are vector-valued functions, allowed to be non-linear w.r.t. the state $x_{S_1}$; (ii) matrices 
\begin{gather*}
{C}_{\beta_{S_1}} \in \mathbb{R}^{N_{S_1} \times N_{S_1} },\qquad  
\hat{C}_{\beta_{S_1}} \in \mathbb{R}^{N_{S_1} \times N_{S_2} }, \qquad C_{\beta_{S_2}} \in \mathbb{R}^{N_{S_2} \times (N_{S_2} + N_R)}, \qquad  C_{\beta_{R}} \in \mathbb{R}^{N_{R} \times (N_{S_2} + N_R)}
\end{gather*}
are independent of the state $x$. Critically, given the observable component $x_{S_1}$, the drift functions are linear functions of the hidden component $x_H$. For the model with the above choice of coefficients, the locally Gaussian scheme (\ref{eq:lg_hypo_II}) writes as:  
\begin{align} \label{eq:scheme_linear}
\bar{X}_{i+1} 
=
\begin{bmatrix}
  \bar{X}_{S_1, i+1}  \\[0.2cm] 
  \bar{X}_{S_2, i+1}  \\[0.2cm] 
  \bar{X}_{R, i+1}  
\end{bmatrix}
=  b (\Delta_n, \bar{X}_{S_1, i}, \theta) 
+ A (\Delta_n, \bar{X}_{S_1, i}, \theta) 
\begin{bmatrix}
    \bar{X}_{S_2, i} \\[0.1cm]
    \bar{X}_{R, i}
\end{bmatrix}
+ w_i (\Delta_n, \theta), 
\end{align}
for functions 
$b : (0, \infty) \times \mathbb{R}^{N_{S_1}} \times \Theta \to \mathbb{R}^{N} $, 
$A : (0, \infty) \times \mathbb{R}^{N_{S_1}} \times \Theta \to \mathbb{R}^{N \times (N_{S_2} + N_R)}$ and an $N$-dimensional Gaussian variate $w_i (\Delta_n, \theta)$. Since the 
right-hand side of scheme (\ref{eq:scheme_linear}) is linear w.r.t.~the hidden components $\bar{X}_{S_2, i}, \, \bar{X}_{R, i}$ given the observed component $\bar{X}_{S_1, i}$, one can obtain Kalman filtering and smoothing recursions, and calculate the marginal likelihood for the observations $\bar{X}_{S_1, 0:n}$. We provide the closed form filtering and marginal likelihood calculations in Appendix \ref{appendix:kalman}. We use these tools  in the numerical experiments of parameter inference under the partial observation regime in Section \ref{sec:numerical_studies} that follows.  
\begin{remark} 
\cite{vr:22} studied parameter inference for the QGLE of first-type in (\ref{eq:qgle-I}), 
where they applied an Euler-Maruyama scheme to construct Kalman filtering and smoothing for the rough components $(p_t, s_t)$  given the velocity $p_t$, with values of the latter obtained (via finite-differences) from discrete observations of the position $q_t$. Then, they used  Kalman filtering and smoothing within an Expectation-Maximisation (EM) algorithm to estimate the parameters. However, as we have seen, such a finite-differences approach can induce bias in the estimation of the diffusion parameters. 
\end{remark}
\section{Numerical Studies} 
\label{sec:numerical_studies}
\subsection{Linear SDE in a Partial Observation Regime} \label{sec:num_case_study} 
We illustrate empirically, for an example SDE model, that parameter estimation via the proposed locally Gaussian scheme (\ref{eq:lg_hypo_II}) leads to asymptotically unbiased estimation under the partial observation regime. We also highlight the effect of the drift correction in the properties of the estimators. We again consider the model studied in Section \ref{sec:case_drift}, that is,  
\begin{align*} 
\begin{aligned}
 d q_t & = p_t dt; \\
 d p_t & = s_t dt; \\ 
 d s_t & = - \beta s_t dt + \sigma dB_t,
\end{aligned}
\end{align*}
where $\theta = (\beta, \sigma) \in \Theta = (0, \infty) \times (0, \infty)$ is the parameter vector. In agreement with practice,  we assume that only discrete observations of the smoothest component, $q_{0:n}$, are available, with an equidistant step-size $\Delta_n$. We compute the following two estimators based on two different discretisation schemes:  
\begin{align*} 
\hat{\theta}_{n, j} 
= (\hat{\beta}_{n, j}, \hat{\sigma}_{n, j})
= {\rm argmax}_{\theta \in \Theta} \log p_{j} (\theta ; q_{0:n}), \quad j = 1,2, 
\end{align*}
where $p_{1} (\theta ; q_{0:n})$ is the approximate likelihood of the observations as obtained by use of Kalman filter in the setting of our locally Gaussian scheme (\ref{eq:lg_hypo_II}): 
\begin{align*}
\begin{aligned}
\begin{bmatrix}
\bar{q}_{i+1} \\[0.1cm]
\bar{p}_{i+1} \\[0.1cm]
\bar{s}_{i+1}
\end{bmatrix}
= 
\begin{bmatrix}
\bar{q}_{i} + \bar{p}_i \Delta_n 
+ \bar{s}_i \tfrac{\Delta_n^2}{2} 
- \beta \bar{s}_i \tfrac{\Delta_n^3}{6}  \\[0.3cm]
\bar{p}_{i} + \bar{s}_{i} \Delta_n 
- \beta \bar{s}_i \tfrac{\Delta_n^2}{2} \\[0.3cm]
\bar{s}_{i} - \beta \bar{s}_{i} \Delta_n 
\end{bmatrix}
+ 
\begin{bmatrix}
    \sigma \times \int_{t_i}^{t_{i+1}} \int_{t_i}^u \int_{t_i}^v dB_{w} dv du \\[0.3cm] 
    \sigma \times \int_{t_i}^{t_{i+1}} \int_{t_i}^v dB_{v} du
    \\[0.3cm] 
    \sigma \times \bigl( B_{i+1} - B_i \bigr)
\end{bmatrix}, 
\end{aligned} 
\end{align*}
and $p_{2} (\theta ; q_{0:n})$ is a different approximate likelihood obtained in the setting of the following conditionally Gaussian scheme that omits higher-order correction terms of order $\mathcal{O}(\Delta_n^2)$ and $\mathcal{O}(\Delta_n^3)$ from the stochastic Taylor expansion of the drift function of component $p$ and $q$: 
\begin{align} \label{eq:incorrect}
\begin{aligned}
\begin{bmatrix}
\widetilde{q}_{i+1} \\[0.1cm]
\widetilde{p}_{i+1} \\[0.1cm]
\widetilde{s}_{i+1}
\end{bmatrix}
= 
\begin{bmatrix}
\widetilde{q}_{i} + \widetilde{p}_i \Delta_n \\[0.3cm]
\widetilde{p}_{i} + \widetilde{s}_{i} \Delta_n  \\[0.3cm]
\widetilde{s}_{i} - \beta \widetilde{s}_{i} \Delta_n 
\end{bmatrix}
+ 
\begin{bmatrix}
    \sigma \times \int_{t_i}^{t_{i+1}} \int_{t_i}^u \int_{t_i}^v dB_{w} dv du \\[0.3cm] 
    \sigma \times \int_{t_i}^{t_{i+1}} \int_{t_i}^v dB_{v} du
    \\[0.3cm] 
    \sigma \times \bigl( B_{i+1} - B_i \bigr)
\end{bmatrix}. 
\end{aligned} 
\end{align}
We generate $50$ independent realisations of the dataset $q_{0:n}$ by sub-sampling trajectories obtained from scheme (\ref{eq:lg_hypo_II}) with a small step-size $10^{-4}$. We have chosen the scheme because it is expected to have a better accuracy than other classical schemes (such as Euler-Maruyama scheme) due to the higher order stochastic Taylor expansion of drift functions. We consider the following three high-frequency scenarios for the data: 
\begin{itemize}
\item[] \textbf{Set I}. $n = 5 \cdot 10^5, \, \Delta_n = 10^{-3}, \, T (= n \Delta_n) =500$. 
\item[] \textbf{Set II}. $n = 10^6, \, \Delta_n = 5 \cdot10^{-4}, \, {T=500} 
$. 
\item[] \textbf{Set III}. $n = 10^7, \, \Delta_n = 10^{-3}, \, T = 10^4$.  
\end{itemize}
\begin{figure}
     \centering
     \begin{subfigure}[b]{0.32\textwidth}
         \centering
         \includegraphics[width=\textwidth]{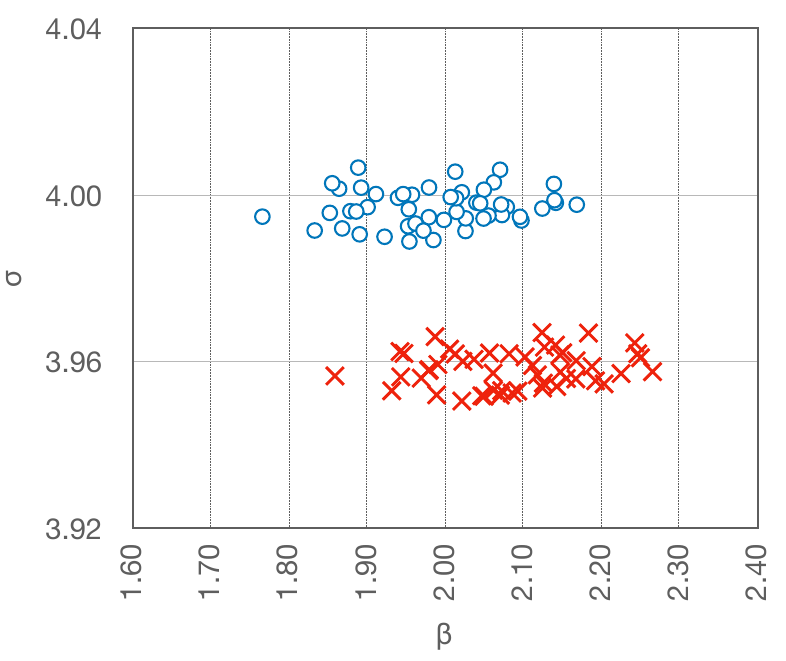}
        \caption{Set I.} %
         \label{fig:set_1}
     \end{subfigure}
     \hfill
     \begin{subfigure}[b]{0.32\textwidth}
         \centering
         \includegraphics[width=\textwidth]{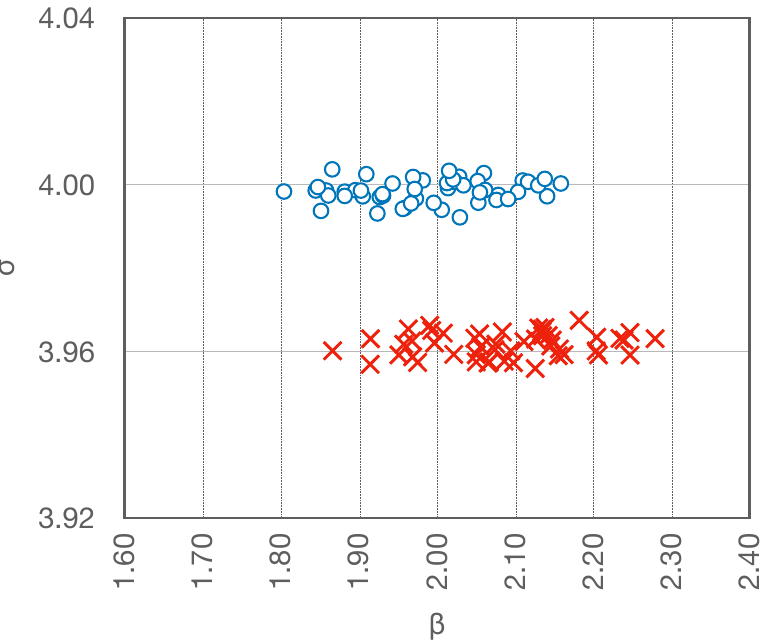}
        \caption{Set II.} 
         \label{fig:set_2}
     \end{subfigure} 
     \hfill
     \begin{subfigure}[b]{0.32\textwidth}
         \centering
         \includegraphics[width=\textwidth]{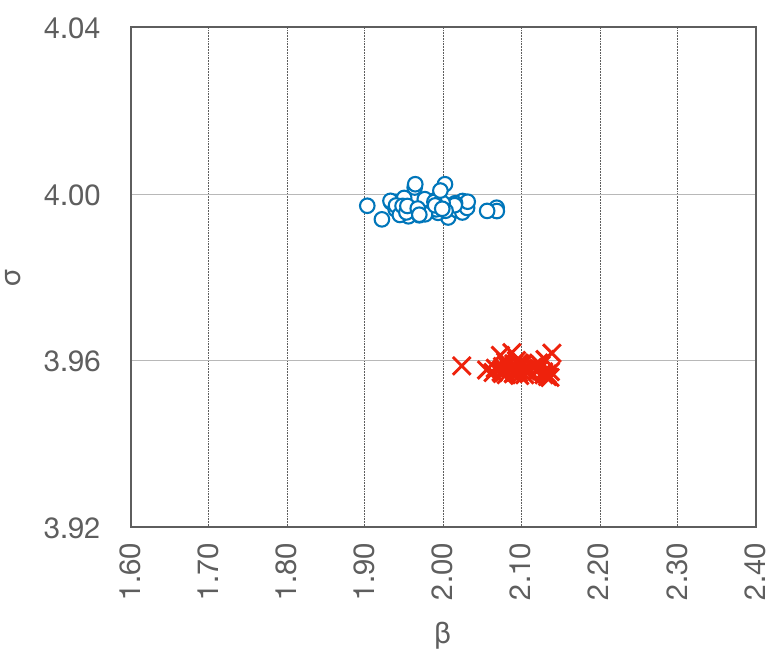}
         \caption{Set III.} 
         \label{fig:set_3}
     \end{subfigure}
        \caption{Parameter estimates from ${50}$ independent realisations of partial  observations: Each blue $\color{blue}{\circ}$ and red $\color{red}{\times}$ in the three figures represents one realisation of $\hat{\theta}_{n, 1}$ ({using (\ref{eq:lg_hypo_II})}) and $\hat{\theta}_{n, 2}$ ({using incorrect discretisation (\ref{eq:incorrect})}), respectively. The true value is $\trueparam = (\truebeta, \truesigma) = (2.0, 4.0)$. }
        \label{fig:mle}
\end{figure} 
The true parameter value is set to $\trueparam = (\truebeta, \truesigma) = (2.0, 4.0)$, and the Nelder-Mead method is applied to optimise the marginal likelihoods {with the initial guess set to $(\beta_0, \sigma_0) = (3.0, 3.0)$}.   
In Figure \ref{fig:mle}, we plot the ${50}$ realisations of the two different estimators. Table \ref{table:mle} summarizes the mean and standard deviation of {relative error of estimate, i.e., $(\hat{\theta}_{n, j} - \trueparam)/\trueparam, \, j = 1,2$}, from the ${50}$ repetitions. 
First, we observe that the estimates of $\hat{\theta}_{n ,1}$ (using our scheme (\ref{eq:lg_hypo_II})) are centred at the true value in all scenarios, thus in this case we have an empirical illustration of an asymptotically unbiased estimation in the partial observation setting. Secondly, it is clear from the figures and the table that the mean of estimates of $\hat{\theta}_{n,2}$ (estimator based on the conditionally Gaussian scheme without appropriate drift correction (\ref{eq:incorrect})) is shifted from the true value, and seems to be centred at $(\beta, \sigma) = (2.10, 3.960)$. Thus $\hat{\theta}_{n,2}$ induces an asymptotic bias in the partial observation regime, in agreement with the case study in Section \ref{sec:case_drift} in the complete observation case. 
Notably, there is a clear separation between the two estimators of  $\sigma$, {and the range of realisations of $\hat{\sigma}_{n, 2}$ does not cover the true value $\truesigma = 4.0$ for all scenarios. Comparing Set I and III in Figure \ref{fig:mle}, we observe that the variance of drift parameter estimation decreases as $T$ increases, and then the bias in $\hat{\beta}_{n, 2}$ becomes more clear in Set III.} We stress here that the bias in $\hat{\theta}_{n,2}$ is not removed with increasing $n$ or decreasing $\Delta_n$.  Also, one can still observe the bias even if the datasets are obtained with other numerical schemes, e.g., Euler-Maruyama scheme, rather than scheme (\ref{eq:lg_hypo_II}).    
%
%
\begin{table}
\centering
\caption{Mean and standard deviation (in parenthesis) of {$(\hat{\theta}_{n, j} - \trueparam)/\trueparam, \, j = 1,2$ from 
$50$} trajectories of partial observations.}
\begin{tabular}{@{}lcccc@{}}
\toprule
\\[-0.3cm]
\multirow{2}{*}{Set \hspace{0.2cm}}
& \multicolumn{2}{c}{
Proposed scheme (\ref{eq:lg_hypo_II}) \hspace{0.2cm}}
& \multicolumn{2}{c}{
Incorrect scheme (\ref{eq:incorrect})} \\ \cmidrule(l){2-3} \cmidrule(l){4-5}  
\\[-0.3cm]
& $(\hat{\beta}_{n, 1} - \truebeta)/\truebeta$
& $(\hat{\sigma}_{n, 1} - \truesigma)/\truesigma$
& $(\hat{\beta}_{n, 2} - \truebeta)/\truebeta$
& $(\hat{\sigma}_{n, 2} - \truesigma)/\truesigma$
\\[0.1cm] 
\midrule 
I.  
& -0.0052 (0.0463)	& -0.0007 (0.0011) 
& 0.0441 (0.0477) &-0.0104	(0.0011)  \\
II.   
& -0.0066	(0.0458)  & -0.0003	 (0.0007) 
& 0.0438 (0.0488) & -0.0096 (0.0007) \\ 
III.  
& -0.0074 (0.018) &	-0.0007	(0.0005) 
& 0.0496 (0.012) &	-0.0104	(0.0003)
\\
\bottomrule
\end{tabular}
\label{table:mle}
\end{table}
\subsection{Quasi-Markovian Generalised Langevin Equations} 
\subsubsection{Scalar Extended State} 
\label{sec:num_qgle}
We consider the QGLE describing one-dimensional positional domain: 
\begin{align} \label{eq:qgle_num_1} 
\begin{aligned} 
d q_t  & = p_t dt; \\[0.1cm]
d p_t  & = \bigl( - {U'}(q_t) + \lambda s_t  \bigr) dt;   \\[0.1cm]
d s_t  & =  (- \lambda p_t - \alpha  s_t) dt   + \sigma  d B_t,  \quad (q_0,p_0,s_0) \in \mathbb{R}^3, 
\end{aligned} 
\end{align}
where $\alpha > 0, \, \sigma > 0, \, \lambda \in \mathbb{R}\setminus \{0\}$ and $U : \mathbb{R} \to \mathbb{R}$. In this experiment, we consider the following two  choices of potential $U$: 
\begin{align*} 
q \mapsto U_{\mrm{HO}} (q) = D \times \tfrac{q^2}{2}, 
\qquad 
q \mapsto U_{\mrm{DW}} (q) = D \times \tfrac{q^2}{2} 
+ \sin \Bigl( \tfrac{1}{4} + 2q \Bigr),
\end{align*}
where $D > 0$ is a parameter. The function $U_{\mrm{DW}}$ (used in experiments in the work \cite{lei:22}) represents an uneven double well potential, under which model (\ref{eq:qgle_num_1}) is non-linear. We generate {50} independent datasets by sub-sampling trajectories produced by the discretisation scheme (\ref{eq:lg_hypo_II}) with a small step-size $10^{-4}$ so that obtained observations correspond to $n = 2 \times 10^5, \, \Delta_n = 10^{-3}, \, T = n\Delta_n= 200$. For the complete observation regime, we compute the contrast estimator (\ref{eq:MLE}) for each given dataset. For experiments of partial observations, we use the trajectories of the position $q_t$ only, extracted from the complete observations, and compute the MLE by maximising the marginal likelihood obtained from the Kalman recursion formula under the locally Gaussian scheme (\ref{eq:lg_hypo_II}), as shown in Section \ref{sec:kalman_filter}. To minimise 
the relevant target functions 
we use the adaptive moments (Adam) optimiser with the following algorithmic specifications: (step-size) = $0.2$,  (exponential decay rate for the first moment estimates) = $0.9$, 
(exponential decay rate for the first moment estimates) = $0.9$,  (exponential decay rate for the second moment estimates) = $0.999$, (additive term for numerical stability) = $10^{-8}$ and (number of iteration) = $2,000$. {The true parameters are set to $\trueparam = (D^{\dagger}, \lambda^{\dagger}, \alpha^{\dagger}, \sigma^{\dagger}) = (1.0, 2.0, 4.0, 4.0)$. 
Also, the initial guesses for the parameter are set to $(D_0, \lambda_0, \alpha_0, \sigma_0) = (3.0, 3.0, 3.0, 3.0)$.} 
{We summarise the mean and standard deviation of $(\mathrm{MLE} - \trueparam)/\trueparam$from the 50} independent trajectories in Table~\ref{table:qgle_second}. We notice that the results for the complete observation regime are in agreement with the analytical results in Theorem \ref{thm:clt}. For instance, convergence to the true values appears to be faster for parameters $(D, \lambda)$ in the smooth component $p_t$ (recall the convergence rate $\sqrt{\Delta_n / n}$ for such parameters in the CLT of Theorem \ref{thm:clt}). Besides, under the partial observation regime, the estimates seem to be centred around the true parameter as well, with standard deviations that are larger than the ones in the case of complete observations (as expected). Thus, parameter inference carried out via the proposed locally Gaussian scheme (\ref{eq:lg_hypo_II}) appears in this case to provide unbiased estimates in the partial observation regime.         
\begin{table} 
\caption{Parameter estimation of the QGLE (\ref{eq:qgle_num_1}). 
Mean and standard deviation (in brackets) of {$(\mathrm{MLE} - \trueparam)/\trueparam$ from $50$ trajectories of observations}.}
\centering
\begin{tabular}{@{}ccccc@{}}
\toprule 
\\[-0.2cm]
\multicolumn{1}{c}{\multirow{2}{*}{Potential}} 
& \multicolumn{1}{c}{\multirow{2}{*}{Parameter}} 
& \multicolumn{2}{c}{Relative error of estimates}   
\\[0.2cm] \cmidrule(l){3-4} 
\multicolumn{1}{c}{}   
& \multicolumn{1}{c}{}    
& \multicolumn{1}{c}{Complete observation} 
& Partial observation 
\\[0.2cm] 
\midrule 
\multirow{4}{*}{$U_{\mrm{HO}}$} 
& $D$  & -0.0000 (0.0001) &	-0.0037	 (0.1038)   \\
& $\lambda$  & 0.0000 (0.0001) & 0.0100	 (0.0533)  \\
& $\alpha$   & 0.0057	(0.0540) &	0.0091	(0.0546)  \\
& $\sigma$   &  -0.0000	 (0.0009) & -0.0079	 (0.0524)   \\[0.4cm] 
\multirow{4}{*}{$U_{\mrm{DW}}$} 
& $D$  &  0.0000 (0.0001) &	-0.0027	(0.1139) \\
& $\lambda$   & 0.0000	(0.0001) &	0.0201	(0.0638)  \\
& $\alpha$    & 0.0022	(0.0520) &	0.0059	(0.0530)   \\
& $\sigma$    & 0.0001	(0.0010) &	-0.0166	(0.0613)    \\
\bottomrule
\end{tabular} \label{table:qgle_second} 
\end{table}
\begin{figure}
    \centering
     \begin{subfigure}[b]{0.40\textwidth}
         \centering
     \includegraphics[width=\textwidth]{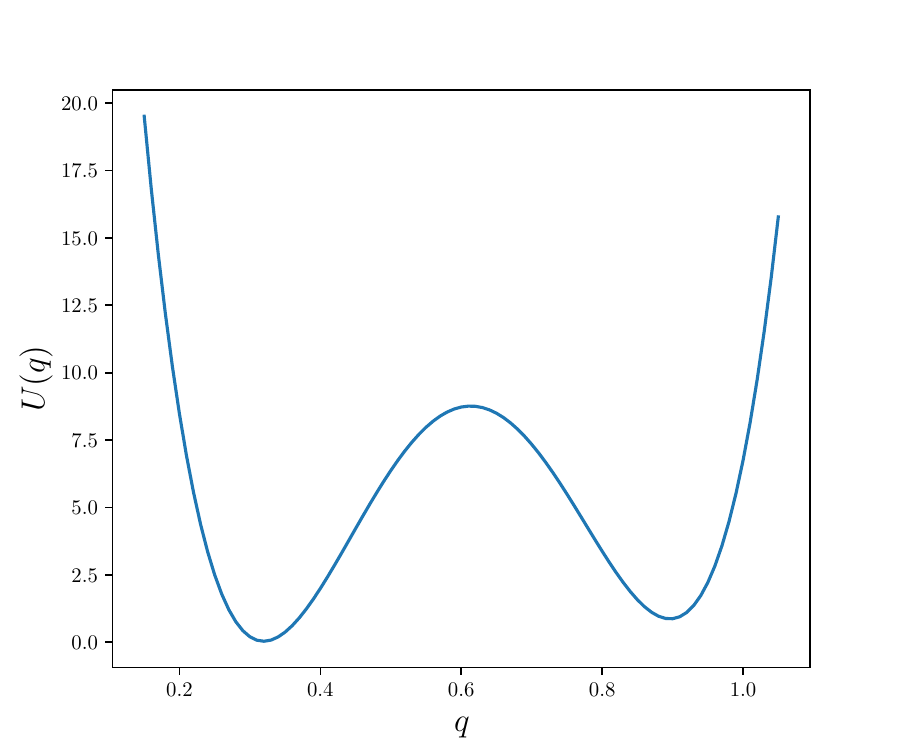}
     \caption{Free energy $U(q)$} %
     \label{fig:free_energy}
     \end{subfigure}
     \centering
    \centering
     \begin{subfigure}[b]{0.40\textwidth}
         \centering
     \includegraphics[width=\textwidth]{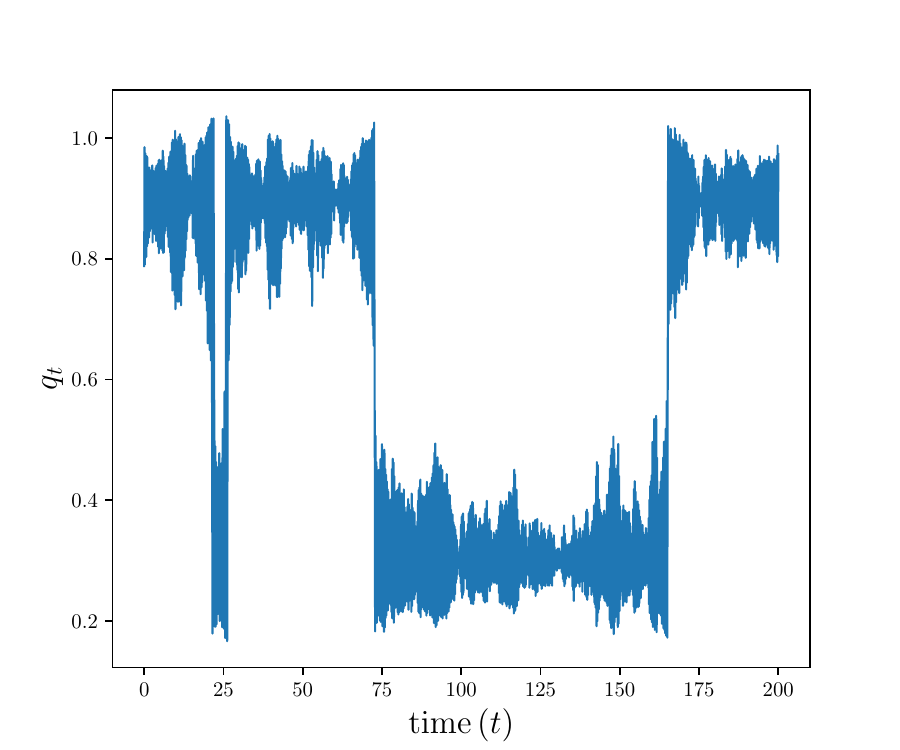}
     \caption{Trajectory of $q_t$} %
     \label{fig:sample_path}
     \end{subfigure}
     \centering
    \caption{
    Left panel (\ref{fig:free_energy}):
    The free energy used in the experiment. 
    Right panel (\ref{fig:sample_path}): A trajectory of the observable coordinate $q_t$ from the QGLE (\ref{eq:qgle_second}).
         }
    \label{fig:gle_model}
\end{figure} 
\subsubsection{Multivariate Extended State}
\label{sec:num_qgle_2}
We consider a QGLE with one-dimensional coordinates and multivariate extended variable, motivated from the work of  \citep{ay:21} that studies protein-folding kinetics via a Quasi-Markovian GLE (\ref{eq:qgle-II}) and showcases that a QGLE accurately reproduces simulations of molecular dynamics (MD) that involve memory effects in the friction. In their investigation, a one-dimensional reaction coordinate, $q_t$, given as the sum of the separations between native contacts, is modelled via the following QGLE:
\begin{align}
\begin{aligned} \label{eq:qgle_second}
d q_t & = \tfrac{1}{m} \times p_t \, dt; \\[0.1cm] 
d p_t & = - {U'} (q_t) dt + \sum_{l = 1}^L s_{l, t} dt; \\[0.1cm] 
d s_{l, t} & = - \tfrac{1}{\tau_l} \times s_{l, t} \,  dt  -  \tfrac{c_l}{\tau_l} \times p_t \, dt + \tfrac{\sqrt{2 \beta^{-1} c_l}}{\tau_l} \, d B_{l, t},  \quad  s_{l, 0} \sim \mathscr{N} (0, \beta^{-1}), \quad  \, 1 \le l \le L, 
\end{aligned}
\end{align}
where $m, \beta > 0$ denote the mass and the inverse thermal energy respectively, $\{ c_l, \tau_l\}_{1 \le l \le L}$ are the unknown parameters taking positive values, for $L\ge 1$, and  $U : \mathbb{R} \to \mathbb{R}$, the folding free energy landscape for proteins, is specified as $q \mapsto U(q) = - \beta^{-1} \log \nu (q)$ with $\nu (\cdot)$ being the equilibrium probability density function. QGLE (\ref{eq:qgle_second}) corresponds to the non-Markovian GLE (\ref{eq:gle}) with the memory kernel given as a so-called \emph{Prony series}:
\begin{align} \label{eq:memory_kernel}
K(t) = \sum_{l = 1}^L \frac{c_l}{\tau_l} 
\times \exp \Bigl(- \frac{t}{\tau_l} \Bigr), \quad t \geq 0.  
\end{align}
\cite{ay:21} constructed QGLE (\ref{eq:qgle_second}) with $L = 5$ by determining the parameters via a least squares method so that the memory kernel (\ref{eq:memory_kernel}) fits the one extracted numerically from the observed time-series of $q$. 

\begin{figure}
    \centering
     \includegraphics[width=\textwidth]{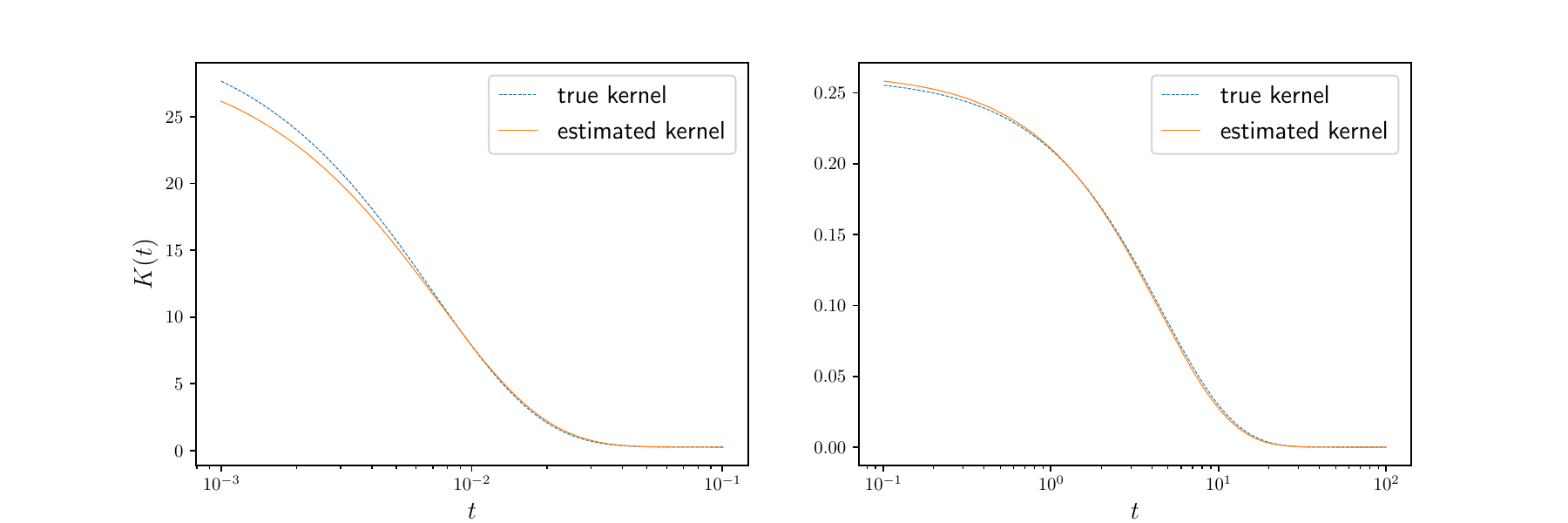}
    \caption{
    Memory kernel (\ref{eq:memory_kernel}) computed with the true parameter (true kernel) and with the mean value of MLEs (estimated kernel). 
    }
    \label{fig:memory_kernel}
\end{figure}

In our experiment, we estimate the unknown parameters by maximising the marginal likelihood given the partial observations $q_{0:n}$. 
{Since the mass ($m$) of a particle can typically be measured, we assume $m$ to be known and set $m = 1$, in agreement with numerical experiments in the literature \citep{dit:19, poke:09}}. We set $L = 2$, and specify the free energy function as $q \mapsto U (q) = a (q - q_{\mathrm{min}})^2 (q - q_{\mrm{max}})^2 + b q^3$ with constants $(q_{\mathrm{min}}, q_{\mathrm{max}}, a, b) = (0.30, 0.90, 1200, 0.001)$. We set $\beta^{-1} = 2.949$ and select the true parameters as $\theta^{\dagger} = (c_1^{\dagger}, \tau_1^{\dagger}, c_2^{\dagger}, \tau_2^{\dagger}) = (0.22, 0.007, 1.2, 4.6)$, as such a choice closely reproduces the shape of the memory kernel estimated in \citep{ay:21}. We generate $50$ independent trajectories of $q$ on the time inverval $[0,1500]$ by applying scheme (\ref{eq:lg_hypo_II}) to QGLE (\ref{eq:qgle_second}) with step-size $10^{-4}$. We discard the observations up to time $500$ and sub-sample the datasets $q_{0:n}$ in equilibrium, so that  $n = 10^6$, $\Delta_n = 10^{-3}$, $T = 1000$. 
In Figure~\ref{fig:gle_model}, we plot the shape of the free energy $U$ and one trajectory of the component $q$ from the QGLE (\ref{eq:qgle_second}) in the chosen setting. Notice that QGLE (\ref{eq:qgle_second}) is a conditionally Gaussian non-linear system (given the component~$q$),  thus
upon adoption of the locally Gaussian discretisation (\ref{eq:lg_hypo_II}),
the marginal likelihood can be calculated via the Kalman filter shown in Section \ref{sec:kalman_filter}. We use the Nelder-Mead method to optimise the marginal likelihoods with the initial value $\theta_0 = (0.1, 0.01, 1.0, 10.0)$. 
{Table \ref{table:protein_folding} summarises the mean and standard deviation of relative error of estimates 
$\hat{\theta_n} = (\hat{c}_{1,n}, \hat{\tau}_{1,n}, \hat{c}_{2,n}, \hat{\tau}_{2,n})$. The estimation results are overall acceptable though the estimator $\hat{\tau}_{1, n}$ is slightly biased and centred around $0.0073$ -- we have found empirically that such bias disappears upon chosing a smaller step-size $\Delta_n$, e.g.~$\Delta_n = 2 \times 10^{-4}$, so that $n \Delta_n^2$ becomes closer to $0$. 
}
We typically observe that the standard derviation of estimates for $(c_1, \tau_1)$ is much smaller than that for $(c_2, \tau_2)$. In Figure \ref{fig:memory_kernel}, we plot the memory kernel (\ref{eq:memory_kernel}) with the parameters set equal to the mean value of the $50$ MLEs to observe the level of agreement between the estimated memory kernel and the reference kernel, the latter computed with the true parameter values. The relative absolute errors between the true and estimated memory kernels are within $0.1$ across periods $t \in [0.001, 10]$.  
\begin{table}
\caption{
    {Parameter estimation of the QGLE (\ref{eq:qgle_second}). Mean and standard deviation (in brackets) of 
    $(\hat{\theta}_n - \trueparam)/\trueparam$ 
    from $50$ trajectories of $q_{0:n}$. The value of true parameter is $\theta^{\dagger} = (c_1^{\dagger}, \tau_1^{\dagger}, c_2^{\dagger}, \tau_2^{\dagger}) = (0.22, 0.007, 1.2, 4.6)$.
    }
    } 
    \label{table:protein_folding}
\centering 
\begin{tabular}{c c c c} 
\toprule 
\\[-10pt] 
$ \bigl( \hat{c}_{1,n}-c_1^\dagger \bigr)/c_1^\dagger$  
& $\bigl( \hat{\tau}_{1,n}-\tau_1^\dagger \bigr)/\tau_1^\dagger$  
& $\bigl( \hat{c}_{2,n}-c_2^\dagger \bigr)/c_2^\dagger$  
& $\bigl( \hat{\tau}_{2,n}-\tau_2^\dagger \bigr)/\tau_2^\dagger$   
\\[0.2cm]
\hline 
\\[-0.3cm]
-0.0099	(0.0090)
& 0.0481 (0.0046)
& -0.0364 (0.2441)
& -0.0423 (0.1795)
\\ 
\bottomrule   
\end{tabular}
\end{table} 
%
%
\section{Conclusions and Future Directions}
\label{sec:end}
We have studied parameter inference procedures for the highly degenerate class of SDEs 
that includes a wide range of practical models, e.g., quasi-Markovian generalised equations (QGLEs), 
epidemiological models with time-varying parameters \citep{sir:22, dur:13},  non-linear continuous-time autoregressive models \citep{ts:00} and the classical Lorentz system upon consideration of noise effects \citep{co:21}. We have introduced the locally Gaussian time-discretisation scheme (\ref{eq:lg_hypo_II}) and provided analytical/numerical results showcasing that parameter estimation based upon such scheme sidesteps biases that would arise under alternative schemes. 
The approach followed in this work for establishing our results for class  (\ref{eq:hypo-II}) are expected to also 
guide extensions to more general classes of degenerate diffusions, for which iterated Lie brackets of \emph{any} order, i.e., $[\tilde{V}_0, [\tilde{V}_0,\ldots, [\tilde{V}_0, V_k]]]$, $1 \le k \le d$, are required for Hormander's condition to hold. Here, we draw upon the understanding obtained via the study of classes (\ref{eq:hypo-I}) and (\ref{eq:hypo-II}) to summarise key arguments for carrying out unbiased parameter estimation for general hypo-elliptic systems.  
First, in a partial observation regime, use of a degenerate discretisation (e.g.~Euler-Maruyama) or equivalently of finite-differences to impute latent components will induce bias at estimates of diffusion coefficient parameters (recall the case study in Section \ref{sec:case_2}). {Also, note that finite-differences can be used only when the smooth component is determined as $d X_{S, t} = X_{R, t} dt$. Thus, a natural approach for parametric inference of degenerate diffusions is to develop a non-degenerate conditionally Gaussian scheme for the full coordinates}, with the Gaussian noise obtained via high-order stochastic Taylor expansion of the drift functions. 
A lot of care should be given at the deterministic terms of the expansion to be included into the scheme, to avoid emergence of biases in estimates of drift parameters (recall the analytical study in Section \ref{sec:case_drift} and the numerical results in Section \ref{sec:num_case_study}).   
We summarise below the above-designated roadmap for the construction of `correct' time-discretisation schemes for general classes of hypo-elliptic diffusions. 
%
%
\begin{itemize}
\item[Step 1.]  
For the rough component, $X_{R}$, the Euler-Maruyama scheme is applied. 
%
\item[Step 2.]  For the smooth coordinates in the model, one recursively applies stochastic Taylor expansion to drift functions so that Gaussian variates, in the form of iterated integrals involving Brownian motions, e.g.~of the form $\textstyle{\int_{t_i}^{t_{i+1}} B_s ds}$, $\textstyle{\int_{t_i}^{t_{i+1}} \int_{t_i}^{u} B_s ds du}$, appear in all smooth coordinates. This process is completed once the covariance-variance of the Gaussian approximation is positive definite. 
\item[Step 3.] For a smooth component containing Gaussian noise of size $\mathcal{O} ( \Delta_n^{(2k - 1)/2} )$, $k \ge 2$, the scheme should include all deterministic terms from the stochastic  Taylor expansion up to size  $\mathcal{O} (\Delta_n^k)$. 
\end{itemize}
Indicatively, Table \ref{table:lg_II} summarises the size of determistic and noisy parts of the locally Gaussian scheme (\ref{eq:lg_hypo_II}) for class (\ref{eq:hypo-II}). 
%

Our work in this paper leads to further research in several directions. 
In the CLT of the main analytical result for the parameter estimator (Theorem \ref{thm:clt}), the step-size $\Delta_n$ is required to satisfy $\Delta_n = o (n^{-1/2})$. An open problem for hypo-elliptic diffusions is the construction of estimators giving a CLT under a weaker condition $\Delta_n = o (n^{-1/p})$, $p \ge 3$. 
We expect that such a general estimator for degenerate diffusion models can be produced, with accompanying theory then following the strategy used in our proofs in this work, as we mentioned in Remark \ref{rem:pf_main}, that is, the proof of consistency does not require a {condition of $\Delta_n = o(n^{-1/2})$.}  In a different direction, the effectiveness of the developed locally Gaussian scheme is yet to be studied under a low-frequency observation setting, i.e.~with the step-size $\Delta$ assumed fixed and not small enough, in which case a number, say $M$, of inner sub-steps are introduced by the user. Under such a setting, the discretisation error of the true (intractable) density over the period of size $\Delta$ typically diminishes as $M$ increases. In the case of elliptic diffusions, explicit rates of convergence to zero are provided in  \cite{go:08, iguchi:21-2}. Finally, 
in this work, in the practical scenario of partial observations, we have investigated the behaviour of discretisation schemes via case studies and numerical experiments. Analytical theory would be quite instructive in this setting. Techniques used in the context of hidden Markov models (see e.g.~\cite{douc:14}) are expected to be valuable in such a pursuit.  
\begin{table}
\centering
\caption{Size (in $\Delta_n$) of the terms appearing in the locally Gaussian scheme (\ref{eq:lg_hypo_II}).}
\label{table:lg_II}
\begin{tabular}{lcc}
 \toprule
  {Component}
 & {Gaussian part}
 & {Deterministic part}
  \\ 
 \midrule
 \\[-0.4cm] 
 $\bar{X}^{(\mrm{II})}_{S_1, {i+1}}$ 
 &   $\mathcal{O} (\Delta_n^{5/2}), \quad 
 \Bigl( {\textstyle \int_{t_i}^{t_{i+1}} \int_{t_i}^u B_{s} ds du} \Bigr)$      &   $\mathcal{O} (\Delta_n^3)$      
 \\[0.4cm]
 $ \bar{X}^{(\mrm{II})}_{S_2, {i+1}}$  
 & $ \mathcal{O} (\Delta_n^{3/2}), \quad 
 \Bigl( {\textstyle \int_{t_i}^{t_{i+1}} B_{s} ds} \Bigr)$                
 & $\mathcal{O} (\Delta_n^2)$       
 \\[0.4cm]
 $\bar{X}^{(\mrm{II})}_{R, {i+1}}$                           
 & $\mathcal{O} (\Delta_n^{1/2}), \quad \bigl( B_{t_{i+1}} - B_{t_{i}} \bigr)$
 & $\mathcal{O} (\Delta_n) $ 
 \\[0.2cm] 
\bottomrule
\end{tabular}
\end{table}
\section*{Funding}
\noindent Yuga Iguchi acknowledges support from the Additional Funding Programme for Mathematical Sciences, delivered by EPSRC (EP/V521917/1) and the Heilbronn Institute for Mathematical Research.
\section*{Acknowledgements}
We thank the associate editor and two anonymous referees for their helpful comments, which improved the quality of the article. 
\begin{appendices}
\section*{Appendix}
\section{Preliminaries} \label{appendix:aux}
In Section \ref{app:notation} we present some notation  used in the Appendix. In Section \ref{app:aux} we introduce six auxiliary results needed in the proof of our main theorems (Theorems~\ref{thm:consistency},  \ref{thm:clt}) in Section~\ref{sec:main}.  
\subsection{Notation}
\label{app:notation}
For $0 = t_0 < \cdots < t_n$, with equi-distant step-size $\Delta_n$, we write $\sample{i}$ for the observation at time $t_i$ of the solution of the hypo-elliptic SDE (\ref{eq:hypo-II}) under the true parameter value~$\trueparam$,  defined upon the filtered probability space $(\Omega, \mathcal{F}, \{\mathcal{F}_t\}_{t \geq 0}, \mathbb{P})$. We denote by $\truedist$ the invariant distribution of process (\ref{eq:hypo-II}) under $\trueparam$. In agreement with the structure of class (\ref{eq:hypo-II}), we often represent $x \in \mathbb{R}^N$ and $\theta \in \Theta \subset \mathbb{R}^{N_{\theta}}$ as
\begin{gather*}
x = (x_{S_1}, x_{S_2}, x_R) \in \mathbb{R}^{N_{S_1}} \times \mathbb{R}^{N_{S_2}} \times \mathbb{R}^{N_R}, \quad  
x_S \equiv (x_{S_1}, x_{S_2}); \\[0.2cm]
\theta = (\beta_{S_1}, \beta_{S_2}, \beta_R, \sigma) \in 
\Theta_{{\beta_{S_1}}} \times \Theta_{{\beta_{S_2}}} \times \Theta_{{\beta_{R}}} \times \Theta_{{\sigma}}, \quad 
\beta_S \equiv (\beta_{S_1}, \beta_{S_2}),
\end{gather*}
For $\varphi (\cdot , \theta) : \mathbb{R}^N \to \mathbb{R}$, $\theta \in \Theta$, bounded up to second derivatives, we define the differential operators $\mathcal{L}$ and $\mathcal{L}_j, \; 1 \le  j  \le d$: 
\begin{gather*}
\mathcal{L} \varphi (x ,  \theta) 
= \sum_{i=1}^N V_0^i (x, \theta) \frac{\partial \varphi} {\partial x_i}(x, \theta)  
+ \tfrac{1}{2}  \sum_{i_1, i_2 = 1}^N \sum_{k=1}^d  V_k^{i_1} (x, \theta) V_k^{i_2} (x, \theta)  \frac{\partial^2 \varphi }{\partial x_{i_1} \partial x_{i_2} } (x,\theta);  \\
\mathcal{L}_j \varphi (x , \theta)
= \sum_{i=1}^N V_j^i (x, \theta) \frac{\partial \varphi}{\partial x_i}(x , \theta), \quad 1 \le j \le d.  
\end{gather*} 
%
Application of the above differential operators is extended to vector-valued functions in the apparent way, via separate consideration of each scalar component. We recall some notation used in the definition of the contrast function $\ell_n (\theta)$ in (\ref{eq:contrast}).
We have that 
\begin{align*} 
\mu (\Delta, x, \theta) = \bigl[ 
\mu_{S_1} (\Delta, x, \theta)^\top , \, 
\mu_{S_2}  (\Delta, x, \theta)^\top , \, 
\mu_{R} (\Delta, x, \theta)^\top 
\bigr]^\top,  
\end{align*}
where 
\begin{align*} 
\begin{bmatrix}
{\mu}_{S_1}  (\Delta, x, \theta) \\[0.2cm]
{\mu}_{S_2}  (\Delta, x, \theta)  \\[0.2cm]
{\mu}_{R}  (\Delta, x, \theta) 
\end{bmatrix} 
=
\begin{bmatrix}
x_{S_1} + V_{S_1,0} (x_S, \beta_{S_1}) \Delta 
+ \mathcal{L} V_{S_1, 0} (x, \theta) \tfrac{\Delta^2}{2}
+ \mathcal{L}^2 V_{S_1, 0} (x, \theta) \tfrac{\Delta^3}{6}   
\\[0.3cm]
x_{S_2} + V_{S_2,0} (x, \beta_{S_2}) \Delta 
+ \mathcal{L} V_{S_2, 0} (x, \theta) \tfrac{\Delta^2}{2}  \\[0.3cm] 
x_R + V_{R, 0} (x, \beta_R) \Delta
\end{bmatrix}. 
\end{align*} 
When $\Delta = 1$, we simply write 
\begin{align*}
\mu (x, \theta) \equiv \mu (1, x, \theta). 
\end{align*}
For $x = (x_{S_1}, x_{S_2}, x_R) \in \mathbb{R}^N \equiv \mathbb{R}^{N_{S_1}} \times \mathbb{R}^{N_{S_2}} \times \mathbb{R}^{N_{R}}, \, y = (y_{S_1}, y_{S_2}, y_R) \in \mathbb{R}^N$ , $\Delta > 0$ and $\theta \in \Theta$, we define
\begin{align} \label{eq:m} 
m  (\Delta, x, y , \theta)  
= \left[
\,  
\frac{y_{S_1}^\top - {\mu}_{S_1} (\Delta, x; \theta)^\top}{\sqrt{\Delta^{5}}}, \; 
\frac{y_{S_2}^\top - {\mu}_{S_2} (\Delta, x; \theta)^\top}{\sqrt{\Delta^{3}}}, \;  
\frac{y_{R}^\top - {\mu}_{R} (\Delta, x; \theta)^\top }{\sqrt{\Delta}} \,
\right]^\top. 
\end{align}
%
%
%
We write, for $1 \le i \le n$, 
\begin{gather}  \label{eq:m_simple}
  m_{i} (\Delta, \theta)  
\equiv m (\Delta, \sample{i-1}, \sample{i}, \theta). 
\end{gather}
We use $\Sigma (\Delta, x, \theta)$ to represent the covariance of one step of the local Gaussian scheme (\ref{eq:lg_hypo_II}) for the hypo-elliptic SDE (\ref{eq:hypo-II}), given step-size $\Delta>0$, initial point $x\in\mathbb{R}^N$ and parameter~$\theta$. 
We often write 
\begin{align*}
\Sigma (x, \theta) \equiv \Sigma (1, x, \theta).
\end{align*}
We express the inverse of $\Sigma (x, \theta)$
as: 
\begin{align} 
   \Sigma^{-1} (x, \theta) 
   = 
   \Lambda (x, \theta)
   & =
   \begin{bmatrix}  \label{eq:inv_Sigma_2} 
    \Lambda_{S_1 S_1} (x, \theta) 
    & \Lambda_{S_1 S_2} (x, \theta)
    & \Lambda_{S_1 R} (x, \theta)  \\[0.2cm] 
    \Lambda_{S_2 S_1} (x, \theta) 
    & \Lambda_{S_2 S_2} (x, \theta)
    & \Lambda_{S_2 R} (x, \theta)  \\[0.2cm]  
    \Lambda_{R S_1} (x, \theta) 
    & \Lambda_{R S_2} (x, \theta)
    & \Lambda_{R R} (x, \theta) 
   \end{bmatrix},    
\end{align}  
where each block matrix is specified as
\begin{align*}
\Lambda_{\iota_1 \iota_2} (x, \theta) \in \mathbb{R}^{N_{\iota_1} \times N_{\iota_2}}, 
\quad \iota_1, \iota_2  \in \{S_1, S_2, R\}. 
\end{align*}
%
We emphasise here that $\Sigma (x, \theta)$ and its inverse $\Lambda (x, \theta)$ depend on $x$ and $(\beta_{S}, \sigma)$ but not on the drift parameter $\beta_R$ in the rough component, and this is critical in the proof of consistency of $\hat{\beta}_{R, n}$. Thus, we sometimes write 
$\Sigma (x, (\beta_S, \sigma))$ and $\Lambda (x, (\beta_S, \sigma ))$ to highlight the parameter dependency.
%
We recall the definition of the matrices 
$a_R(x, \sigma) \in \mathbb{R}^{N_R \times N_R}$, \, 
$a_{S_1}(x, \theta) \in \mathbb{R}^{N_{S_2} \times N_{S_2}}$, \, 
$a_{S_2} (x, \theta) \in \mathbb{R}^{N_{S_1} \times N_{S_1}}$ as: 
\begin{align*}
 a_R (x , \sigma) & =  
 \sum_{k = 1}^d  V_{R,k} (x , \sigma)  V_{R, k} (x , \sigma)^\top, \quad 
a_{S_2} (x, \theta) =
\sum_{k = 1}^d 
\Bigl( \partial_{x_R}^\top V_{S_2, 0} (x, \beta_{S_2})  V_{R,k} (x , \sigma) 
\Bigr) 
\Bigl( \partial_{x_R}^\top V_{S_2, 0} (x, \beta_{S_2}) V_{R,k} (x , \sigma) \Bigr)^\top; \\[0.2cm] 
a_{S_1} (x, \theta) & = 
\sum_{k = 1}^d \, 
\Bigl( 
\partial_{x_{S_2}}^\top V_{S_1, 0} (x_S, \beta_{S_1})
\partial_{x_R}^\top V_{S_2, 0} (x, \beta_{S_2}) V_{R,k} (x , \sigma)  
\Bigr) 
\Bigl( 
\partial_{x_{S_2}}^\top V_{S_1, 0} (x_S, \beta_{S_1})
\partial_{x_R}^\top V_{S_2, 0} (x, \beta_{S_2}) V_{R,k} (x , \sigma)  
\Bigr)^\top. 
\end{align*} 
We define the mappings 
\begin{align*}
\eta_{S_1}:\mathbb{R}^N \times \Theta_{\beta_{S_1}} \to \mathbb{R}^{N_{S_1}}, \quad 
\eta_{S_2}:\mathbb{R}^N \times \Theta_{\beta_{S_2}} \to \mathbb{R}^{N_{S_2}}, \quad   
\eta_{R}:\mathbb{R}^N \times \Theta_{\beta_{R}}  \to \mathbb{R}^{N_{R}} 
\end{align*}
as 
\begin{gather*}
\eta_{S_1} (x, \beta_{S_1})  = V_{S_1, 0} (x_S, \truebeta_{S_1}) - V_{S_1, 0} (x_S, \beta_{S_1}), \quad  
\eta_{S_2} (x, \beta_{S_2}) = V_{S_2, 0} (x, \truebeta_{S_2}) - V_{S_2, 0} (x, \beta_{S_2}); \\[0.2cm] 
\eta_{R} (x, \beta_{R}) = V_{R, 0} (x, \truebeta_{R}) - V_{R, 0} (x, \beta_{R}).   
\end{gather*}
We write, with a slight abuse of notation, for $0 \le i \le n$,
\begin{align*}
\eta_{S_1, i} (\Delta, \beta_{S_1})
\equiv \tfrac{\eta_{S_1} (\sample{i}, \, \beta_{S_1})}{\Delta}, 
\quad \eta_{S_2, i} (\Delta, \beta_{S_2})
 \equiv \tfrac{ \eta_{S_2} (\sample{i} , {\beta}_{S_2})}{\Delta}. 
\end{align*}
We denote by $\mathcal{S}$ the space of functions $f : [0, \infty) \times \mathbb{R}^N \times \Theta \to \mathbb{R}$ so that there are  constants $C, q > 0$ such that 
$| f (\Delta, x, \theta) |  \le C \Delta\, ( 1 + |x|^q)$ for any $(\Delta, x, \theta) \in [0, \infty) \times \mathbb{R}^N \times \Theta$. For an $M_1 \times M_2$ matrix $A$, with $M_1,M_2\ge 1$, we write each matrix entry as $[A]_{ij}$ for $1 \le i \le M_1, \, 1 \le j \le M_2$.  
An expectation under the probability law $\mathbb{P}_\theta$ is written as $\mathbb{E}_\theta$. 
We write $\textstyle \partial_{u}  = \big[ {\partial }/{\partial u_{1}}, \ldots, {\partial }/{\partial u_{n} } \big]^{\top}, \, \partial^{2}_{u}   = \partial_{u} \partial_{u}^{\top} \equiv \big({\partial^{2}}/{\partial u_{i} \partial u_{j} } \big)_{i,j=1}^{n}$ for the standard differential operators acting upon maps $\mathbb{R}^{n}\to \mathbb{R}$, $n\ge 1$. 
We also write $\partial^u_\alpha = {\partial^l}/{\partial u_{\alpha_1} \cdots \partial u_{\alpha_l}}$ for a multi-index $\alpha \in \{1, \ldots, n \}^l$, $l \in \mathbb{N}$. For a function $g : \mathbb{R}^n \to \mathbb{R}^m, \, n, m \in \mathbb{N}$, we write: 
\begin{gather*}
\partial_u g (u)^\top = \bigl[ \tfrac{\partial}{\partial u_i} g^j (u) \bigr]_{\substack{1 \le i \le n, \, 1 \le j \le m}}, \quad 
\partial_u^\top g (u) = \bigl( \partial_u g (u)^\top  \bigr)^\top;  \\[0.2cm]  
\partial^u_\alpha g (u)^\top = \bigl[ \partial^u_\alpha g^1 (u), \ldots, \partial^u_\alpha g^m (u) \bigr], \, \quad \partial^u_\alpha g (u) = \bigl( \partial^u_\alpha g (u)^\top \bigr)^\top.  
\end{gather*}
%
\subsection{Auxiliary Results}
\label{app:aux}
\begin{lemma} \label{lemma:aux_1}
Let $Y_{t_i}$, $U$ be random variables, with $Y_{t_i}$ being $\mathcal{F}_{t_i}$-measurable. 
If
\begin{align*}
  \sum_{i=1}^n \trueE [\,Y_{t_{i}}\,|\,\filtration{i-1}\,]  \probconv  U,  \quad 
  \sum_{i=1}^n \trueE\big[\,(Y_{t_{i}})^2\,|\,\filtration{i-1}\,\big]  \probconv  0,
\end{align*} 
then $\textstyle{\sum_{i=1}^n Y_{t_{i}} \probconv U}$.
\end{lemma} 
\begin{proof}[Proof.]
See Lemma 9 in \cite{genon:93}. 
\end{proof} 
\begin{lemma} \label{lemma:ergodic_thm}
Let $f : \mathbb{R}^N \times \Theta \to \mathbb{R}$ be differentiable w.r.t.~$(x, \theta) \in \mathbb{R}^N \times \Theta$ with derivatives of polynomial growth in $x$ uniformly in $\theta$. Under conditions (\ref{assump:coeff})--(\ref{assump:finite_moment}), it holds that, if \limit, then
\begin{align*}
   \tfrac{1}{n} \sum_{i=1}^n f (\sample{i-1}, \theta) \probconv 
   \int f (x, \theta) \truedist (dx),
\end{align*}
uniformly in $\theta \in \Theta$. 
\end{lemma}
\begin{proof}[Proof.]
This is a multivariate version of Lemma 8 in \cite{kess:97}, so we omit the proof.  
\end{proof}
\begin{lemma} \label{lemma:canonical_conv}
Let $1 \leq j_1, j_2 \leq N$ and assume that $f : \mathbb{R}^N \times \Theta \to \mathbb{R}$ is as in Lemma \ref{lemma:ergodic_thm}. Under conditions (\ref{assump:coeff})--(\ref{assump:finite_moment}), it holds that, if \limit, then
\begin{gather}
\tfrac{1}{n} \sum_{i=1}^n 
f (\sample{i-1}, \theta) \, 
m^{j_1}_i (\Delta_n, \trueparam)  \, 
m^{j_2}_i (\Delta_n, \trueparam) \probconv \int  f (x, \theta) \bigl[ {\Sigma} (x, \trueparam) \bigr]_{j_1 j_2} \truedist (dx);  \label{eq:canonical_conv1} \\
\tfrac{1}{n  \sqrt{\Delta_n}} \sum_{i=1}^n   f (\sample{i-1}, \theta)  m^{j_1}_i  (\Delta_n, \trueparam)  \probconv 0, 
\label{eq:canonical_conv2} 
\end{gather}
uniformly in $\theta \in \Theta$. 
\end{lemma}
\begin{proof}[Proof.]
Notice that for any $1 \le j_1, j_2 \le N$ and $0 \le i \le n$, 
\begin{align} 
& \mathbb{E}_{\trueparam}
\bigl[m^{j_1}_{i} (\Delta_n, \trueparam) 
\, | \,  \mathcal{F}_{t_{i-1}} \bigr]
= R_{j_1} (\sqrt{\Delta_n^3}, \sample{i-1}, \trueparam); 
\label{eq:mean_expectation} \\[0.3cm]
& 
\mathbb{E}_{\trueparam} 
\bigl[ 
m^{j_1}_i (\Delta_n, \trueparam) \, 
m^{j_2}_i (\Delta_n, \trueparam)  
\, |  \, \mathcal{F}_{t_{i-1}} \bigr]  
= \bigl[ \Sigma (\sample{i-1}, \trueparam) \bigr]_{j_1 j_2} + R_{j_1 j_2} (\Delta_n, \sample{i-1}, \trueparam),  \label{eq:squared_mean_expectation} 
\end{align} 
where $R_{j_1 j_2}, \, R_{j_1} \in \mathcal{S}$. Applying Lemmas \ref{lemma:aux_1}, \ref{lemma:ergodic_thm} with formulae (\ref{eq:mean_expectation}) and (\ref{eq:squared_mean_expectation}), we obtain (\ref{eq:canonical_conv1}) and (\ref{eq:canonical_conv2}) in the same way as in the proof of Lemma 12 in \cite{iguchi:22}. 
\end{proof}  
\begin{lemma} \label{lemma:matrix}
Assume that condition \ref{assump:hypo}-II holds. We have that for any $(x, \theta) \in \mathbb{R}^N \times \Theta$: \begin{align*}
    \Lambda_{S_1 S_1} (x, \theta)  \, 
    \partial_{x_{S_2}}^\top V_{S_1, 0} (x_S, \beta_{S_1}) 
    = -2 \Lambda_{S_1 S_2} (x, \theta).
\end{align*}    
\end{lemma} 
\begin{proof}
We write $\Sigma (x, \theta)$, $(x, \theta) \in \mathbb{R}^N \times \Theta$, in the form of the block matrix: 
\begin{align} \label{eq:Sigma}
\Sigma (x, \theta)
= 
\begin{bmatrix}
\Sigma_{S_1 S_1} (x, \theta)
& \widetilde{\Sigma} (x, \theta) \\[0.2cm]
\widetilde{\Sigma} (x, \theta)^\top
& \hat{\Sigma} (x, \theta)
\end{bmatrix},
\end{align}
where we have set: 
\begin{align} \label{eq:Sigma_blocks}
\widetilde{\Sigma} (x, \theta) 
= \Bigl[  
    \Sigma_{S_1 S_2} (x, \theta), \, 
    \Sigma_{S_1 R} (x, \theta) 
    \Bigr],
\quad 
\hat{\Sigma} (x, \theta) =
\begin{bmatrix}
\Sigma_{S_2 S_2} (x, \theta) & \Sigma_{S_2 R} (x, \theta)  \\[0.2cm]
\Sigma_{R S_2} (x, \theta)  & \Sigma_{RR} (x, \theta)
\end{bmatrix}. 
\end{align}
Notice that under condition \ref{assump:hypo}--II, matrix $\hat{\Sigma} (x, \theta)$ is invertible for any $(x, \theta) \in 
\mathbb{R}^N \times \Theta$. We write the inverse of $\hat{\Sigma}(x, \theta)$ as: 
\begin{align*}
\hat{\Sigma}^{-1} (x, \theta)
=
\hat{\Lambda} (x, \theta) 
= 
\begin{bmatrix}
 \hat{\Lambda}_{S_2 S_2} (x, \theta)
 & \hat{\Lambda}_{S_2 R} (x, \theta) \\[0.2cm]
 \hat{\Lambda}_{R S_2} (x, \theta)
 & \hat{\Lambda}_{R R} (x, \theta)
\end{bmatrix}.  
\end{align*}
 %
 Recall the notation for the inverse of $\Sigma (x, \theta)$ in (\ref{eq:inv_Sigma_2}). Using the inverse formula for a block matrix, we obtain: 
 \begin{align} \label{eq:inv_1}
     \Lambda_{S_1 S_2} (x, \theta)
     = - \Lambda_{S_1 S_1} (x, \theta) \Xi (x, \theta), 
 \end{align}
 where we have set 
\begin{align} \label{eq:Xi}
\Xi (x, \theta) =  \Sigma_{S_1 S_2} (x, \theta) \hat{\Lambda}_{S_2 S_2} (x, \theta) 
 + \Sigma_{S_1 R} (x, \theta) 
 \hat{\Lambda}_{R S_2} (x, \theta). 
\end{align}
%
From the block matrix representation of $\Sigma (\Delta, x, \theta)$ in (\ref{eq:covariance_lg}), we obtain 
\begin{gather*}
\Sigma_{S_1 S_2} (x, \theta) = \tfrac{3}{8} \partial_{x_{S_2}}^\top V_{S_1, 0} (x_S, \beta_{S_1}) \Sigma_{S_2 S_2} (x, \theta), \quad 
\Sigma_{S_1 R} (x, \theta)  = \tfrac{1}{3} \partial_{x_{S_2}}^\top V_{S_1, 0} (x_S, \beta_{S_1}) \Sigma_{S_2 R} (x, \theta).  
\end{gather*}
We then have    
 \begin{align}  
 \Xi (x, \theta)  
 & =  \tfrac{3}{8} \partial_{x_{S_2}}^\top V_{S_1, 0} (x_S, \beta_{S_1}) \Sigma_{S_2 S_2} (x, \theta) \hat{\Lambda}_{S_2 S_2} (x, \theta) + \tfrac{1}{3} \partial_{x_{S_2}}^\top V_{S_1, 0} (x_S, \beta_{S_1}) \Sigma_{S_2 R} (x, \theta) \hat{\Lambda}_{R S_2} (x, \theta) \nonumber \\[0.2cm]
 & = \tfrac{1}{24} \partial_{x_{S_2}}^\top V_{S_1, 0} (x_S, \beta_{S_1}) \Sigma_{S_2 S_2} (x, \theta) \hat{\Lambda}_{S_2 S_2} (x, \theta) 
 + \tfrac{1}{3} \partial_{x_{S_2}}^\top V_{S_1, 0} (x_S, \beta_{S_1}) \nonumber \\[0.2cm]  
 & = \tfrac{1}{2}  \partial_{x_{S_2}}^\top V_{S_1, 0} (x_S, \beta_{S_1}).  \label{eq:inv_2}
 \end{align}
 In the above calculation we have used: 
 \begin{gather*}
{\Sigma}_{S_2 S_2} (x, \theta) 
\hat{\Lambda}_{S_2 S_2} (x, \theta) 
+ {\Sigma}_{S_2 R} (x, \theta) 
\hat{\Lambda}_{R S_2} (x, \theta)  
= I_{N_{S_2} \times N_{S_2}}; \\[0.1cm]
\hat{\Lambda}_{S_2 S_2} (x, \theta) 
= \bigl({\Sigma}_{S_2 S_2} (x, \theta) 
- \Sigma_{S_2 R} (x, \theta) 
\Sigma^{-1}_{RR} (x, \theta) 
\Sigma_{R S_2} (x, \theta)  \bigr)^{-1} = 4 {\Sigma}_{S_2 S_2}^{-1} (x, \theta), 
 \end{gather*}
 where  matrix $\Sigma_{S_2 S_2} (x, \theta)$ is invertible under  condition \ref{assump:hypo}--II.  
 Thus, from (\ref{eq:inv_1}) and (\ref{eq:inv_2}), we obtain 
 \begin{align} \label{eq:L_S1S2}
 \Lambda_{S_1 S_2} (x, \theta) 
 = - \tfrac{1}{2} \Lambda_{S_1 S_1} (x, \theta) \partial_{x_{S_2}}^\top V_{S_1, 0} (x_S, \beta_{S_1}), 
 \end{align} 
and the proof is now complete.  
\end{proof}
\begin{lemma} \label{lemma:lambda}
Assume that condition \ref{assump:hypo}-II holds. We have  that, for any $(x, \theta) \in \mathbb{R}^N \times \Theta$: 
\begin{gather}
\Lambda_{S_1 S_1} ( x,  \theta ) = 720 \, a_{S_1}^{-1} (x, \theta); \label{eq:L_S1S1} \\[0.1cm] 
\Lambda_{S_2 S_2} (x, \theta) = 12 a_{S_2}^{-1} (x, \theta)
- \tfrac{1}{2} 
\Lambda_{S_2 S_1} (x, \theta) 
\partial_{x_{S_2}}^\top V_{S_1, 0} (x_S,  \beta_{S_1}).  \label{eq:L_S2S2} 
\end{gather}
\end{lemma}
\begin{proof}%
\noindent
({\it{Proof of (\ref{eq:L_S1S1})}}). 
First, we note that the matrices 
$\Sigma_{S_1 S_1}(x, \theta), \,  \Sigma_{S_2 S_2}(x, \theta), \, \Sigma_{RR }(x, \theta)$
are invertible for any $(x, \theta) \in \mathbb{R}^N \times \Theta$ under condition \ref{assump:hypo}-II. Due to the block expression of matrix $\Sigma (x, \theta)$ in (\ref{eq:Sigma}), we have: 
\begin{align*}
\Lambda_{S_1 S_1} (x, \theta) 
= \Bigl( \Sigma_{S_1 S_1} (x, \theta) - 
\widetilde{\Sigma} (x, \theta) \hat{\Lambda} (x, \theta) 
\widetilde{\Sigma} (x, \theta)^\top
\Bigr)^{-1}, 
\end{align*}
where $\hat{\Lambda} (x, \theta)$, the inverse of matrix $\hat{\Sigma}(x, \theta)$ given in (\ref{eq:Sigma_blocks}), has the following block expression:
\begin{align*}
    \hat{\Lambda} (x, \theta) 
    = 
    \begin{bmatrix}
        \hat{\Lambda}_{S_2 S_2} (x, \theta)  &  \hat{\Lambda}_{S_2 R}  (x, \theta)  \\[0.1cm] 
        \hat{\Lambda}_{S_2 R}  (x, \theta)^\top  &  \hat{\Lambda}_{R R}  (x, \theta)  
    \end{bmatrix}
\end{align*}
where we have set: 
\begin{gather*}
\hat{\Lambda}_{S_2 S_2} (x, \theta) = 4 \Sigma_{S_2 S_2}^{-1} (x, \theta), 
\quad 
\hat{\Lambda}_{S_2 R} (x, \theta) = - 4 \Sigma^{-1}_{S_2 S_2} (x, \theta) \Sigma_{S_2 R} (x, \theta) \Sigma_{RR}^{-1} (x, \sigma);  \\[0.1cm] 
\hat{\Lambda}_{R R} (x, \theta) = \Sigma_{RR}^{-1} (x, \sigma) 
+ 4 \Sigma_{RR}^{-1} (x, \sigma) \Sigma_{R S_2} (x, \theta) \Sigma^{-1}_{S_2 S_2} (x, \theta) \Sigma_{S_2 R} (x, \theta) 
\Sigma_{RR}^{-1} (x, \sigma). 
\end{gather*}
We then have: 
\begin{align*}
\hat{\Lambda} (x, \theta) \widetilde{\Sigma} (x, \theta)^\top 
= 
\begin{bmatrix} 
\tfrac{1}{2} \bigl( \partial_{x_{S_2}}^\top V_{S_1, 0} ( x_S, \beta_{S_1} ) \bigr)^\top  \\[0.1cm]
- \tfrac{1}{12} 
\bigl( 
\partial_{x_{S_2}}^\top V_{S_1, 0} ( x_S, \beta_{S_1} )  
\partial_{x_{R}}^\top V_{S_2, 0} ( x, \beta_{S_2} )  
\bigr)^\top 
\end{bmatrix}, 
\end{align*}
where we used (\ref{eq:inv_2}) for the upper block matrix, while the lower one is obtained via: 
\begin{align*}
& \hat{\Lambda}_{S_2 R}  (x, \theta)^\top \Sigma_{S_2 S_1} (x, \theta) 
+ \hat{\Lambda}_{R R}  (x, \theta) \Sigma_{R S_1} (x, \theta)
\\[0.2cm]
& = - 4 \Sigma_{RR}^{-1} (x, \sigma)  \Sigma_{R S_2} (x, \theta) \Sigma^{-1}_{S_2 S_2} (x, \theta) \Sigma_{S_2 S_1} (x, \theta) + \Sigma_{RR}^{-1} (x, \sigma) \Sigma_{R S_1} (x, \theta)  \\[0.2cm] 
& \quad \quad \quad + 4 \Sigma_{RR}^{-1} (x, \sigma) \Sigma_{R S_2} (x, \theta) \Sigma^{-1}_{S_2 S_2} (x, \theta) \Sigma_{S_2 R} (x, \theta) 
\Sigma_{RR}^{-1} (x, \sigma) \Sigma_{R S_1} (x, \theta) 
\\[0.2cm]
& = - \tfrac{3}{4} \bigl( \partial_{x_{R}}^\top V_{S_2, 0} ( x, \beta_{S_2} ) \bigr)^\top 
\bigl( \partial_{x_{S_2}}^\top V_{S_1, 0} ( x_S, \beta_{S_1} ) \bigr)^\top  + \tfrac{1}{6} \bigl( \partial_{x_{R}}^\top V_{S_2, 0} ( x, \beta_{S_2} ) \bigr)^\top 
\bigl( \partial_{x_{S_2}}^\top V_{S_1, 0} ( x_S, \beta_{S_1} ) \bigr)^\top  \\[0.1cm]
& \quad \quad \quad \quad + \tfrac{1}{2} \bigl( \partial_{x_R}^\top V_{S_2, 0} (x, \beta_{S_2})  \bigr)^\top \Sigma^{-1}_{S_2 S_2} (x, \theta) \Sigma_{S_2 S_2} (x, \theta) \bigl( \partial_{x_{S_2}}^\top V_{S_1, 0} (x_S, \beta_{S_1}) \bigr )^\top \\[0.2cm] 
& = -\tfrac{1}{12} 
\bigl( 
\partial_{x_{S_2}}^\top V_{S_1, 0} ( x_S, \beta_{S_1} )
\partial_{x_{R}}^\top V_{S_2, 0} ( x, \beta_{S_2} ) 
\bigr)^\top.  
\end{align*}
Thus, we obtain: 
\begin{align*}
\Bigl( \Sigma_{S_1 S_1} (x, \theta) - 
\widetilde{\Sigma} (x, \theta) \hat{\Lambda} (x, \theta) 
\widetilde{\Sigma} (x, \theta)^\top \Bigr)^{-1}
= \bigl( \tfrac{1}{36} \Sigma_{S_1 S_1} (x, \theta) \bigr)^{-1} 
= 720 \, a_{S_1}^{-1} (x, \theta), 
\end{align*}
and now the proof of (\ref{eq:L_S1S1}) is complete. 
\\

\noindent
({\it{Proof of (\ref{eq:L_S2S2})}}). From the block expression of the matrix $\Sigma (x, \theta)$ in (\ref{eq:Sigma}), we obtain:
\begin{align*}
\begin{bmatrix}
    \Lambda_{S_2 S_2} (x, \theta) &  \Lambda_{S_2 R} (x, \theta) 
    \\[0.1cm]
    \Lambda_{R S_2} (x, \theta) &  \Lambda_{R R} (x, \theta)   
\end{bmatrix}
 = \hat{\Lambda} (x, \theta)  
 + \hat{\Lambda} (x, \theta) \widetilde{\Sigma} (x, \theta)^\top
 \Lambda_{S_1 S_1} (x, \theta) 
 \widetilde{\Sigma} (x, \theta)  \hat{\Lambda} (x, \theta).  
\end{align*}
Thus, we have: 
\begin{align*}
\Lambda_{S_2 S_2} (x, \theta)
& = \hat{\Lambda}_{S_2 S_2} (x, \theta) 
+ \Xi (x, \theta)^\top \Lambda_{S_1 S_1} (x, \theta) 
 \Xi (x, \theta)  \\[0.1cm]
& = 4 \Sigma^{-1}_{S_2 S_2} (x, \theta) 
+ \tfrac{1}{4} \bigl( \partial_{x_{S_2}}^\top V_{S_1, 0} (x_S, \beta _{S_1}) \bigr)^\top  \Lambda_{S_1 S_1} (x, \theta) 
\partial_{x_{S_2}}^\top V_{S_1, 0} ( x_S, \beta _{S_1} )  \\[0.1cm]
& = 12 a_{S_2}^{-1} (x, \theta) - \tfrac{1}{2} \Lambda_{S_2 S_1} (x, \theta)  \partial_{x_{S_2}}^\top V_{S_1, 0} ( x_S, \beta _{S_1}), 
\end{align*}
where $\Xi (x, \theta)$ is defined as in (\ref{eq:Xi}), and we made use of (\ref{eq:inv_2}), (\ref{eq:L_S1S2}). The proof of (\ref{eq:L_S2S2}) is now complete. 
\end{proof}
\begin{lemma} \label{lemma:matrix_2}
Assume that conditions \ref{assump:hypo}-II and (\ref{assump:coeff}) hold. For any $\theta = (\beta_S, \beta_R, \sigma) \in \Theta$ and $1 \le i \le N$, we have that
\begin{align*}
 \Lambda (\sample{i-1},  (\beta_S, \sigma))
 \bigl( m_i (\Delta_n, \theta) - m_i (\Delta_n, (\beta_S, \truebeta_R, \sigma)) \bigr)
 = \sqrt{\Delta_n}
 \begin{bmatrix}
 \mathbf{0}_{N_{S_1}} \\[0.1cm]
 \mathbf{0}_{N_{S_2}} \\[0.1cm]
 b_R
 \end{bmatrix}
 + R (\sqrt{\Delta_n^3}, \sample{i-1}, \theta),
\end{align*}
for $\mathbb{R}^N$-valued function $R$ with $R^j \in \mathcal{S}, \, 1 \le j \le N$, where the $N_R$-dimensional vector $b_R$ is specified as: 
\begin{align*}
  b_R \equiv  a_R^{-1} (\sample{i-1}, \sigma) 
  \, \eta_R (\sample{i-1}, \beta_R). 
\end{align*}
\end{lemma}
\begin{proof}
We have 
\begin{align*} 
m_{i} (\Delta_n, \theta) - m_{i} (\Delta_n, (\beta_S, \truebeta_R, \sigma)) & =  
\begin{bmatrix}
\tfrac{\sqrt{\Delta_n}}{6} \, \partial_{x_{S_2}}^\top  
V_{S_1, 0} (X_{S, {i-1}}, \beta_{S_1})  
\partial_{x_R}^\top  V_{S_2,0} (\sample{i-1}, \beta_{S_2})
\, \eta_{R} (\sample{i-1}, \beta_R) \\[0.2cm]
\tfrac{\sqrt{\Delta_n}}{2} \, \partial_{x_R}^\top V_{S_2,0} (x, \beta_{S_2})  
\, \eta_{R} (\sample{i-1}, \beta_R)
\\[0.2cm]
\sqrt{\Delta_n} \,  \eta_{R} (\sample{i-1}, \beta_R) 
\end{bmatrix} + R (\sqrt{\Delta_n^3}, \sample{i-1}, \theta)
\nonumber \\[0.3cm] 
& = 
\sqrt{\Delta_n} \begin{bmatrix}
\Sigma_{S_1 R} \bigl(\sample{i-1}, \theta \bigr) 
\\[0.2cm]
\Sigma_{S_2 R} \bigl( \sample{i-1}, \theta \bigr)   \\[0.2cm] 
\Sigma_{R R} \bigl( \sample{i-1}, \sigma \bigr) 
\end{bmatrix} 
a_{R}^{-1} (\sample{i-1}, \sigma) \eta_{R} (\sample{i-1}, \beta_R) 
+ R (\sqrt{\Delta_n^3}, \sample{i-1}, \theta),
\end{align*}
for $\mathbb{R}^N$-valued function $R$ with $R^j \in \mathcal{S}, \, 1 \le j \le N$, where we used for $(x, \theta) \in \mathbb{R}^N \times \Theta$, 
\begin{align*} 
  \mathcal{L} V_{S_2, 0} (x, \theta) 
  & = \partial_{x_R}^\top V_{S_2, 0} (x, \beta_{S_2}) V_{R, 0} (x, \beta_R) + v_{S_2} (x, \beta_S, \sigma); 
  \\[0.2cm]
  \mathcal{L}^2 V_{S_1, 0} (x, \theta)  
  & = 
  \partial_{x_{S_2}}^\top V_{S_1, 0} (x_S, \beta_{S_1})  \partial_{x_R}^\top V_{S_2, 0} (x, \beta_{S_2}) V_{R, 0} (x, \beta_R) + v_{S_1} (x, \beta_S, \sigma), 
\end{align*}
with some functions $v_{S_2} : \mathbb{R}^N \times \Theta_{\beta_S} \times \Theta_\sigma \to \mathbb{R}^{N_{S_2}}$ 
and $v_{S_1} : \mathbb{R}^N \times \Theta_{\beta_S} \times \Theta_\sigma \to \mathbb{R}^{N_{S_1}}$ that are independent of $\beta_R \in \Theta_{\beta_R}$.   
Thus, it follows that 
\begin{align*}
 \Lambda \bigl(\sample{i-1}, (\beta_S, \sigma) \bigr)
 \bigl(  m_{i} (\Delta_n, \theta) - m_{i} (\Delta_n, (\beta_S, \truebeta_R, \sigma))  \bigr)  
 = \sqrt{\Delta_n}
 \begin{bmatrix}
 b_{S_1} \\
 b_{S_2} \\
 b_{R}
 \end{bmatrix}
 a_R^{-1} (\sample{i-1}, \sigma) \eta_R (\sample{i-1}, \beta_R)
 + \widetilde{R} (\sqrt{\Delta_n^3}, \sample{i-1}, \theta), 
\end{align*}
for $\widetilde{R}^j \in \mathcal{S}, \, 1 \le j \le N$,  where we have set: 
\begin{align*}
\begin{bmatrix}
 b_{S_1} \\[0.1cm]
 b_{S_2} \\[0.1cm]
 b_{R}
 \end{bmatrix}
 =  
 \Lambda \big(\sample{i-1}, \theta \bigr) \cdot 
 \begin{bmatrix}
     \Sigma_{S_1 R} \bigl( \sample{i-1},  \theta \bigr) \\[0.2cm] 
     \Sigma_{S_2 R} \bigl( \sample{i-1}, \theta \bigr) \\[0.2cm]  
     \Sigma_{R R} \bigl( \sample{i-1},  \sigma \bigr) 
 \end{bmatrix}
 = 
 \begin{bmatrix}
 \mathbf{0}_{N_{S_1} \times N_{R}} \\[0.1cm] 
 \mathbf{0}_{N_{S_2} \times N_{R}} \\[0.1cm] 
 I_{N_{R} \times N_{R}} 
 \end{bmatrix}, 
\end{align*}
since it holds $\Lambda (x, \theta) \Sigma (x, \theta) = I_{N \times N}$ for each $(x, \theta) \in \mathbb{R}^N \times \Theta$. 
The proof is now complete. 
\end{proof}
\section{Proof of Proposition \ref{prop:positive_definite}} \label{appendix:positive_definite}
{
The first equation within condition \ref{assump:hypo}-II immediately leads to: 
\begin{align*}
 \inf_{(x, \sigma) \in \mathbb{R}^N \times \Theta} 
 | a_R (x, \sigma) | 
 > 0.  
\end{align*}
Furthermore, the first and second equations of condition \ref{assump:hypo}-II yield:   
\begin{align*} 
& 
\inf_{(x, \theta) \in \mathbb{R}^N \times \Theta}
\  \inf_{\substack{ \xi \in \mathbb{R}^{N_{S_2}} \\ 
\mrm{s.t.} \, \| \xi \| = 1}} 
\sum_{k = 1}^d 
\Bigl\langle   
\mathrm{proj}_{N_{S_{1}}+1,N_{S_1} + N_{S_2}}
\bigl\{ \bigl[ \widetilde{V}_0, V_{k} \bigr](x, \theta) \bigr\} 
, \ \xi 
\Bigr\rangle^2 > 0 . 
\end{align*} 
Noticing that $\mathrm{proj}_{N_{S_{1}}+1,N_{S_1} + N_{S_2}} \bigl\{ \bigl[ \widetilde{V}_0, V_{k} \bigr](x, \theta) \bigr\} = \partial_{x_R}^\top V_{S_2, 0} (x, \beta_{S_2})  V_{R,k} (x , \sigma)$, we obtain: 
\begin{align*}
\inf_{(x, \theta) \in \mathbb{R}^N \times \Theta} 
| a_{S_2} (x, \theta) | > 0. 
\end{align*} 
Similarly, under condition \ref{assump:hypo}-II,  
\begin{align*} 
\inf_{(x, \theta) \in \mathbb{R}^N \times \Theta}
\  
\inf_{ \substack{ 
\xi \in \mathbb{R}^{N_{S_1}} \\
\mrm{s.t.} \, \| \xi \| = 1}} \sum_{k = 1}^d 
\Bigl\langle   
\mathrm{proj}_{1,N_{S_1}}
\bigl\{ \bigl[ \widetilde{V}_0, [\widetilde{V}_0, V_{k}] \bigr](x, \theta) \bigr\} 
, \ \xi \Bigr\rangle^2 > 0.  
\end{align*}
Since $\mathrm{proj}_{1,N_{S_1}}
\bigl\{ \bigl[ \widetilde{V}_0, [\widetilde{V}_0, V_{k}] \bigr](x, \theta) \bigr\} = \partial_{x_{S_2}}^\top V_{S_1, 0} (x_S, \beta_{S_1})
\partial_{x_R}^\top V_{S_2, 0} (x, \beta_{S_2}) V_{R,k} (x , \sigma)$, we have that:  
\begin{align*}
\inf_{(x, \theta) \in \mathbb{R}^N \times \Theta_\sigma} 
| a_{S_1} (x, \theta) | > 0. 
\end{align*}  
Finally, from the expression of the covariance $\Sigma (\Delta, x, \theta)$ in (\ref{eq:covariance_lg}), its determinant is given as: 
\begin{align*} 
| \Sigma (\Delta, x, \theta) | 
=  \tfrac{\Delta^9}{8640} \, 
| a_{R} (x, \sigma) | \, 
| a_{S_1} (x, \theta) |\, 
| a_{S_2} (x, \theta) |.   
\end{align*} 
Thus, (\ref{eq:det_cov}) holds and the proof is complete. 
}
%
%
%
%
\section{Proof of Main Results} \label{sec:pf} 
In this section
we prove the main results, i.e.~Theorem \ref{thm:consistency}, \ref{thm:clt} in Section \ref{sec:main} of  the main text. The proofs make use of some technical results from Appendix \ref{appendix:pf_tech}. 
\subsection{Proof of Theorem \ref{thm:consistency} -- Consistency} 
\label{sec:pf_consistency}
To show consistency, we study the limit of the contrast function $\ell_{n}  (\theta)$, defined in (\ref{eq:contrast}), that involves terms such as
\begin{align*}
\tfrac{X_{S_1, i+1}- \mu_{S_1} (\Delta_n, X_{i}, \theta)}{\sqrt{\Delta_n^{5}}}, \ \ 
\tfrac{X_{S_2, {i+1}} - \mu_{S_2} (\Delta_n, X_{i} , \theta)}{\sqrt{\Delta_n^{3}}}, 
\ \ 1 \le i \le n-1, 
\end{align*}
where $\{ X_i \}_{i = 0, \ldots, n} $ are discrete-time observations under the true model (\ref{eq:hypo-II}) with parameter $\trueparam$. Then, the stochastic Taylor expansion for $X_{S, {i+1}}$ yields  
\begin{align} 
\begin{aligned} \label{eq:diff_terms}
\tfrac{X_{S_1, i+1} - \mu_{S_1} (\Delta_n, X_{i}, \theta)}{\sqrt{\Delta_n^{5}}} 
 & = \tfrac{ 
   V_{S_1, 0} ( X_{S, i}, \truebeta_{S_1}) 
   - V_{S_1, 0} (X_{S, i}, \beta_{S_1}) }{\sqrt{\Delta_n^3}}  
   + R_{S_1}(\Delta_n, X_{i} , \theta);  \\[0.2cm]
 \tfrac{X_{S_2, {i+1}} - \mu_{S_2} ( \Delta_n, X_{i} , \theta)}{\sqrt{\Delta_n^{3}}} 
 & = \tfrac{ 
   V_{S_2, 0} ( X_{i}, \truebeta_{S_2}) 
   - V_{S_2, 0} (X_{i}, \beta_{S_2})  }{\sqrt{\Delta_n}}  
   + R_{S_2} (\Delta_n, X_{i}, \theta), 
\end{aligned}
\end{align}
where $R_{S_1}^{j_1}, \, R_{S_2}^{j_2} \in \mathcal{S}$ for $1 \le j_1 \le N_{S_1}, \, 1 \le j_2 \le N_{S_2}$. Careful steps are needed to control the first terms in the right-hand sides of (\ref{eq:diff_terms}) as $\Delta_n \to 0$, within the proof of consistency. Our proof proceeds with the following strategy which extends arguments used in \cite{iguchi:22}:
%
%
\begin{itemize}
    \item[Step 1.] We prove consistency, along with a convergence rate, for the estimator $\hat{\beta}_{S_1, n}$. That is, if \limit, then
    $
    \hat{\beta}_{S_1, n} \probconv \beta_{S_1}^\dagger.
    $
    In particular, we show the rate: 
    \begin{align} \label{eq:step1}
    \tfrac{1}{\sqrt{\Delta_n^3}} \bigl(\hat{\beta}_{S_1, n} - \beta_{S_1}^\dagger \bigr) \probconv 0. 
    \end{align}
    \item[Step 2.] Making use of the convergence rate in (\ref{eq:step1}), we prove consistency, along with a convergence rate, for the estimator $\hat{\beta}_{S_2, n}$. That is,  if \limit, then $\hat{\beta}_{S_2, n} \probconv \beta_{S_2}^\dagger$. In particular, we show the rate:
    \begin{align} \label{eq:step2}
        \tfrac{1}{\sqrt{\Delta_n}} \bigl(\hat{\beta}_{S_2, n} - \beta_{S_2}^\dagger \bigr) \probconv 0. 
    \end{align} 
    \item[Step 3.] Making use of the rates in (\ref{eq:step1}) and (\ref{eq:step2}), we prove consistency for the estimators $(\hat{\beta}_{R, n}, \hat{\sigma}_{n})$. That is, if \limit, then $(\hat{\beta}_{R, n}, \hat{\sigma}_{n}) \probconv (\truebeta_R, \truesigma)$.
\end{itemize}
%
We emphasise that in our proof of consistency, the condition $\Delta_n = o (n^{-1/2})$ is not required 
while \cite{glot:20} assumed the condition throughout the proof of consistency in the case of the degenerate diffusion class (\ref{eq:hypo-I}). Typically, in order to show the consistency of $\hat{\beta}_{R, n}$, \cite{glot:20} exploited the rates of convergence 
\begin{align*}
\sqrt{\tfrac{n}{\Delta_n}} ( \beta_{S, n} - \truebeta_{S}) \probconv 0, \quad 
\sqrt{n} (\sigma_n - \truesigma) \probconv 0 
\end{align*}
that are derived under the condition $\Delta_n = o (n^{-1/2})$. In contrast, in our strategy, the rates of convergence (\ref{eq:step1}) and (\ref{eq:step2}) are obtained without requiring $\Delta_n = o (n^{-1/2})$, and are put into effective use to avoid explosion of terms such as
\begin{align*}
\tfrac{V_{S_1, 0} (X_{S, {i-1}}, \truebeta_{S_1}) 
- V_{S_1, 0} (X_{S, {i-1}}, \hat{\beta}_{S_1, n})}{\sqrt{\Delta_n^3}}, \quad 
\tfrac{V_{S_2, 0} ( X_{i}, \truebeta_{S_2}) - V_{S_2, 0} (X_{i}, \hat{\beta}_{S_2, n})}{\sqrt{\Delta_n}} 
\end{align*}
as $\Delta_n \to 0$ only, with the help of some results derived from straightforward matrix calculations, i.e., Lemmas \ref{lemma:matrix}, \ref{lemma:lambda} and \ref{lemma:matrix_2} in Appendix \ref{app:aux}.
\begin{remark} \label{rem:a_inv}
{Due to (\ref{eq:inf_dets}) in Proposition \ref{prop:positive_definite}, we have under conditions \ref{assump:hypo}-II and (\ref{assump:coeff}) that: 
\begin{align*}
a_{S_1, i_1 j_1}^{-1} \in C_p^\infty (\mathbb{R}^N\times \Theta; \mathbb{R}), 
\quad 
a_{S_2, i_2 j_2}^{-1} \in C_p^\infty (\mathbb{R}^N\times \Theta; \mathbb{R}),  
\quad 
a_{R, i_3 j_3}^{-1} \in C_p^\infty (\mathbb{R}^N\times \Theta; \mathbb{R}),  
\end{align*}
for $1 \le i_1, j_1 \le N_{S_1}$, $1 \le i_2, j_2 \le N_{S_2}$,  and $1 \le i_3, j_3 \le N_{R}$. Thus, one is able to apply Lemmas \ref{lemma:ergodic_thm}, \ref{lemma:canonical_conv} for the (scaled) contrast function that includes the above inverse matrices in the proofs below. 
} 
\end{remark}
\subsubsection{Step 1}  \label{sec:step1}
Consistency of the estimator $\hat{\beta}_{S_1, n}$ is deduced from the following result.%
\begin{lemma} \label{lemm:step1}
Assume that conditions \ref{assump:hypo}-II and (\ref{assump:coeff})--(\ref{assump:finite_moment}) hold. If \limit, then,
\begin{align} \label{eq:conv_step1}
\tfrac{\Delta_n^3}{n} \ell_{n} (\theta) 
\probconv 
\int 
\eta_{S_1} (x, \beta_{S_1})^\top 
\Lambda_{S_1 S_1} (x, \theta) \, 
\eta_{S_1} (x, \beta_{S_1}) 
\truedist (dx),             
\end{align} 
uniformly in $\theta \in \Theta$. 
\end{lemma}
\noindent The proof is given in Appendix \ref{appendix:pf_step1}. Lemma \ref{lemm:step1} implies the consistency of $\hat{\beta}_{S_1, n}$ via the following discussion. {We write the right-hand side of (\ref{eq:conv_step1}) as $\mathscr{Q}(\theta), \, \theta = (\beta_{S_1}, \beta_{S_2}, \beta_R, \sigma) \in \Theta$. Under condition \ref{assump:hypo}-II, it follows from Proposition \ref{prop:positive_definite} and (\ref{eq:L_S1S1}) in Lemma \ref{lemma:lambda} that for all non-zero $\xi \in \mathbb{R}^{N_{S_1}}$, 
\begin{align*}
\inf_{(x, \theta) \in \mathbb{R}^N \times \Theta} \, 
\xi^\top \Lambda_{S_1 S_1} (x, \theta) \, \xi > 0.  
\end{align*}
Thus, it holds under the identifiablity condition (\ref{assump:ident}) that for any $\varepsilon > 0$ there exists a constant $\delta > 0$ such that 
\begin{align*}
\mathbb{P}_{\trueparam} 
\Bigl( \, \bigl| \hat{\beta}_{S_1, n}  - \truebeta_{S_1} \bigr| > \varepsilon 
\Bigr) 
\le 
\mathbb{P}_{\trueparam} 
\Bigl( \, 
\bigl| \mathscr{Q} (\hat{\theta}_{n}) 
\bigr|  > \delta 
\Bigr). 
\end{align*}
Then, Lemma \ref{lemm:step1} together with the definition of the estimator $\hat{\theta}_n$ yields: 
\begin{align}
\mathbb{P}_{\trueparam} 
\Bigl( \, 
\bigl| \mathscr{Q} (\hat{\theta}_{n}) 
\bigr|  > \delta 
\Bigr) 
& \le 
\mathbb{P}_{\trueparam} 
\Bigl( \, 
\bigl| \tfrac{\Delta_n^3}{n} \ell_{n} 
\bigl( \truebeta_{S_1}, \hat{\beta}_{S_2, n}, \hat{\beta}_{R, n}, \hat{\sigma}_{n} \bigr)   
- \tfrac{\Delta_n^3}{n} \ell_{n} 
\bigl( \hat{\theta}_n \bigr)  
+  \mathscr{Q} (\hat{\theta}_{n}) 
\bigr|  > \delta 
\Bigr) \nonumber  \\ 
& \le 
\mathbb{P}_{\trueparam} 
\Bigl( \, 
\sup_{\theta \in \Theta} 
\, 
\bigl| \tfrac{\Delta_n^3}{n} \ell_{n} 
\bigl( \truebeta_{S_1}, {\beta}_{S_2}, {\beta}_{R}, \sigma_n \bigr)   
- \tfrac{\Delta_n^3}{n} \ell_{n} 
\bigl( \theta \bigr)  
+  \mathscr{Q} ( \theta ) 
\bigr|  > \delta 
\Bigr) \to 0 
\nonumber
\end{align} 
as \limit, which leads to the consistency of $\hat{\beta}_{S_1, n}$. 
}
%
%
\\

We now prove the rate of convergence in (\ref{eq:step1}). Considering the Taylor expansion of $\partial_{\beta_{S_1}} \ell_{n} (\hat{\theta}_{n})$ around the derivative
$\partial_{\beta_{S_1}} \ell_{n} ( \truebeta_{S_1}, \hat{\beta}_{S_2,n}, \hat{\beta}_{R, n}, \hat{\sigma}_{n})$ with an appropriate scaling factor, we obtain   
\begin{align*} 
\mathscr{A}_{S_1, n}(  
 \truebeta_{S_1}, 
 \hat{\beta}_{S_2, n}, \hat{\beta}_{R, n}, \hat{\sigma}_{n} ) 
= \mathscr{B}_{S_1, n}
( 
 \hat{\theta}_{n} ) 
 \; \times 
 \tfrac{1}{\sqrt{\Delta_n^3}} (\hat{\beta}_{S_1, n} - \truebeta_{S_1}),  
\end{align*}
where we have set, for 
$\theta = (\beta_{S_1}, \theta^{- \beta_{S_1}} ) \in \Theta$ with $\theta^{- \beta_{S_1}} \equiv (\beta_{S_2}, \beta_R, \sigma)$, 
\begin{align*} 
\mathscr{A}_{S_1, n} ( \theta )  
=
- \tfrac{\sqrt{\Delta_n^{3}}}{n} 
\partial_{\beta_{S_1}} \ell_{n} 
( 
\theta ), \quad   
\mathscr{B}_{S_1, n} ( \theta )  
= \tfrac{\Delta_n^{3}}{n} 
\int_0^1  
\partial_{\beta_{S_1}}^2 \ell_{n}  
\bigl( \truebeta_{S_1} + \lambda ( {\beta}_{S_1} - \truebeta_{S_1}),  
\theta^{- \beta_{S_1}} \bigr) \, d \lambda.  
\end{align*}
For the matrix $\mathscr{B}_{S_1, n}( 
\theta)$, we have the following result: 
\begin{lemma} \label{lemma:step1_B}
Assume that conditions \ref{assump:hypo}-II and (\ref{assump:coeff})--(\ref{assump:ident}) hold. 
If \limit, then 
%
\begin{align*}
\mathscr{B}_{S_1, n}
( \hat{\beta}_{S_1,n},
\theta^{- \beta_{S_1}}  ) 
\probconv
\int 
\partial_{\beta_{S_1}} \bigl( V_{S_1, 0} (x_S, \truebeta_{S_1}) \bigr)^\top 
\Lambda_{S_1 S_1} \bigl(x, (\truebeta_{S_1},\theta^{- \beta_{S_1}}) \bigr)
\partial_{\beta_{S_1}}^\top V_{S_1, 0} (x_S, \truebeta_{S_1}) \, \truedist (dx) ,  
\end{align*}
uniformly in $\theta^{- \beta_{S_1}} \equiv (\beta_{S_2}, \beta_R, \sigma)  \in \Theta_{\beta_{S_2}} \times \Theta_{\beta_{R}} \times \Theta_{\sigma}$. 
\end{lemma}
\noindent We give the proof in Appendix \ref{appendix:pf_step1_B}. We now check that $\mathscr{A}_{S_1, n} ( \truebeta_{S_1}, \theta^{- \beta_{S_1}} ) \probconv \mathbf{0}_{N_{\beta_{S_1}}}$ as \limit, uniformly in $\theta^{- \beta_{S_1}} = (\beta_{S_2}, \beta_R, \sigma)$. 
We have: 
 \begin{align*}
 \mathscr{A}_{S_1, n} (  
 \truebeta_{S_1}, \theta^{- \beta_{S_1}} ) 
 =  \widetilde{\mathscr{A}}_{S_1, n} (  
 \truebeta_{S_1}, \theta^{- \beta_{S_1}} ) 
 + \tfrac{1}{n} \sum_{i = 1}^n 
R \big(\sqrt{\Delta_n}, \sample{i-1}, (\truebeta_{S_1}, \theta^{- \beta_{S_1}})\big), 
 \end{align*}
for $R^{j} \in \mathcal{S}, \, 1 \le j \le N_{\beta_{S_1}}$, where we have set, for $\theta = (\beta_{S_1}, \beta_{S_2}, \beta_R, \sigma) \in \Theta$, 
\begin{align*}
 \widetilde{\mathscr{A}}_{S_1, n} (  
 \theta ) = 
 \tfrac{1}{n \sqrt{\Delta_n}} \sum_{i = 1}^n 
 \partial_{\beta_{S_1}} \bigl(  V_{S_1, 0} (
 \sample{S, i-1}, \beta_{S_1}) \bigr)^\top \, 
 \Phi (\sample{i-1}, \theta) \, 
 \eta_{S_2} (\sample{i-1}, \beta_{S_2})
\end{align*}
 with $\Phi : \mathbb{R}^N \times \Theta \to \mathbb{R}^{N_{S_1} \times N_{S_2}}$ defined as:
 \begin{align} \label{eq:phi_0}
    \Phi (x ,\theta) = \Lambda_{S_1 S_1} (x, \theta)  \, 
    \partial_{x_{S_2}}^\top V_{S_1, 0} (x_S, \beta_{S_1}) 
    + 2 \Lambda_{S_1 S_2} (x, \theta). 
 \end{align}
From Lemmas \ref{lemma:ergodic_thm} and \ref{lemma:canonical_conv} in Appendix \ref{appendix:aux}, we immediately have that if \limit, then  
\begin{align*}
\tfrac{1}{n} \sum_{i = 1}^n R 
\bigl( \sqrt{\Delta_n}, \sample{i-1}, (\truebeta_{S_1},  \,  \theta^{- \beta_{S_1}}) \bigr) \probconv \mathbf{0}_{N_{\beta_{S_1}}},   
\end{align*} 
uniformly in $\theta^{- \beta_{S_1}}$. Furthermore, we have from Lemma \ref{lemma:matrix} that, for any $\theta \in \Theta$, $\widetilde{\mathscr{A}}_{S_1, n}( \theta ) = \mathbf{0}_{N_{\beta_{S_1}}}$ with probability~1, since $\Phi (x ,\theta) = \mathbf{0}_{N_{S_1}\times N_{S_2}}$ for any $(x, \theta) \in \mathbb{R}^N \times \Theta$. Hence, we obtain $ \mathscr{A}_{S_1, n} ( \truebeta_{S_1}, \theta^{- \beta_{S_1}}) \probconv \mathbf{0}_{N_{\beta_{S_1}}}$, 
and now  convergence (\ref{eq:step1}) holds. 
\subsubsection{Step 2}
Making use of convergence (\ref{eq:step1}), we obtain the following result whose proof is postponed to Appendix \ref{appendix:pf_step_2}.  
\begin{lemma} \label{lemm:step_2}
Assume that conditions \ref{assump:hypo}-II and (\ref{assump:coeff})--(\ref{assump:ident}) hold. 
If \limit, then 
\begin{align*}
 & \tfrac{\Delta_n}{n} \ell_{n} 
 ( \hat{\beta}_{S_1,n}, \beta_{S_2}, \beta_R, \sigma
)  \probconv 
 12 \times \int
  \eta_{S_2} (x, \beta_{S_2})^\top 
  a_{S_2}^{-1} \bigl(x, (\truebeta_{S_1}, \beta_{S_2}, \sigma) \bigr) 
  \, \eta_{S_2} (x, \beta_{S_2}) 
  \, \truedist (dx),
\end{align*}
uniformly in $(\beta_{S_2}, \beta_R, \sigma) \in \Theta_{\beta_{S_2}} \times \Theta_{\beta_{R}} \times  \Theta_\sigma$. 
\end{lemma}
\noindent This result leads to the consistency of $\hat{\beta}_{S_2, n}$ following an argument similar to the one used in \textbf{Step 1} to show consistency of $\hat{\beta}_{S_1, n}$. 
\\ 

To prove convergence (\ref{eq:step2}), we apply a Taylor expansion on the contrast function to get: 
\begin{align*} 
\mathscr{A}_{S, n} (  
 \truebeta_{S_1}, 
 \truebeta_{S_2}, \hat{\beta}_{R, n}, \hat{\sigma}_{n} ) 
= \mathscr{B}_{S, n} ( 
 \hat{\theta}_{n} ) 
 \; \times 
 \begin{bmatrix}
 \tfrac{1}{\sqrt{\Delta_n^3}} (\hat{\beta}_{S_1, n} - \truebeta_{S_1}) \\[0.3cm]  
 \tfrac{1}{\sqrt{\Delta_n}} (\hat{\beta}_{S_2, n} - \truebeta_{S_2})
\end{bmatrix}, 
\end{align*}
where we have set, for $\theta = (\beta_{S}, \beta_R, \sigma )\in \Theta$ with $\beta_S \equiv ( \beta_{S_1}, \beta_{S_2})$,  
\begin{align*}
 & \mathscr{A}_{S, n} (  \theta )  
 =
 \begin{bmatrix}
  - \tfrac{\sqrt{\Delta_n^3}}{n} \partial_{\beta_{S_1}} \ell_{n} ( \theta ) \\[0.2cm]
  - \tfrac{\sqrt{\Delta_n}}{n} \partial_{\beta_{S_2}} \ell_{n} ( \theta )
 \end{bmatrix}, \qquad 
 \mathscr{B}_{S,n} \bigl(  \theta \bigr)  
 = \int_0^1 M_{\beta_S, n} \,  \partial_{\beta_S}^2 \ell_{n}
 \bigl( \truebeta_S + \lambda ( \beta_S - \truebeta_S ), \beta_R, \sigma \bigr) \, M_{\beta_S, n} d \lambda,
\end{align*}
and $M_{\beta_S, n} \in \mathbb{R}^{N_{\beta_S} \times N_{\beta_S}}$ is defined as: 
\begin{align*}
M_{\beta_S, n} = \mrm{diag}
\Bigl( 
   \Bigl[ \, 
    \underbrace{ \sqrt{\tfrac{\Delta_n^3}{n}},
    \ldots, 
    \sqrt{\tfrac{\Delta_n^3}{n}}}_{N_{\beta_{S_1}}}, \ \ 
    \underbrace{\sqrt{\tfrac{\Delta_n}{n}},
    \ldots, 
    \sqrt{\tfrac{\Delta_n}{n}}}_{N_{\beta_{S_2}}} 
    \Bigr]^\top \Bigr). 
\end{align*}
%
%
Convergence (\ref{eq:step2}) is immediately deduced from the following result.
\begin{lemma} \label{lemma:step_2}
Assume that conditions \ref{assump:hypo}-II and (\ref{assump:coeff})--(\ref{assump:ident}) hold. If \limit, then 
\begin{align}
& \mathscr{A}_{S, n} ( \truebeta_{S}, \beta_R, \sigma)  
\probconv \mbf{0}_{N_{\beta_S}};  \label{eq:step2_A}   \\[0.2cm] 
& \mathscr{B}_{S, n} \bigl( \hat{\beta}_{S, n}, \beta_R, \sigma \bigr)  
\probconv 2 
\times 
\mrm{diag} \Bigl(\mathscr{B}_{S_1 S_1} ( \truebeta_S, \beta_R, \sigma ),
\, \mathscr{B}_{S_2 S_2} ( \truebeta_S, \beta_R, \sigma )
\Bigr),  
\label{eq:step2_B}  
%
%
\end{align}
uniformly in $(\beta_R, \sigma) \in \Theta_{\beta_R} \times \Theta_\sigma$, where $\truebeta_S \equiv (\truebeta_{S_1}, \truebeta_{S_2})$, $\hat{\beta}_{S, n} \equiv (\hat{\beta}_{S_1, n}, \hat{\beta}_{S_2, n})$ and we have set: 
%
%
\begin{align*}
\mathscr{B}_{S_1 S_1} (\theta) & = 
720 \int \partial_{\beta_{S_1}} \bigl( V_{S_1, 0} (x_S, \beta_{S_1 }) \bigr)^\top  
a_{S_1}^{-1} (x, \theta) 
\, \partial_{\beta_{S_1}}^\top V_{S_1, 0} (x_S, \beta_{S_1}) \, \truedist (dx);  \\[0.1cm]
\mathscr{B}_{S_2 S_2} (\theta) & = 
12 \int  \partial_{\beta_{S_2}} \bigl( V_{S_2, 0} (x, \beta_{S_2}) \bigr)^\top 
a_{S_2}^{-1} (x, \theta) 
\partial_{\beta_{S_2}}^\top V_{S_2, 0} (x, \beta_{S_2}) 
\, \truedist (dx),    
\end{align*} 
for $x \in \mathbb{R}^N, \, \theta = (\beta_S, \beta_R, \sigma) \in \Theta$, where $\beta_S = (\beta_{S_1}, \beta_{S_2})$.   
\end{lemma}
\noindent We give the proof in Appendix \ref{appendix:pf_step_4}. 

\subsubsection{Step 3}
Finally, we prove the consistency of estimators $(\hat{\beta}_{R, n}, \hat{\sigma}_{n})$. Working with the rates of convergence (\ref{eq:step1}) and (\ref{eq:step2}), we obtain the following result leading to the consistency of $\hat{\sigma}_{n}$: 
\begin{lemma} \label{lemma:step3_1}
Assume that conditions \ref{assump:hypo}-II and (\ref{assump:coeff})--(\ref{assump:ident}) hold. If \limit, then
\begin{align*}
\tfrac{1}{n} \, \ell_{n} \bigl( \hat{\beta}_{S, n},
\beta_R, \sigma \bigr) \nonumber 
\probconv 
\int \Bigl\{ 
\mrm{tr} \bigl( \Lambda (x, (\truebeta_S, \sigma)) 
\Sigma (x, (\truebeta_S, \truesigma))  \bigr)
+ \log \bigl| \Sigma (x,  (\truebeta_S, \sigma) ) \bigr| \Bigr\} \truedist (dx),  
\end{align*}
uniformly in $(\beta_R, \sigma) \in \Theta_{\beta_R} \times \Theta_{\sigma}$. 
\end{lemma}
\noindent We provide the proof in Appendix \ref{appendix:pf_step3_1}. 
To show the consistency of $\hat{\beta}_{R,n}$, we consider, for 
$\beta_R \in \Theta_{\beta_R}$, 
\begin{align*} 
\mathscr{L} ( \beta_R )  
:= \tfrac{1}{n \Delta_n}  \ell_{n} 
( \hat{\beta}_{S_1,n},  \hat{\beta}_{S_2, n},  \beta_R,  \hat{\sigma}_{n} ) 
- \tfrac{1}{n \Delta_n}  \ell_{n} 
( \hat{\beta}_{S_1, n},  \hat{\beta}_{S_2,  n},  \truebeta_R,  \hat{\sigma}_{n} ). 
\end{align*}
The consistency of estimator $\hat{\beta}_{R, n}$ is obtained via the following result whose proof is given in Appendix~\ref{appendix:pf_step3_2}.
\begin{lemma} \label{lemma:step3_2}
Assume that conditions \ref{assump:hypo}-II and (\ref{assump:coeff})--(\ref{assump:ident}) hold. If \limit, then
\begin{align*}
\mathscr{L} (\beta_R) \probconv 
\int 
\eta_{R} (x, \beta_R)^\top   
a_R^{-1} \bigl( x, \truesigma \bigr)
\, \eta_R (x, \beta_R) \, 
\truedist (dx),
\end{align*} 
%
%
uniformly in  $\beta_R \in \Theta_{\beta_R}$.
\end{lemma}
\noindent The proof of consistency for the contrast estimator $\hat{\theta}_n$ is now complete.
\subsection{Proof of Theorem \ref{thm:clt} -- Asymptotic Normality}  \label{sec:pf_asymptotic_normality}
We consider the Taylor expansion of the contrast function 
$\ell_{n} ( \theta )$: 
\begin{align*}
 \mathscr{C}_{n} ( \trueparam ) 
 = \int_0^1 
  \mathscr{I}_{n} \bigl(  \trueparam  + \lambda 
  ( \hat{\theta}_{n} - \trueparam )  \bigr) d \lambda \times
  {M}_{n} ( \hat{\theta}_{n} -  \trueparam )
\end{align*} 
where we have set, for $\theta \in \Theta$,  
\begin{gather*}
  \mathscr{C}_{n} (\theta) =
  - M_n^{-1} \,  \partial_\theta 
  \ell_{n}  ( \theta ), 
 \quad  
  \mathscr{I}_{n} (\theta) =
   M_{n}^{-1} \, \partial^2_\theta \ell_{n}
   (\theta ) \, M_n^{-1}, \quad 
 M_n = \mathrm{diag} ({v}_{n}),
\end{gather*} 
with the $N_\theta$-dimensional vector ${v}_{n}$ defined as: 
\begin{align*}
    v_n = 
    \Bigl[
    \, 
    \underbrace{\sqrt{\tfrac{n}{\Delta_n^3}},
    \ldots, 
    \sqrt{\tfrac{n}{\Delta_n^3}}}_{N_{\beta_{S_1}}}, \ \ 
    \underbrace{\sqrt{\tfrac{n}{\Delta_n}},
    \ldots, 
    \sqrt{\tfrac{n}{\Delta_n}}}_{N_{\beta_{S_2}}},  \ \  
     \underbrace{
     \sqrt{n \Delta_n},
     \ldots, 
     \sqrt{n \Delta_n}}_{N_{\beta_R}}, \ \ 
     \underbrace{
     \sqrt{n},
     \ldots, 
     \sqrt{n}
     }_{N_{\sigma}} \, 
    \Bigr]^\top.  
\end{align*}
The asymptotic normality immediately holds from the following two results -- their proofs are shown in Appendices \ref{appendix:pf_slln} and \ref{appendix:pf_clt}. 
\begin{lemma} \label{lemma:slln}
Assume that conditions \ref{assump:hypo}-II and (\ref{assump:coeff})--(\ref{assump:ident}) hold. 
If \limit, then 
\begin{align*}
 \mathscr{I}_{n} \bigl(  \trueparam  + \lambda ( \hat{\theta}_{n} - \trueparam)  \bigr) \probconv 2 \Gamma (\trueparam), 
\end{align*} 
uniformly in $\lambda \in [0,1]$, where the matrix $\Gamma(\trueparam)$ is defined as in (\ref{eq:precision_matrix}) in the main text.
\end{lemma}  
\begin{lemma} \label{lemma:clt}
Assume that conditions \ref{assump:hypo}-II and (\ref{assump:coeff})--(\ref{assump:ident}) hold.
If \limit, with $\Delta_n = o (n^{-1/2})$, then 
\begin{align*}
   \mathscr{C}_{n} (\trueparam) \distconv \mathscr{N} \bigl( \mathbf{0}_{N_\theta}, 4 \Gamma (\trueparam) \bigr). 
\end{align*} 
\end{lemma}
\noindent The proof of Theorem \ref{thm:clt} is now complete. 
\section{Proof of Technical Results}  \label{appendix:pf_tech} 
\subsection{Proof of Lemma \ref{lemm:step1}}  
We have that $\textstyle \tfrac{\Delta_n^3}{n} \ell_n (\theta)  = \sum_{1 \le i \le 4} \mathscr{E}_{i} (\theta)$, $\theta = (\beta_{S_1}, \beta_{S_2}, \beta_R, \sigma) \in \Theta$, where we have set:  
\begin{align*}
& 
\mathscr{E}_1 (\theta)  
= \tfrac{1}{n} \sum_{i = 1}^n 
\eta_{S_1} (X_{{i-1}}, \beta_{S_1})^\top  
\Lambda_{S_1 S_1} (\sample{i-1}, \theta) 
\, \eta_{S_1} (X_{{i-1}}, \beta_{S_1}); \\  
%
%
& \mathscr{E}_2 (\theta)  
= \tfrac{1}{n} \sum_{i = 1}^n 
 \sum_{1 \le j_1, j_2 \le N} 
 R_1^{j_1 j_2} (\Delta_n^{q_1}, \sample{i-1}, \theta)
 \, 
 m_{i}^{j_1} (\Delta, \trueparam) 
 \, 
 m_{i}^{j_2} (\Delta, \trueparam); 
 \\ 
& \mathscr{E}_3 (\theta)  
 = \tfrac{1}{n} \sum_{i = 1}^n 
 \sum_{1 \le j \le N} 
 R_2^j (\Delta_n^{q_2}, \sample{i-1}, \theta)
 \, m_{i}^{j} (\Delta, \trueparam);  \\ 
& \mathscr{E}_4 (\theta) 
 =\tfrac{1}{n} \sum_{i = 1}^n 
  R_3 (\Delta_n^{q_3}, \sample{i-1}, \theta),
\end{align*}
for some functions $R_1^{j_1 j_2}, R_2^j, R_3 \in \mathcal{S}$ and constants $q_1, q_2, q_3 \ge 1$. From Lemmas \ref{lemma:ergodic_thm}, \ref{lemma:canonical_conv}, we immediately have that as \limit, 
\begin{align*}
 \mathscr{E}_1 (\theta) 
  \probconv 
\int 
\eta_{S_1} (x , \beta_{S_1})^\top  
\Lambda_{S_1 S_1} (x, \theta) 
\, \eta_{S_1} (x , \beta_{S_1}) 
\truedist (dx), \qquad 
\mathscr{E}_k ( \theta )  \probconv 0, \ \ 2 \le k \le 4, 
\end{align*}
uniformly in $\theta = (\beta_{S_1}, \beta_{S_2}, \beta_R, \sigma) \in \Theta$, and the proof is now complete. 
\label{appendix:pf_step1}

\subsection{Proof of Lemma \ref{lemma:step1_B}}  
\label{appendix:pf_step1_B}
We define $\mathscr{F} : \Theta \to \mathbb{R}^{N_{\beta_{S_1}} \times N_{\beta_{S_1}}}$ as: 
\begin{align*}
 \mathscr{F} (\theta ) =
 \tfrac{\Delta_n^3}{n}
  \partial_{\beta_{S_1}}^2 \ell_n 
  \left( \theta \right), \quad \theta \in \Theta.  
\end{align*} 
%
$\mathscr{F} (\theta)$ can be expressed as $\textstyle{\mathscr{F}(\theta) = \sum_{1 \le k \le 6} \mathscr{F}_k (\theta)}$, where we have set, for $1 \le j_1, j_2 \le N_{\beta_{S_1}}$ with multi-index $\mathbf{j} = (j_1, j_2)$, 
\begin{align*}
[\mathscr{F}_1 (\theta)]_{j_1 j_2} 
  & = 
 \tfrac{2}{n} \sum_{i = 1}^n 
 \bigl( \partial^{\beta_{S_1}}_{j_1} V_{S_1, 0} (\sample{S, i-1}, \beta_{S_1}) \bigr)^\top 
 \Lambda_{S_1 S_1} (\sample{i-1}, \theta) \partial^{\beta_{S_1}}_{j_2}  
 V_{S_1, 0} (\sample{S, i-1}, \beta_{S_1});  \\ 
 [\mathscr{F}_2 (\theta)]_{j_1 j_2} 
 & =  \tfrac{1}{n}  
 \sum_{i=1}^n  
 \eta_{S_1} (\sample{i-1}, \beta_{S_1})^\top
 \partial_{\mathbf{j}}^{\beta_{S_1}}  
 \Lambda_{S_1 S_1}(\sample{i-1}, \theta) \, 
 \eta_{S_1 } (\sample{i-1}, \beta_{S_1});  
 \\ 
 [\mathscr{F}_3 (\theta)]_{j_1 j_2} 
 & =
 - \tfrac{2}{n} \sum_{i = 1}^n 
 \eta_{S_1} (\sample{i-1}, \beta_{S_1})^\top
\partial_{\mathbf{j}}^{\beta_{S_1}}   
\bigl\{  
\Lambda_{S_1 S_1} (\sample{i-1}, \theta) 
 V_{S_1, 0} (\sample{S, i-1}, \beta_{S_1}) 
\bigr\};  \\
 [\mathscr{F}_4 (\theta)]_{j_1 j_2}
 & = \tfrac{1}{n} \sum_{i = 1}^n 
 \sum_{1 \le k_1, k_2 \le N} 
 R_{k_1 k_2}^{j_1 j_2} (\Delta_n^{q_1}, \sample{i-1}, \theta)
 \, 
 m_{i}^{k_1} (\Delta_n, \trueparam) 
 \, 
 m_{i}^{k_2} (\Delta_n, \trueparam); 
 \\ 
 [\mathscr{F}_5 (\theta)]_{j_1 j_2}
 & = \tfrac{1}{n} \sum_{i = 1}^n 
 \sum_{1 \le k \le N} R_{k}^{j_1 j_2} (\Delta_n^{q_2}, \sample{i-1}, \theta) \, m_{i}^{k} (\Delta_n, \trueparam); 
 \\ 
 [\mathscr{F}_6 (\theta)]_{j_1 j_2}
 & = \tfrac{1}{n} \sum_{i = 1}^n  R^{j_1 j_2} (\Delta_n^{q_3}, \sample{i-1}, \theta),
\end{align*} 
for some functions $R^{j_1 j_2}_{k_1 k_2}, R^{j_1 j_2}_{k}, R^{j_1 j_2} \in \mathcal{S}$ and constants $q_1, q_2, q_3 > 0$. It follows from Lemmas \ref{lemma:ergodic_thm}, \ref{lemma:canonical_conv} and the consistency of estimator $\hat{\beta}_{S_1, n}$ that if \limit, 
%
\begin{align*}
 & \Bigl[ \mathscr{F}_1 \bigl( \truebeta_{S_1} - \lambda (\hat{\beta}_{S_1,n} - \truebeta_{S_1} ), \beta_{S_2}, \beta_R, \sigma  \bigr) \Bigr]_{j_1 j_2} \\
 & \qquad  \probconv 2 \int \bigl( \partial^{\beta_{S_1}}_{j_1} V_{S_1, 0} (x_S, \truebeta_{S_1}) \bigr)^\top 
 \Lambda_{S_1 S_1} \bigl(x, (\truebeta_{S_1}, \beta_{S_2}, \sigma) \bigr) \partial^{\beta_{S_1}}_{j_2}  
 V_{S_1, 0} (x_S, \truebeta_{S_1}) 
 \truedist (dx); \\[0.2cm] 
 & \Bigl[ \mathscr{F}_k \bigl( \truebeta_{S_1} - \lambda (\hat{\beta}_{S_1,n} - \truebeta_{S_1} ), \beta_{S_2}, \beta_R, \sigma  \bigr) \Bigr]_{j_1 j_2}
 \probconv 0, \qquad  2 \le k \le 6, \\  
\end{align*} 
uniformly in $(\beta_{S_2}, \beta_R, \sigma) \in \Theta_{\beta_{S_2}} \times \Theta_{\beta_{R}} \times\Theta_{\sigma}$ and $\lambda \in [0,1]$. The proof is now complete. 
\subsection{Proof of Lemma \ref{lemm:step_2}} 
\label{appendix:pf_step_2}
We write $\theta^{- S_1}
\equiv \bigl( \beta_{S_2}, \beta_R, \sigma \bigr) \in \Theta_{\beta_{S_2}} \times \Theta_{\beta_R} \times \Theta_\sigma$. It holds that 
\begin{align*}
\tfrac{\Delta_n}{n} \, \ell_{n} 
 ( \hat{\beta}_{S_1, n}, \,  \theta^{-S_1} ) 
= \sum_{1 \le k \le 7} \mathscr{G}_k (\hat{\beta}_{S_1, n}, \,  \theta^{- S_1}), 
\end{align*} 
where we have set, for $\theta \in \Theta$, 
\begin{align*} 
  \mathscr{G}_1 (\theta) 
 & =  \tfrac{\Delta_n}{n} \sum_{i = 1}^n 
  \eta_{S_1, i-1} (\sqrt{\Delta_n^3}, \beta_{S_1})^\top 
  \Lambda_{S_1 S_1} (\sample{i-1}, \theta)  
  \eta_{S_1, i-1} (\sqrt{\Delta_n^3}, \beta_{S_1}); \\[0.2cm]
\mathscr{G}_2  (\theta) 
& =  \tfrac{1}{n} \sum_{i = 1}^n 
\sum_{1 \le j \le N} 
R_j (\sqrt{\Delta_n},  \sample{i-1}, \theta ) 
\eta_{S_1, i-1}^{j} (\sqrt{\Delta_n^3}, \beta_{S_1});  \\[0.2cm]  
\mathscr{G}_3 (\theta) 
& = \tfrac{1}{n} \sum_{i = 1}^n 
\sum_{1 \le j_1, j_2 \le N} 
R_{j_1 j_2} ({\Delta_n},  \sample{i-1}, \theta ) \eta_{S_1, i-1}^{j_1} (\sqrt{\Delta_n^3}, \beta_{S_1})
m_{i}^{j_2} (\Delta_n, \trueparam );  \\[0.2cm] 
 \mathscr{G}_4 (\theta) 
& = \tfrac{1}{n} \sum_{i = 1}^n
\Biggl\{ 
\tfrac{1}{4} 
\Bigl( 
\mathcal{L} V_{S_1, 0} (\sample{i-1}, \trueparam)
- \mathcal{L} V_{S_1, 0} (\sample{i-1}, \theta) 
\Bigr)^\top 
\Lambda_{S_1 S_1} (\sample{i-1}, \theta)
\, \Bigl( 
\mathcal{L} V_{S_1, 0} (\sample{i-1}, \trueparam)
- \mathcal{L} V_{S_1, 0} (\sample{i-1}, \theta) 
\Bigr) \\
& \qquad \qquad  
+ \tfrac{1}{2} 
\Bigl( 
\mathcal{L} V_{S_1, 0} (\sample{i-1}, \trueparam)
- \mathcal{L} V_{S_1, 0} (\sample{i-1}, \theta) 
\Bigr)^\top 
\Lambda_{S_1 S_2} (\sample{i-1}, \theta) 
\, \eta_{S_2} (\sample{i-1}, \beta_{S_2})
\\
& \qquad \qquad  
+ \tfrac{1}{2}
\eta_{S_2} (\sample{i-1}, \beta_{S_2})^\top 
\Lambda_{S_2 S_1} (\sample{i-1}, \theta) 
\, \Bigl( 
\mathcal{L} V_{S_1, 0} (\sample{i-1}, \trueparam)
- \mathcal{L} V_{S_1, 0} (\sample{i-1}, \theta) 
\Bigr)  
\\
& \qquad \qquad  
+ \eta_{S_2} (\sample{i-1}, \beta_{S_2})^\top
\Lambda_{S_2 S_2} (\sample{i-1}, \theta)  
\eta_{S_2} (\sample{i-1}, \beta_{S_2})
\Biggr\}; \\ 
\mathscr{G}_5 (\theta) 
& = \tfrac{1}{n} \sum_{i = 1}^n 
\sum_{1 \le j \le N} 
\widetilde{R}_{j} ({\Delta_n},  \sample{i-1}, \theta ) \,
m_{i}^{j} (\Delta_n, \trueparam);  \\[0.2cm]  
\mathscr{G}_6 (\theta) 
& = \tfrac{\Delta_n}{n} \sum_{i = 1}^n  
\Bigl\{ m_{i} (\Delta_n, \trueparam)^\top \Lambda (\sample{i-1}, \theta) 
\, m_{i} (\Delta_n, \trueparam )
+ \log \bigl| \Sigma (\sample{i-1}, \theta) \bigr| \Bigr\},  
 %
\end{align*} 
for some functions $R_j, R_{j_1 j_2}, \, \widetilde{R}_{j} \in \mathcal{S}$. Note that $\eta_{S_1, i-1}^j (\sqrt{\Delta_n^3}, \hat{\beta}_{S_1, n})$ can be expressed as:  
\begin{align*}
 \eta_{S_1, i-1}^j (\sqrt{\Delta_n^3}, \hat{\beta}_{S_1, n})
 = \sum_{1 \le k \le N_{\beta_{S_1}}} \tfrac{\eta_{S_1}^{j, [k]} (\sample{S, i-1})}{
 | \beta_{S_1}^{\dagger, k} - \hat{\beta}_{S_1, n}^k | } 
 \times \Bigl| \tfrac{\beta_{S_1}^{\dagger, k} - \hat{\beta}_{S_1, n}^k}{\sqrt{\Delta_n^3}} \Bigr|,  
\end{align*}
where we have set: 
\begin{align*} 
\eta_{S_1}^{j, [k]} (\sample{S, i-1}) \equiv V_{S_1, 0}^j ( \sample{S, i-1},  \bar{\beta}_{S_1, n}^{[k-1]}) 
- V_{S_1, 0}^j ( \sample{S, i-1},  \bar{\beta}_{S_1, n}^{[k]}),  \quad  1 \le k \le N_{\beta_{S_1}}, 
\end{align*}
with the notation:
\begin{gather*}
\bar{\beta}_{S_1, n}^{[\ell]} = (\hat{\beta}_{S_1, n}^{1}, \ldots, \hat{\beta}_{S_1, n}^{\ell},
{\beta}_{S_1}^{\dagger, \ell+1}, \ldots,  {\beta}_{S_1}^{\dagger, N_{\beta_{S_1}}}), \quad 1 \le \ell \le N_{\beta_{S_1}} -1;\\ 
\bar{\beta}_{S_1, n}^{[0]} = \truebeta_{S_1}, \quad  
\bar{\beta}_{S_1, n}^{\, [N_{\beta_{S_1}}]} = \hat{\beta}_{S_1, n}.  
\end{gather*} 
From  convergence (\ref{eq:step1}), Lemma \ref{lemma:ergodic_thm} and condition (\ref{assump:bd_deriv}) it follows that if \limit, then 
\begin{align}
\begin{aligned}
\tfrac{1}{n} \sum_{i = 1}^n f (\sample{i-1}, \theta) \eta_{S_1, i-1}^j (\sqrt{\Delta_n^3}, \hat{\beta}_{S_1, n})  
\probconv 0,  \label{eq:eta_conv}
\end{aligned}
\end{align} 
uniformly in $\theta \in \Theta = (\beta_{S_2}, \beta_R, \sigma) \in \Theta_{\beta_{S_2}} \times \Theta_{\beta_R} \times \Theta_\sigma$ for any $f: \mathbb{R}^N \times \Theta \to \mathbb{R}$ satisfying the same property in Lemma~\ref{lemma:ergodic_thm}. 
%
%
Thus, we obtain 
\begin{align*} 
 \mathscr{G}_k ( \hat{\beta}_{S_1, n}, \theta^{-S_1} )  \probconv 0, \qquad   1 \le k \le 3, 
\end{align*}
uniformly in $\theta^{-S_1}$. Also, we immediately have from Lemmas \ref{lemma:ergodic_thm}, \ref{lemma:canonical_conv} that
\begin{gather*} 
 \mathscr{G}_k  ( \hat{\beta}_{S_1, n}, \theta^{-S_1}) 
  \probconv  0, \qquad  k = 5,6, 
\end{gather*} 
as \limit, uniformly in $\theta^{-S_1}$. Finally, we consider the fourth term $\mathscr{G}_4$. Noticing that  
\begin{align*} 
    & \mathcal{L} V_{S_1, 0} (\sample{i-1}, \trueparam)
    - \mathcal{L} V_{S_1, 0} (\sample{i-1}, \theta) \\ 
    & \quad = 
    \partial_{x_{S_1}}^\top V_{S_1, 0} (\sample{S, i-1}, \truebeta_{S_1}) \, 
    \eta_{S_1} (\sample{i-1}, \beta_{S_1})  
    + \Bigl( \partial_{x_{S_1}}^\top V_{S_1, 0} (\sample{S, i-1}, \truebeta_{S_1}) - \partial_{x_{S_1}}^\top V_{S_1, 0} (\sample{S, i-1}, \beta_{S_1}) \Bigr) 
    V_{S_1, 0} (\sample{S, i-1}, \beta_{S_1}) 
    \\
    & \qquad 
    + \partial_{x_{S_2}}^\top V_{S_1, 0} (\sample{S, i-1}, \truebeta_{S_1}) \, \eta_{S_2} (\sample{i-1}, \beta_{S_2}) 
    + \Bigl( \partial_{x_{S_2}}^\top V_{S_1, 0} (\sample{S, i-1}, \truebeta_{S_1}) - \partial_{x_{S_2}}^\top V_{S_1, 0} (\sample{S, i-1}, \beta_{S_1}) \Bigr) 
    V_{S_2, 0} (\sample{i-1}, \beta_{S_2}), 
\end{align*}
we obtain from Lemmas \ref{lemma:ergodic_thm}, \ref{lemma:canonical_conv} and consistency of $\hat{\beta}_{S_1, n}$ that 
\begin{align*}
    \mathscr{G}_4 ( \hat{\beta}_{S_1, n}, \theta^{-S_1} ) 
    \probconv \int 
    \eta_{S_2} (x, \beta_{S_2})^\top  \,   
    \widetilde{\mathscr{G}}_4 \bigl( x,  
    (\truebeta_{S_1}, \beta_{S_2}, \beta_R, \sigma) \bigr) 
    \, \eta_{S_2} (x, \beta_{S_2})  \, \truedist (dx),
\end{align*}
as \limit, uniformly in $\theta^{-S_1}$, where we have set: for $x = (x_S, x_R) \in \mathbb{R}^{N_S} \times \mathbb{R}^{N_R}$ and $\theta = (\beta_{S_1}, \beta_{S_2}, \beta_R, \sigma) \in \Theta$, 
\begin{align*}
    \widetilde{\mathscr{G}}_4 \bigl( x,  \theta \bigr) 
    & = \tfrac{1}{4} \bigl( \partial_{x_{S_2}}^\top V_{S_1, 0} (x_S, \beta_{S_1}) \bigr)^\top 
    \Lambda_{S_1 S_1} (x, \theta) 
    \, \partial_{x_{S_2}}^\top V_{S_1, 0} (x_S, \beta_{S_1}) \\[0.2cm]  
    & \quad
    + \tfrac{1}{2} \bigl( \partial_{x_{S_2}}^\top V_{S_1, 0} (x_S, \beta_{S_1}) \bigr)^\top 
    \Lambda_{S_1 S_2} (x, \theta) 
    + \tfrac{1}{2} \Lambda_{S_2 S_1} (x, \theta) 
    \,\partial_{x_{S_2}}^\top V_{S_1, 0} (x_S, \beta_{S_1})
    + \Lambda_{S_2 S_2} (x, \theta).   
\end{align*} 
Lemmas \ref{lemma:matrix}, \ref{lemma:lambda} yield: 
\begin{align*}
    \widetilde{\mathscr{G}}_4 \bigl( x,  \theta \bigr) 
    & = \tfrac{1}{2} \Lambda_{S_2 S_1} (x, \theta) 
    \,\partial_{x_{S_2}}^\top V_{S_1, 0} (x_S, \beta_{S_1})
    + \Lambda_{S_2 S_2} (x, \theta)  = 12 a_{S_2}^{-1} (x, \theta). 
\end{align*}
The proof is now complete. 
\subsection{Proof of Lemma \ref{lemma:step_2}} 
\label{appendix:pf_step_4} 
\noindent 
(\textit{Proof of (\ref{eq:step2_A})}). 
Since the proof of $\tfrac{\sqrt{\Delta_n^3}}{n} \partial_{\beta_{S_1}} \ell_{n} (\truebeta_S, \beta_R, \sigma) \probconv \mathbf{0}_{N_{\beta_{S_1}}}$ is identical with that in Section \ref{sec:step1}, we will only show that, if \limit, then
\begin{align} \label{eq:K}
  \mathscr{K} (\truebeta_S, \beta_R, \sigma)  
  := \tfrac{\sqrt{\Delta_n}}{n} \partial_{\beta_{S_2}} 
  \ell_{n} (\truebeta_S, \beta_R, \sigma) 
  \probconv \mathbf{0}_{N_{\beta_{S_2}}}, 
\end{align}
uniformly in $(\beta_R, \sigma) \in \Theta_{\beta_R} \times \Theta_\sigma$. It holds that 
$\textstyle \mathscr{K} (\truebeta_S, \beta_R, \sigma) =  \sum_{1 \le l \le 3} \mathscr{K}_l (\truebeta_S, \beta_R, \sigma)$,
where we have set, for $\theta \in \Theta$, 
\begin{align*}
  \mathscr{K}_1 ( \theta )  
  & = \tfrac{1}{n} 
  \sum_{i = 1}^n \sum_{1 \le j \le N} 
  R_j (1, \sample{i-1}, \theta) \, m_{i}^j (\Delta_n, \trueparam);
  \\[0.2cm]
  \mathscr{K}_2 ( \theta )  
  & = \tfrac{1}{n} \sum_{i = 1}^n \sum_{1 \le j_1, j_2 \le N} 
  R_{j_1 j_2} (\sqrt{\Delta_n}, \sample{i-1}, \theta) 
  \, m_{i}^{j_1} (\Delta_n, \trueparam) 
  \, m_{i}^{j_2} (\Delta_n, \trueparam);  \\[0.2cm]  
  \mathscr{K}_3 ( \theta )  
  & = \tfrac{1}{n} \sum_{i = 1}^n 
  R (\sqrt{\Delta_n}, \sample{i-1}, \theta), 
\end{align*}
for $R, R_j, R_{j_1 j_2} \in \mathcal{S}$. From Lemmas \ref{lemma:ergodic_thm} and \ref{lemma:canonical_conv}, we immediately have that $ \mathscr{K}_i ( \truebeta_S, \beta_R, \sigma )  \probconv \mathbf{0}_{N_{\beta_{S_2}}}, \; 1 \le i \le 3$, as \limit, thus (\ref{eq:K}) holds. 
\\ 

\noindent  
(\textit{Proof of (\ref{eq:step2_B})}). 
We define 
\begin{align*}
\mathcal{Q} (\theta) = M_{\beta_{S}, n} \, \partial^2_{\beta_S} 
\ell_n (\theta) \, M_{\beta_S, n} 
=
\begin{bmatrix} 
\mathcal{Q}_{S_1 S_1}  ( \theta )  
& \mathcal{Q}_{S_1 S_2}  ( \theta )  \\[0.2cm]
\mathcal{Q}_{S_2 S_1}  ( \theta )  
& \mathcal{Q}_{S_2 S_2}  ( \theta )  
\end{bmatrix}, 
\end{align*} 
where we have set: 
\begin{align*}
\mathcal{Q}_{S_1 S_1}  ( \theta )  
& = \tfrac{\Delta_n^3}{n} \partial^2_{\beta_{S_1}}  \ell_n (\theta), \quad 
\mathcal{Q}_{S_1 S_2} ( \theta )  
= \tfrac{\Delta_n^2}{n} \partial_{\beta_{S_1}} \partial_{\beta_{S_2}}^\top \ell_n (\theta); \\[0.2cm] 
\mathcal{Q}_{S_2 S_1}  ( \theta ) 
& =  \mathcal{Q}_{S_1 S_2} ( \theta )^\top , \quad \mathcal{Q}_{S_2 S_2} ( \theta )  
= \tfrac{\Delta_n}{n} \partial^2_{\beta_{S_2}}  
\ell_n (\theta),  
\end{align*} 
for $\theta = (\beta_S, \beta_R, \sigma) \in \Theta$. From the proof of Lemma \ref{lemma:step1_B} in Appendix \ref{appendix:pf_step1_B} and the consistency of the estimator $\hat{\beta}_{S,n}$, we have that if \limit, then 
\begin{align} \label{eq:Q_S1S1}
  \mathcal{Q}_{S_1 S_1}
\bigl( 
 \truebeta_S + \lambda ( \hat{\beta}_{S, n} - \truebeta_S), \beta_R, \sigma 
\bigr) \probconv 2 \int 
\partial_{\beta_{S_1}} \bigl( V_{S_1, 0} (x_S, \truebeta_{S_1}) \bigr)^\top
\Lambda_{S_1 S_1} \bigl(x, (\truebeta_S, \sigma) \bigr) 
\partial_{\beta_{S_1}}^\top V_{S_1, 0} (x_S, \truebeta_{S_1}) \, \truedist (dx),
\end{align}
uniformly in $(\beta_R, \sigma) \in \Theta_{\beta_R} \times \Theta_\sigma$ and $\lambda \in [0,1]$. 
We will check the convergence of the two matrices $\mathcal{Q}_{S_1 S_2} (\theta)$ and $\mathcal{Q}_{S_2 S_2} (\theta)$. 
\\ 

\noindent We have 
$\textstyle{ 
\mathcal{Q}_{S_1 S_2} (\theta) 
=\sum_{1 \le k \le 4}  \mathcal{Q}_{S_1 S_2, k} (\theta) }$, where we set for $1 \le j_1 \le N_{\beta_{S_1}}$, $1 \le j_2 \le N_{\beta_{S_2}}$, 
\begin{align*}
\bigl[ \mathcal{Q}_{S_1 S_2, 1} (\theta)  \bigr]_{j_1 j_2} 
& = \tfrac{2}{n} \sum_{i = 1}^n 
   \bigl( \partial^{\beta_{S_1}}_{j_1} 
   V_{S_1, 0} (\sample{S, i-1}, \beta_{S_1}) \bigr)^\top 
   \times 
   \biggl\{ 
   \Lambda_{S_1 S_2} \bigl(\sample{i-1}, \theta \bigr)
   \partial^{\beta_{S_2}}_{j_2} V_{S_2, 0} ( \sample{i-1}, \beta_{S_2})  \\[0.1cm]
& \qquad \qquad \qquad \qquad  
  + \tfrac{1}{2} 
   \Lambda_{S_1 S_1} \bigl(\sample{i-1}, \theta \bigr)
   \partial^{\beta_{S_2}}_{j_2} \mathcal{L} V_{S_1, 0} ( \sample{i-1}, \beta_{S}) 
   \biggr\}; \\[0.3cm] 
 \bigl[ \mathcal{Q}_{S_1 S_2, 2} (\theta) \bigr]_{j_1 j_2} 
 & = \tfrac{1}{n} \sum_{i = 1}^n 
  \sum_{\substack{1 \le k_1, k_2 \le N_{S_1}}} 
  R_{j_1 j_2}^{k_1 k_2} (1, \sample{i-1}, \theta) 
  \eta_{S_1, i-1}^{k_1} (\Delta_n, \beta_{S_1})
  \eta_{S_1}^{k_2} (\sample{i-1}, \beta_{S_1})
  ; \\[0.2cm]  
 \bigl[ \mathcal{Q}_{S_1 S_2, 3} (\theta) \bigr]_{j_1 j_2} 
 & = \tfrac{1}{n} \sum_{i = 1}^n 
  \sum_{\substack{1 \le k_1 \le N_{S_1} \\ 1 \le k_2 \le N_{S} }}  
  \widetilde{R}_{j_1 j_2}^{k_1 k_2} (1, \sample{i-1}, \theta) 
  \eta_{S_1}^{k_1} (\sample{i-1}, \beta_{S_1}) b^{k_2} (\sample{i-1}, \beta_S)
  ; \\[0.2cm] 
  \bigl[ \mathcal{Q}_{S_1 S_2, 4} (\theta) \bigr]_{j_1 j_2} 
 & = \tfrac{1}{n} \sum_{i = 1}^n \, 
 \biggl\{  \sum_{1 \le k_1, k_2 \le N} 
 \bar{R}_{j_1j_2}^{k_1 k_2} (\Delta_n, \sample{i-1}, \theta) m_{i}^{k_1} (\Delta_n, \trueparam) 
 m_{i}^{k_2} (\Delta_n, \trueparam) \\[0.2cm]
 & \qquad \qquad \quad
+ \sum_{1 \le k \le N} \bar{R}_{j_1 j_2}^k (\Delta_n, \sample{i-1}, \theta) 
 m_{i}^k (\Delta_n, \trueparam ) 
 + \bar{R} (\Delta_n, \sample{i-1}, \theta) 
 \biggr\}, 
\end{align*}
for some functions $R_{j_1 j_2}^{k_1 k_2}, \, \widetilde{R}_{j_1 j_2}^{k_1 k_2}, \, \bar{R}_{j_1 j_2}^{k_1 k_2}, \, \bar{R}_{j_1 j_2}^{k}, \,  \bar{R} \in \mathcal{S}$, and the function $b: \mathbb{R}^N \times \Theta_{\beta_S} \to \mathbb{R}^{N_S}$ is defined as:
\begin{align*}
b (x, \beta_S ) =
\begin{bmatrix}  
\tfrac{1}{2} \mathcal{L} V_{S_1, 0} (x, \beta_{S})
- \tfrac{1}{2} \mathcal{L} V_{S_1, 0} (x, \truebeta_{S}) 
\\[0.2cm]
V_{S_2, 0} (x, \beta_{S_2}) - V_{S_2, 0} (x, \truebeta_{S_2}) 
\end{bmatrix}, 
\end{align*} 
Notice that for any $\theta \in \Theta$, 
\begin{align*}
[\mathcal{Q}_{S_1 S_2, 1}(\theta)]_{j_1 j_2} = 0, \, \quad 1 \le j_1 \le N_{\beta_{S_1}}, \quad 1 \le j_2 \le N_{\beta_{S_2}}, 
\end{align*}
because it follows from Lemma \ref{lemma:matrix} that
\begin{align*}
\Lambda_{S_1 S_2} \bigl(x, \theta \bigr)
\partial^{\beta_{S_2}}_{j_2} V_{S_2, 0} ( x, \beta_{S_2}) 
+ \tfrac{1}{2} 
\Lambda_{S_1 S_1} \bigl(x, \theta \bigr)
\partial^{\beta_{S_2}}_{j_2} \mathcal{L} V_{S_1, 0} (x, \beta_{S}) 
= \Phi (x, \theta) \, \partial^{\beta_{S_2}}_{j_2} V_{S_2, 0} ( x, \beta_{S_2})  = \mathbf{0}_{N_{S_1}}, 
\end{align*}
for any $x \in \mathbb{R}^N$, $1 \le j_2 \le N_{\beta_{S_2}}$, 
where $\Phi (x, \theta)$ is defined as in (\ref{eq:phi_0}). 
From Lemmas \ref{lemma:ergodic_thm}, \ref{lemma:canonical_conv} and the consistency of $\hat{\beta}_{S, n}$ with 
convergence (\ref{eq:step1}), we obtain: 
\begin{align*}
%
\bigl[  \mathcal{Q}_{S_1 S_2, k} \bigl( \truebeta_S + \lambda (\hat{\beta}_{S, n} - \truebeta_S), \beta_R, \sigma \bigr) 
 \bigr]_{j_1 j_2} \probconv 0, \qquad  2 \le k \le 4,
\end{align*} 
uniformly in $(\beta_R, \sigma) \in \Theta_{\beta_R} \times \Theta_{\sigma}$ and $\lambda \in [0,1]$. 

Finally, we consider the term $\mathcal{Q}_{S_2 S_2} (\theta)$. It holds that 
$\textstyle{\mathcal{Q}_{S_2 S_2} (\theta) = \sum_{1 \le k \le 4} \mathcal{Q}_{S_2 S_2, k}(\theta)}$, where we set, for $1 \le j_1, j_2 \le N_{\beta_{S_2}}$, 
\begin{align*}    
 \bigl[ \mathcal{Q}_{S_2 S_2, 1}(\theta) \bigr]_{j_1 j_2} 
 &  = \tfrac{1}{n} \sum_{i = 1}^n \biggl\{
     \sum_{1 \le k  \le  N_{S_1}}   
     {R}_{j_1 j_2}^{k} (1, \sample{i-1}, \theta)) 
     \eta_{S_1, i-1}^{k} (\Delta_n, \beta_{S_1})  \\[0.2cm] 
     & \quad \quad 
     + \sum_{1 \le k_1, k_2  \le  N_{S_1}} R_{j_1 j_2}^{k_1 k_2} (1, \sample{i-1}, \theta)   
     \eta_{S_1, i-1}^{k_1} (\Delta_n, \beta_{S_1}) 
     \eta_{S_1, i-1}^{k_2} (\Delta_n, \beta_{S_1})
     \\[0.3cm] 
     & \quad \quad  + \sum_{\substack{1 \le k_1 \le N_{S_1} \\ 1 \le k_2 \le N}} \widetilde{R}_{j_1 j_2}^{k_1 k_2} (1, \sample{i-1}, \theta) 
     \eta_{S_1, i-1}^{k_1} (\sqrt{\Delta_n}, \beta_{S_1}) 
     m_{i}^{k_2} (\Delta_n, \trueparam) \biggr\};
     \\[0.2cm]
\bigl[\mathcal{Q}_{S_2 S_2, 2} (\theta) \bigr]_{j_1 j_2} 
& = \tfrac{2}{n} \sum_{i = 1}^n  \Bigl\{
\tfrac{1}{4} \bigl( \partial_{j_1}^{\beta_{S_2}} \mathcal{L} V_{S_1, 0} (\sample{i-1}, \theta) 
\bigr)^\top 
\Lambda_{S_1 S_1} (\sample{i-1}, \theta) 
\partial_{j_2}^{\beta_{S_2}} \mathcal{L} V_{S_1, 0} (\sample{i-1}, \theta) \\
& \qquad \qquad + 
\tfrac{1}{2} 
\bigl( \partial_{j_1}^{\beta_{S_2}} \mathcal{L} V_{S_1, 0} (\sample{i-1}, \theta)  \bigr)^\top
\Lambda_{S_1 S_2} (\sample{i-1}, \theta) 
\, \partial_{j_2}^{\beta_{S_2}} V_{S_2, 0} (\sample{i-1}, \beta_{S_2}) \\[0.2cm] 
& \qquad \qquad + 
\tfrac{1}{2} 
\bigl( \partial_{j_1}^{\beta_{S_2}} V_{S_2, 0} (\sample{i-1}, \beta_{S_2})
\bigr)^\top
\Lambda_{S_2 S_1} (\sample{i-1}, \theta) 
\partial_{j_2}^{\beta_{S_2}} \,  
\mathcal{L} V_{S_1, 0} (\sample{i-1}, \theta) \\[0.2cm]  
& \qquad \qquad + 
\bigl( \partial^{\beta_{S_2}}_{j_1} V_{S_2, 0} (\sample{i-1}, \beta_{S_2}) \bigr)^\top \Lambda_{S_2 S_2} 
\bigl( \sample{i-1}, \theta \bigr) 
\, \partial^{\beta_{S_2}}_{j_2} V_{S_2, 0} (\sample{i-1}, \beta_{S_2}) 
\Bigr\};  \\[0.1cm] 
\bigl[\mathcal{Q}_{S_2 S_2, 3} (\theta) \bigr]_{j_1 j_2} 
& = \tfrac{1}{n} \sum_{i = 1}^n
\biggl\{  \sum_{1 \le k_1, k_2 \le N_{S_2}} 
\overline{R}_{j_1 j_2}^{k_1 k_2} (1, \sample{i-1}, \theta) 
\eta_{S_2}^{k_1} (\sample{i-1}, \beta_{S_2}) 
\eta_{S_2}^{k_2} (\sample{i-1}, \beta_{S_2}) 
\\[0.1cm] 
& \qquad \qquad + 
\sum_{1 \le k \le N_{S_2}} \widetilde{R}_{j_1 j_2}^{k} (1, \sample{i-1}, \theta))
\eta_{S_2}^{k} (\sample{i-1}, \beta_{S_2}) \biggr\};
\\[0.3cm]
\bigl[ \mathcal{Q}_{S_2 S_2, 4}(\theta) \bigr]_{j_1 j_2}
& = 
\tfrac{1}{n} \sum_{i = 1}^n \biggl\{ 
 \sum_{1 \le k_1, k_2 \le N} 
 \hat{R}_{j_1 j_2}^{k_1 k_2} (\Delta_n, \sample{i-1}, \theta) 
 m_{i}^{k_1} (\Delta_n, \trueparam)
 m_{i}^{k_2} (\Delta_n, \trueparam) 
 \\[0.2cm]
 & + \sum_{1 \le k \le N} \overline{R}_{j_1 j_2}^{k} (\sqrt{\Delta_n}, \sample{i-1}, \theta) m_{i}^k (\Delta_n, \trueparam) 
 + R_{j_1 j_2} (\sqrt{\Delta_n}, \sample{i-1}, \theta) \biggr\},  
\end{align*}
for some functions 
$R_{j_1 j_2}^k, \, {R}_{j_1 j_2}^{k_1 k_2}, \, \widetilde{R}_{j_1 j_2}^{k_1 k_2}, \, \overline{R}_{j_1 j_2, }^{k_1 k_2}, \, \widetilde{R}_{j_1 j_2}^k, \, \hat{R}_{j_1 j_2, }^{k_1 k_2}, \, \overline{R}_{j_1 j_2}^{k}, \, R_{j_1 j_2}  \in \mathcal{S}$. 
Note that due to Lemma~\ref{lemma:matrix}, 
\begin{gather*}
[Q_{S_2 S_2, 2} (\theta)]_{j_1 j_2}  
= \tfrac{2}{n} \sum_{i = 1}^n 
\tfrac{1}{2} 
\bigl( \partial_{j_1}^{\beta_{S_2}} V_{S_2, 0} (\sample{i-1}, \beta_{S_2})
\bigr)^\top
\Lambda_{S_2 S_1} (\sample{i-1}, \theta) 
\, \partial_{j_2}^{\beta_{S_2}} \,  
\mathcal{L} V_{S_1, 0} (\sample{i-1}, \theta) \\[0.1cm]
+ \tfrac{2}{n} \sum_{i = 1}^n 
\bigl( \partial^{\beta_{S_2}}_{j_1} V_{S_2, 0} (\sample{i-1}, \beta_{S_2}) \bigr)^\top \Lambda_{S_2 S_2} \bigl( \sample{i-1}, \theta \bigr) \, \partial^{\beta_{S_2}}_{j_2} V_{S_2, 0} (\sample{i-1}, \beta_{S_2}).  
\end{gather*}
We obtain from Lemma \ref{lemma:ergodic_thm}, (\ref{eq:eta_conv}) and consistency of $\hat{\beta}_{S, n}$ that if \limit, then 
\begin{align} 
& \Bigl[\mathcal{Q}_{S_2 S_2, 2} \bigl( \truebeta_{S} + \lambda (\hat{\beta}_{S, n} - \truebeta_S), \beta_R, \sigma  \bigr) \Bigr]_{j_1 j_2}  \nonumber \\[0.2cm]
& \probconv 2 \int \tfrac{1}{2} 
\bigl( \partial_{j_1}^{\beta_{S_2}} V_{S_2, 0} (x, \truebeta_{S_2}) \bigr)^\top
 \Lambda_{S_2 S_1} \bigl(x, (\truebeta_S, \beta_R, \sigma) \bigr) \partial_{j_2}^{\beta_{S_2}} \mathcal{L} V_{S_1, 0} (x, (\truebeta_S, \beta_R, \sigma) )   \truedist (dx)  \nonumber \\[0.1cm]
& \qquad  \qquad   + 2 \int 
\bigl( 
\partial_{j_1}^{\beta_{S_2}} V_{S_2, 0} (x, \truebeta_{S_2}) \bigr)^\top
 \Lambda_{S_2 S_2} \bigl(x, (\truebeta_S, \beta_R, \sigma) \bigr) \partial_{j_2}^{\beta_{S_2}} V_{S_2, 0} (x, \truebeta_{S_2})   \truedist (dx); 
 \label{eq:Q_S2S2} \\[0.1cm] 
 & \Bigl[\mathcal{Q}_{S_2 S_2, k} \bigl( \truebeta_{S} + \lambda (\hat{\beta}_{S,  n} - \truebeta_S), \beta_R, \sigma  \bigr) \Bigr]_{j_1 j_2} 
 \probconv 0, \qquad  k = 1, 3, 4. \nonumber 
\end{align} 
The proof is complete by applying Lemma \ref{lemma:lambda} to (\ref{eq:Q_S1S1}) and (\ref{eq:Q_S2S2}). 
\subsection{Proof of Lemma \ref{lemma:step3_1}} 
\label{appendix:pf_step3_1}
It holds that $\textstyle{\tfrac{1}{n} \ell_n (\theta) = \sum_{1 \le k \le 4} \mathscr{T}_k (\theta)}$, $\theta \in \Theta$, where we have set: 
\begin{align*}
\mathscr{T}_1 (\theta) 
& = \tfrac{1}{n} \sum_{i = 1}^n 
\biggl\{ \sum_{1 \le j_1, j_2 \le N_{S_1}} 
R_{j_1 j_2}^{S_1 S_1} (1, \sample{i-1}, \theta) 
\eta_{S_1, i-1}^{j_1} (\sqrt{\Delta_n^3}, \beta_{S_1})
\eta_{S_1, i-1}^{j_2} (\sqrt{\Delta_n^3}, \beta_{S_1})
\\[0.2cm]
& \qquad  + \sum_{\substack{ 1 \le j_1  \le N_{S_1} \\ 1 \le j_2 \le N_{S_2}} }  R_{j_1 j_2}^{S_1 S_2} (1, \sample{i-1}, \theta) 
\eta_{S_1, i-1}^{j_1} (\sqrt{\Delta_n^3}, \beta_{S_1})
\eta_{S_2, i-1}^{j_2} (\sqrt{\Delta_n}, \beta_{S_2}) \\[0.2cm] 
& \qquad + \sum_{1 \le j_1, j_2 \le N_{S_2}}  
R_{j_1 j_2}^{S_2S_2} (1, \sample{i-1}, \theta) 
\eta_{S_2, i-1}^{j_1} (\sqrt{\Delta_n}, \beta_{S_2})
\eta_{S_2, i-1}^{j_2} (\sqrt{\Delta_n}, \beta_{S_2}) 
\biggr\}; \\[0.2cm]
\mathscr{T}_2 (\theta) 
& = \tfrac{1}{n} \sum_{i = 1}^n 
\biggl\{ 
\sum_{\substack{1 \le j_1 \le N_{S_1} \\ 1 \le j_2 \le N}} {R}_{j_1 j_2}^{S_1} (1, \sample{i-1}, \theta) 
\eta_{S_1, i-1}^{j_1} (\sqrt{\Delta_n^3}, \beta_{S_1} ) 
m_{i}^{j_2} (\Delta_n, \trueparam)  \\[0.2cm]
& \qquad + \sum_{\substack{1 \le j_1 \le N_{S_2} \\ 1 \le j_2 \le N}}{R}_{j_1 j_2}^{S_2} (1, \sample{i-1}, \theta) 
\eta_{S_2, i-1}^{j_1} (\sqrt{\Delta_n}, \beta_{S_2} ) m_{i}^{j_2} (\Delta_n, \trueparam)  \\[0.2cm]
& \qquad + \sum_{1 \le j \le N_{S_1}} {R}_{j}^{S_1} (\sqrt{\Delta_n}, \sample{i-1}, \theta) \eta_{S_1, i-1}^{j} (\sqrt{\Delta_n^3}, \beta_{S_1}) + \sum_{1 \le j \le N_{S_2}}
{R}_{j}^{S_2} (\sqrt{\Delta_n}, \sample{i-1}, \theta) \eta_{S_2, i-1}^{j} (\sqrt{\Delta_n}, \beta_{S_2}) \biggr\}; \\[0.2cm]  
\mathscr{T}_3 (\theta)  
 & = \tfrac{1}{n} \sum_{i = 1}^n
 \biggl\{ 
m_{i}(\Delta_n, \trueparam)^\top
\Lambda \bigl(\sample{i-1}, (\beta_S, \sigma) \bigr) 
m_{i} (\Delta_n, \trueparam) + \log \bigl| \Sigma \bigl( \sample{i-1}, (\beta_S, \sigma) \bigr) \bigr| 
\biggr\};  \\[0.2cm]
\mathscr{T}_4 (\theta)
 & = \tfrac{1}{n} \sum_{i = 1}^n \biggl\{ 
 \sum_{1 \le j \le N} R_j (\sqrt{\Delta_n}, \sample{i-1}, \theta) m_{i}^{j} (\Delta_n, \trueparam) + R (\Delta_n, \sample{i-1}, \theta)  \biggr\}, 
\end{align*} 
where $R_{j_1 j_2}^{S_1 S_1}, \, R_{j_1 j_2}^{S_1 S_2}, \, R_{j_1 j_2}^{S_2 S_2}, \, R_{j_1 j_2}^{S_1}, \, R_{j_1 j_2}^{S_2}, \, R_{j}^{S_1}, \, R_{j}^{S_2}, \, R_{j_1 j_2} \, R_{j}, R  \in \mathcal{S}$. 
From Lemmas \ref{lemma:ergodic_thm}, ~\ref{lemma:canonical_conv} and the convergences (\ref{eq:step1}) and (\ref{eq:step2}), we have: 
\begin{align*} 
& \mathscr{T}_k \bigl( \hat{\beta}_{S, n}, \beta_R, \sigma \bigr) \probconv 0, \qquad  k = 1,2,4; \\[0.1cm]
& \mathscr{T}_3 \bigl( \hat{\beta}_{S, n}, \beta_R, \sigma \bigr) \probconv 
\int \Bigl\{ 
\mrm{tr} \bigl( 
\Lambda \bigl( x, (\truebeta_S, \sigma)  \bigr)
\Sigma \bigl( x, (\truebeta_S, \truesigma) \bigr) \bigr)
+ \log 
\bigl| \Sigma \bigl(x, (\truebeta_S, \sigma) \bigr) 
\bigr| 
\Bigr\} \, \truedist (dx), 
\end{align*} 
as \limit, uniformly in $(\beta_R, \sigma) \in \Theta_{\beta_R} \times \Theta_\sigma$. The proof is now complete. 
\subsection{Proof of Lemma \ref{lemma:step3_2}} 
\label{appendix:pf_step3_2}
The matrix-valued function $\Sigma (x, \theta)$ and its inverse $\Lambda (x, \theta)$ depend on $\beta_S = (\beta_{S_1}, \beta_{S_2}) \in \Theta_{\beta_S}$ and $\sigma \in \Theta_\sigma$ but not on $\beta_R \in \Theta_{\beta_R}$ in terms of the parameter $\theta \in \Theta$. We then define 
\begin{align*}
\mathscr{L} (\theta) = \tfrac{1}{n \Delta_n} \ell_n (\theta) 
 -  \tfrac{1}{n \Delta_n} \ell_n \bigl( \beta_S, \truebeta_R, \sigma \bigr), \quad \theta = (\beta_S, \beta_R, \sigma) \in \Theta,  
\end{align*} 
and $\mathscr{L}(\theta)$ is expressed as 
$\textstyle{ \mathscr{L}(\theta) = \sum_{1 \le k \le 3} \mathscr{U}_k (\theta)}$, where we have set: 
\begin{align*}
\mathscr{U}_{1} (\theta) 
& = \tfrac{1}{n \Delta_n} 
\sum_{i = 1}^n 
\biggl\{ 
\Bigl( m_{i} (\Delta_n, (\beta_S, \truebeta_R, \sigma) ) 
- m_{i} (\Delta_n, \trueparam) \Bigr)^\top 
\Lambda (\sample{i-1}, (\beta_S, \sigma)) \Bigl( 
m_{i} (\Delta_n, \theta) - m_{i} (\Delta_n, (\beta_S, \truebeta_R, \sigma) )  \Bigr) 
\biggr\}; \\[0.2cm]
\mathscr{U}_2 (\theta) 
& = \tfrac{1}{n \Delta_n} \sum_{i = 1}^n \biggl\{ 
\Bigl( m_{i} (\Delta_n, \theta) - m_{i} (\Delta_n, \trueparam) \Bigr)^\top  
\Lambda (\sample{i-1}, (\beta_S, \sigma))
\Bigl( m_{i} (\Delta_n, \theta) 
- m_{i} (\Delta_n, (\beta_S, \truebeta_R, \sigma) )
\Bigr)   
\biggr\}; \\[0.2cm]
\mathscr{U}_3 (\theta) 
& = \tfrac{2}{n \Delta_n} \sum_{i = 1}^n
m_{i} (\Delta_n, \trueparam)^\top
\Lambda (\sample{i-1}, (\beta_S, \sigma))
\bigl( m_{i} (\Delta_n, \theta) - m_{i} (\Delta_n, (\beta_S, \truebeta_R, \sigma) ) \bigr). 
\end{align*} 
We will derive the limit of the terms $\mathscr{U}_k (\theta)$, $1 \le k \le 3$ evaluated at 
$ \theta = ( \hat{\beta}_{S, n}, \,  \beta_R, \, 
 \hat{\sigma}_{n} )$ by utilising 
 Lemma \ref{lemma:matrix_2}: 
%
%
\noindent We first consider the term $\mathscr{U}_1 (\theta)$. Making use of Lemma \ref{lemma:matrix_2}, we have
\begin{align} 
 \begin{aligned}  \label{eq:U_1}
 \mathscr{U}_{1} (\theta) 
 & =  \tfrac{1}{n} \sum_{i = 1}^n 
\biggl\{  \sum_{1 \le j \le N_{S_1}} 
\eta_{S_1, i-1}^j (\sqrt{\Delta_n^3}, \beta_{S_1}) 
R^j_{S_1} (\sqrt{\Delta_n}, \sample{i-1}, \theta) \\[0.2cm]
& \qquad + \sum_{1 \le j \le N_{S_2}} 
\eta_{S_2, i-1}^j (\sqrt{\Delta_n}, \beta_{S_2}) 
R^j_{S_2} (\sqrt{\Delta_n}, \sample{i-1}, \theta)
+ R (\sqrt{\Delta_n}, \sample{i-1}, \theta) \biggr\},
 \end{aligned} 
\end{align} 
where $R_{S_1}^j, \, R_{S_2}^j, \, R \in \mathcal{S}$. From Lemma \ref{lemma:ergodic_thm} and the limits (\ref{eq:step1}), (\ref{eq:step2}), we obtain that if \limit, 
\begin{align*} 
\mathscr{U}_1 \big( \hat{\beta}_{S, n}, \, 
\beta_R, \, 
\hat{\sigma}_{n} \bigr) \probconv 0, 
\end{align*}
uniformly in $\beta_R \in \Theta_{\beta_R}$. For the term $\mathscr{U}_2 (\theta)$, again Lemma \ref{lemma:matrix_2} yields
\begin{align*}
\mathscr{U}_2 (\theta) 
 & = \tfrac{1}{n} \sum_{i = 1}^n  
 \eta_{R} (\sample{i-1}, \beta_R)^\top 
 a_R^{-1} (\sample{i-1}, \sigma)
 \eta_{R} (\sample{i-1}, \beta_R) 
 + \widetilde{\mathscr{U}}_2 (\theta),
\end{align*} 
where $\widetilde{\mathscr{U}}_2 (\theta)$ is given in the form of the right-hand side of formula (\ref{eq:U_1}). We then obtain 
\begin{align*}
\mathscr{U}_2 \big( \hat{\beta}_{S,  n}, \, 
\beta_R, \, 
\hat{\sigma}_{n} \bigr)  
\probconv \int \eta_{R} (x, \beta_R)^\top \, 
a_R^{-1} (x,  \truesigma ) \, 
\eta_{R} (x, \beta_R)  \truedist(dx), 
\end{align*}
as \limit, uniformly in $\beta_R \in \Theta_{\beta_R}$. 
For the third term $\mathscr{U}_3 (\theta)$, it follows from Lemma \ref{lemma:matrix_2} that 
\begin{align*}
 \mathscr{U}_3 (\theta) = 
 & \tfrac{1}{n \sqrt{\Delta_n}} 
  \sum_{i = 1}^n \sum_{N_S +1  \le j \le N} 
  m_{i}^j (\Delta_n, \trueparam) R^j (1, \sample{i-1}, \theta)  +  \tfrac{1}{n} 
  \sum_{i = 1}^n \sum_{1 \le j \le N} 
   m_{i}^j (\Delta_n, \trueparam) 
  \widetilde{R}^j (\sqrt{\Delta_n}, \sample{i-1}, \theta), 
\end{align*}
where $R^j, \, \widetilde{R}^j \in \mathcal{S}$. From Lemma \ref{lemma:canonical_conv}, we have that, if \limit, then 
\begin{align*}
 \mathscr{U}_3 (\hat{\beta}_{S,  n}, \beta_R, \hat{\sigma}_{ n})
 \probconv 0, 
\end{align*} 
uniformly in  $\beta_R \in \Theta_{\beta_R}$. The proof is now complete. 
\subsection{Proof of Lemma \ref{lemma:slln}}
\label{appendix:pf_slln} 
Recall $N_\theta =  N_{\beta} + N_{\sigma}$, where $N_\beta = N_{\beta_S} + N_{\beta_R}$ with $N_{\beta_S} = N_{\beta_{S_1}} + N_{\beta_{S_2}}$. In this section we make use of the notation $\beta = (\beta_S, \beta_R) \in \Theta_{\beta_S} \times \Theta_{\beta_R}$ with $\beta_S = (\beta_{S_1}, \beta_{S_2})$. We note again that the matrices $\Sigma (x, \theta)$ and $\Lambda (x, \theta)$ do not depend on the parameter $\beta_R$. Since we have seen the convergences of the matrices $\mathcal{Q}_{S_1 S_1} (\theta), \, \mathcal{Q}_{S_1 S_2} (\theta), \, \mathcal{Q}_{S_2 S_2} (\theta)$ evaluated at $\theta = \trueparam + \lambda ( \hat{\theta}_{n} - \trueparam), \, \lambda \in [0,1]$ in Appendix \ref{appendix:pf_step_4}, we prove that as \limit, the following convergences hold uniformly in $\lambda \in [0,1]$:
\begin{enumerate}
\item[(a)] For $1 \le i \le N_{\beta}$, $N_{\beta_S} + 1  \le j \le N_{\beta}$, 
\begin{align}
& \bigl[ \mathscr{I}_{n} \bigl( \trueparam + \lambda (\hat{\theta}_{ n} - \trueparam) \bigr) \bigr]_{ij}
\label{eq:slln_conv1}  \\[0.2cm] 
& \probconv 
\begin{cases}
2 \int \bigl( \partial^{\beta}_i V_{R, 0} (x, \truebeta_R) \bigr)^\top\, 
a_R^{-1} (x, \truesigma) \, 
\partial^{\beta}_j V_{R, 0} (x, \truebeta_R) \,  \truedist (dx), 
  & N_{\beta_S} +1 \le i, j \le N_{\theta};  \\[0.3cm]
0, &  (\mrm{otherwise}). 
\end{cases} 
\nonumber 
\end{align}
\item[(b)] For $1 \le i \le N_{\beta}$, $N_{\beta} + 1  \le j \le N_{\theta}$, 
\begin{align}
\bigl[ \mathscr{I}_{n} \bigl( \trueparam + \lambda (\hat{\theta}_{n} - \trueparam) \bigr) \bigr]_{ij}
\probconv 0.   \label{eq:slln_conv2}
\end{align}
 \item[(c)] For $N_{\beta} + 1  \le i, j \le N_{\theta}$, 
 \begin{align}
  \bigl[ \mathscr{I}_{n} \bigl( \trueparam + \lambda (\hat{\theta}_{n} - \trueparam) \bigr) \bigr]_{ij}   
\probconv  \int  \mathrm{tr}  \bigl(
\partial^{\sigma}_i \Sigma  (x, \trueparam ) \, 
\Lambda  (x, \trueparam )  
\partial^{\sigma}_j \Sigma  (x, \trueparam ) \,
\Lambda  (x, \trueparam ) \bigr)  
\truedist (dx),   \label{eq:slln_conv3}  
\end{align}
\end{enumerate}
\subsubsection{Proof of (\ref{eq:slln_conv1})} 
We consider the following three cases separately: 
\begin{enumerate}
\item[(a1)] $1 \le i \le N_{\beta_{S_1}}, \; N_{\beta_S} + 1 \le j \le N_\beta$;
\item[(a2)] $N_{\beta_{S_1}} + 1 \le i \le N_{\beta_S}, \; N_{\beta_S} + 1 \le j \le N_\beta$;
\item[(a3)] $N_{\beta_S} + 1 \le i, j \le N_{\beta}$. 
\end{enumerate}
In  case (a1), we have for $\theta = (\beta_{S_1}, \beta_{S_2}, \beta_R, \sigma ) \in \Theta$,   
\begin{align}  
\bigl[ \mathscr{I}_{n} \bigl( \theta \bigr) \bigr]_{ij} 
& = \tfrac{1}{n} \sum_{m=1}^n H_{ij} (\sample{m-1}, \theta) 
+ \tfrac{1}{n} \sum_{m = 1}^n R (\Delta_n, \sample{m-1}, \theta) \label{eq:case_a1}  \\ 
& + \tfrac{1}{n}  \sum_{m = 1}^n \sum_{ 1 \le k \le N_{S_1} }
R_{k} (1, \sample{m-1}, \theta) 
\bigl( V_{S_1, 0}^{k} (\sample{S, m-1}, \truebeta_{S_1}) - V_{S_1, 0}^{k} (\sample{S, m-1}, \beta_{S_1} ) \bigr),   \nonumber 
\end{align}
for some $R_{k}, R \in \mathcal{S}, \, 1 \le k \le N_{S_1}$, where we have defined $H_{ij} (x, \theta), \, (x, \theta) \in \mathbb{R}^N \times \Theta$ as: 
\begin{align*}
H_{ij} (x, \theta) 
& =
\tfrac{1}{3} 
\bigl( \partial^{\beta}_i V_{S_1, 0} (x_{S_1}, \beta_{S_1}) \bigr)^\top 
\Lambda_{S_1 S_1} ( x, \theta ) 
\, \partial^{\beta}_{j}  \mathcal{L}^2 V_{S_1, 0} (x, \theta) 
+ 
\bigl( \partial^{\beta}_{i} V_{S_1, 0} (x_{S_1}, \beta_{S_1}) \bigr)^\top \Lambda_{S_1 S_2} (x, \theta) 
\, \partial^{\beta}_{j} \mathcal{L} V_{S_2, 0} (x, \theta)  \\[0.2cm] 
& \quad  + 2
\bigl( \partial^{\beta}_i V_{S_1, 0} (x_{S_1}, \beta_{S_1}) \bigr)^\top 
\Lambda_{S_1 R} (x, \theta) 
\, \partial^{\beta}_{j} V_{R, 0} (x, \beta_R). 
\end{align*}
Noticing that for $N_{\beta_S} +  1  \le j  \le N_{\beta}$,  
\begin{gather*}
\partial^\beta_j \mathcal{L} V_{S_2, 0} (x, \theta) 
= \partial_{x_R}^\top V_{S_2, 0} (x, \beta_{S_2}) \partial_j^\beta V_{R, 0} (x, \beta_R), \quad 
\partial^\beta_j \mathcal{L}^2 V_{S_1, 0} (x, \theta) 
= \partial_{x_{S_2}}^\top V_{S_1, 0} (x_S, \beta_{S_1})  
\partial_{x_R}^\top V_{S_2, 0} (x, \beta_{S_2}) 
\partial_j^\beta V_{R, 0} (x, \beta_R),
\end{gather*}
we have $H_{ij} (x, \theta) = 0$ for any $(x, \theta) \in \mathbb{R}^N \times \Theta$ since 
\begin{align*} 
H_{ij} (x, \theta)  
= 2 \bigl( \partial^{\beta}_{i} V_{S_1, 0} (x, \beta_{S_1}) \bigr)^\top  
\widetilde{H}_{ij} (x, \theta) \, 
\partial^{\beta}_j V_{R, 0} (x, \beta_R),   
\end{align*}
where 
\begin{align*}
\widetilde{H}_{ij} (x, \theta)
& = 
\bigl\{ 
\Lambda_{S_1 S_1} (x, \theta) \Sigma_{S_1 R} (x, \theta) 
+ \Lambda_{S_1 S_2} (x, \theta) \Sigma_{S_2 R} (x, \theta)  
+ \Lambda_{S_1 R} (x, \theta ) \Sigma_{R R} (x, \sigma) 
\bigr\} a_R^{-1} (x, \sigma) = \mathbf{0}_{N_{S_1} \times N_{R}}. 
\end{align*}
Thus, due to the consistency of the estimator and Lemma \ref{lemma:ergodic_thm}, we immediately obtain from (\ref{eq:case_a1}) that if \limit, then 
\begin{align*}  
[\mathscr{I}_n \bigl( \trueparam + \lambda (\hat{\theta}_n - \trueparam)  \bigr) ]_{ij} \probconv 0, 
\end{align*} 
uniformly in $\lambda \in [0,1]$ for $1 \le i \le N_{\beta_{S_1}}, \; N_{\beta_S} + 1 \le j \le N_\beta$.  
\\ 

Subsequently, we consider the case (a2). We have 
\begin{align*}
\bigl[ \mathscr{I}_{n} \bigl( \theta \bigr) \bigr]_{ij} 
 & = \tfrac{1}{n} \sum_{m=1}^n \widetilde{H}_{ij, 1} (\sample{m-1}, \theta)
 + \tfrac{1}{n} \sum_{m=1}^n \widetilde{H}_{ij, 2} (\sample{m-1}, \theta) 
+ \tfrac{1}{n} \sum_{m = 1}^n \sum_{1 \le k \le N_{S_1}} 
 R_k  (1, \sample{m-1}, \theta) \eta_{S_1, m-1}^k (\Delta_n, \theta)  \\[0.1cm] 
 & \quad  + \tfrac{1}{n} \sum_{m = 1}^n \sum_{N_{S_1} + 1 \le k \le N_{S}} 
{R}_k  (1, \sample{m-1}, \theta) \bigl( 
V_{S, 0}^{k} ( \sample{m-1}, \truebeta_{S} ) 
-  V_{S, 0}^{k} ( \sample{m-1}, \beta_{S} ) 
\bigr) +  \tfrac{1}{n} \sum_{m = 1}^n R (\Delta_n, \sample{m-1}, \theta),    
\end{align*}
for some $R_k, R \in \mathcal{S}, \, 1 \le k \le N_S$, where we have defined $\widetilde{H}_{ij, k} (x, \theta)$, for $(x, \theta) \in \mathbb{R}^N \times \Theta$, $k = 1,2$, as: 
\begin{align*}
 \widetilde{H}_{ij, 1} (x, \theta) 
& \equiv  \tfrac{1}{6} \bigl( \partial^{\beta}_i \mathcal{L} V_{S_1, 0} (x, \theta) \bigr)^\top 
\Lambda_{S_1 S_1} (x, \theta) \, 
\partial^{\beta}_j \mathcal{L}^2 V_{S_1, 0} (x, \theta) 
+ \tfrac{1}{2} \bigl( \partial^{\beta}_i \mathcal{L} V_{S_1, 0} (x, \theta) \bigr)^\top 
\Lambda_{S_1 S_2} (x, \theta) \,  
\partial^{\beta}_{j} \mathcal{L} V_{S_2, 0} (x, \theta)  \\[0.2cm] 
& \quad  + 
\bigl( \partial^{\beta}_i \mathcal{L} V_{S_1, 0} (x, \theta) \bigr)^\top
\Lambda_{S_1 R} (x, \theta) \, 
\partial^{\beta}_{j} V_{R, 0} (x, \beta_R);  \\[0.3cm] 
\widetilde{H}_{ij, 2} (x, \theta) 
& \equiv  \bigl( \tfrac{1}{3} \partial^{\beta}_i V_{S_2, 0} (x, \beta_{S_2}) \bigr)^\top \Lambda_{S_2 S_1} (x, \theta) \, \partial^{\beta}_j \mathcal{L}^2 V_{S_1, 0} (x, \theta)
+ \bigl( \partial^{\beta}_i V_{S_2, 0}  (x, \beta_{S_2})
\bigr)^\top \Lambda_{S_2 S_2} (x, \theta) 
\, \partial^{\beta}_{j} \mathcal{L} V_{S_2, 0} (x, \theta) \\[0.2cm] 
& \quad  + 
2 \bigl( \partial^{\beta}_i V_{S_2, 0} (x, \beta_{S_2}) \bigr)^\top \Lambda_{S_2 R} (x, \theta) 
\, \partial^{\beta}_{j} V_{R, 0} ( x, \beta_R ).  
\end{align*}
We then have $\widetilde{H}_{ij, k} (x, \theta) = 0$ for any $(x, \theta) \in \mathbb{R}^N \times \Theta$, $k = 1, 2$, from the same argument as in case (a1). Thus, making use of Lemma \ref{lemma:ergodic_thm}, convergence (\ref{eq:step1}) with condition (\ref{assump:bd_deriv}) and the consistency of the estimator, we obtain (\ref{eq:slln_conv1}) in case (a2).
\\ 

Finally, we consider the case (a3). We have 
\begin{align*}
\bigl[ \mathscr{I}_{n} \bigl( \theta \bigr) \bigr]_{ij} 
& =  
\tfrac{1}{n} \sum_{k = 1}^n \bar{H}_{ij, 1}  (\sample{k-1}, \theta)
+ \tfrac{1}{n} \sum_{k = 1}^n \bar{H}_{ij, 2}  (\sample{k-1}, \theta)
+ \tfrac{1}{n} \sum_{k = 1}^n \bar{H}_{ij, 3}  (\sample{k-1}, \theta) + \tfrac{1}{n \sqrt{\Delta_n}} \sum_{k  = 1}^n R (1, \sample{k - 1}, \theta) \, m_{k} (\Delta_n, \trueparam) 
\\[0.1cm]
& \quad + \tfrac{1}{n \sqrt{\Delta_n}} \sum_{k  = 1}^n 
\bigl( m_{k} (\Delta_n, \theta) - m_{k} (\Delta_n, \trueparam) \bigr)^\top 
\Lambda (\sample{k-1}, \theta) 
\, \partial^{\beta}_{(i,j)} v (\sample{k - 1}, \theta)  + \tfrac{1}{n} \sum_{k  = 1}^n \widetilde{R} (\sqrt{\Delta_n}, \sample{k-1}, \theta),
\end{align*}
for $R, \, \widetilde{R} \in \mathcal{S}$, where we have set, for $(x, \theta) \in \mathbb{R}^N \times \Theta$,
\begin{align}
& \begin{aligned} \label{eq:bar_H1} 
 \bar{H}_{ij, 1} (x, \theta)
& \equiv  \tfrac{1}{18} 
\bigl(  \partial^{\beta}_i \mathcal{L}^2 V_{S_1, 0} (x, \theta)  \bigr)^\top 
\Lambda_{S_1 S_1} (x, \theta) 
\, 
\partial^{\beta}_j \mathcal{L}^2 V_{S_1, 0} (x, \theta) \\[0.2cm] 
& \quad + \tfrac{1}{6} 
\bigl( \partial^{\beta}_i \mathcal{L}^2 V_{S_1, 0} (x, \theta) \bigr)^\top 
\Lambda_{S_1 S_2} (x, \theta) 
\, \partial^{\beta}_{j} \mathcal{L} V_{S_2, 0} (x, \theta)  \\[0.2cm] 
& \quad + 
\tfrac{1}{3} \bigl( \partial^{\beta}_i \mathcal{L}^2 V_{S_1, 0} (x, \theta) \bigr)^\top 
\Lambda_{S_1 R} (x, \theta) 
\, \partial^{\beta}_{j} V_{R, 0} (x, \beta_R);  \\[0.3cm] 
\end{aligned} \\
& 
\begin{aligned} \label{eq:bar_H2}
\bar{H}_{ij, 2}  (x, \theta) 
& \equiv \tfrac{1}{6} 
\bigl( \partial^{\beta}_i 
\mathcal{L} V_{S_2, 0} (x, \theta) \bigr)^\top 
\Lambda_{S_2 S_1} (x, \theta) 
\, \partial^{\beta}_j \mathcal{L}^2 V_{S_1, 0} (x, \theta) \\[0.2cm]  
& \quad + \tfrac{1}{2} \bigl( \partial^{\beta}_i \mathcal{L} V_{S_2, 0} (x, \theta) \bigr)^\top 
\Lambda_{S_2 S_2} (x, \theta) 
\, \partial^{\beta}_{j} \mathcal{L} V_{S_2, 0} (x, \theta) \\[0.2cm] 
& \quad  +  \bigl( \partial^{\beta}_i \mathcal{L} V_{S_2, 0}  (x, \theta) \bigr)^\top 
\Lambda_{S_2 R} (x, \theta)
\, \partial^{\beta}_{j} V_{R, 0} (x, \beta_R); \\[0.3cm]
\end{aligned} \\
& 
\begin{aligned} \label{eq:bar_H3} 
\bar{H}_{ij, 3} (x, \theta) 
& \equiv  \tfrac{1}{3} \bigl( \partial^{\beta}_i V_{R, 0} (x, \beta_{R}) \bigr)^\top
\Lambda_{R S_1} (x, \theta) 
\, \partial^{\beta}_j \mathcal{L}^2 V_{S_1, 0} (x, \theta) \\[0.2cm] 
& \quad + \bigl( \partial^{\beta}_i V_{R, 0} (x, \beta_{R})  \bigr)^\top 
\Lambda_{R S_2} (x, \theta) 
\, \partial^{\beta}_{j} \mathcal{L} V_{S_2, 0} (x, \theta) 
\\[0.2cm] 
& \quad  + 
2 \bigl( \partial^{\beta}_i V_{R, 0}  (x, \beta_{R}) \bigr)^\top 
\Lambda_{RR} (x, \theta) 
\, \partial^{\beta}_{j} V_{R, 0} (x, \beta_R),  
\end{aligned} 
\end{align}
and
\begin{align*}
v (x, \theta) \equiv \Bigl[ \tfrac{1}{6} 
\mathcal{L}^2 V_{S_1, 0} (x, \theta)^\top, \, 
\tfrac{1}{2} \mathcal{L} V_{S_2, 0} (x, 
\theta)^\top,  \, 
V_{R, 0} (x, \beta_R)^\top  \Bigr]^\top.    
\end{align*}
Notice that for any $(x, \theta) \in \mathbb{R}^N \times \Theta$, 
\begin{gather*}
\bar{H}_{ij, 1} (x, \theta) = 0, 
\qquad \bar{H}_{ij, 2} (x, \theta) = 0, \qquad 
\bar{H}_{ij, 3} (x, \theta) 
= 2 \bigl( \partial^{\beta}_i V_{R,0} (x, \beta_R) \bigr)^\top a_R^{-1} (x, \sigma) 
\, \partial^{\beta}_j V_{R,0} (x, \beta_R).  
\end{gather*} 
Furthermore, it follows that
\begin{align*}
 & \tfrac{1}{n \sqrt{\Delta_n}} \sum_{k  = 1}^n 
 \bigl( m_{k} (\Delta_n, \theta) - m_{k} (\Delta_n, \trueparam) \bigr)^\top
 \Lambda (\sample{k-1},  \theta)
 \partial^{\beta}_{(i,j)}  v (\sample{k - 1}, \theta)   \\[0.1cm]
 & = \tfrac{1}{n} \sum_{k = 1}^n 
 \bigl( V_{R, 0} (\sample{k -1}, \truebeta_R) 
 - V_{R, 0} (\sample{k -1}, \beta_R) 
 \bigr)^\top 
a_R^{-1} (\sample{k-1}, \sigma) \partial^{\beta}_{(i,j)} 
V_{R, 0} (\sample{k - 1}, \beta_R),  
\end{align*}
where we made use of similar arguments in the proof of Lemma \ref{lemma:matrix_2} in Appendix \ref{app:aux} for the term $\Lambda (\sample{k-1},  \theta)
 \partial^{\beta}_{(i,j)}  v (\sample{k - 1}, \theta)$. Hence, exploiting Lemmas \ref{lemma:ergodic_thm}, \ref{lemma:canonical_conv} and the consistency of estimator $\hat{\theta}_n$, we obtain that: 
\begin{align*}
\Bigl[\mathscr{I}_{n} \bigl(\trueparam + \lambda (\hat{\theta}_n - \trueparam) \bigr) \Bigr]_{ij} 
\probconv 2 \int \bigl( \partial^{\beta}_i V_{R, 0} (x, \truebeta_R) \bigr)^\top \, 
a_R^{-1} (x, \truesigma) \, 
\partial^{\beta}_j V_{R, 0} (x, \truebeta_R) \, \truedist (dx), 
\end{align*}
as \limit, for $N_{\beta} + 1 \le i, j \le N_{\theta}$. 
The proof of (\ref{eq:slln_conv1}) is now complete. 
\subsubsection{Proof of (\ref{eq:slln_conv2})}
We show (\ref{eq:slln_conv2}) when $1 \le i \le N_{\beta_{S_1}}$ and $N_{\beta} + 1 \le j \le N_\theta$. The convergence for the other cases can be deduced from a similar argument used in the proof of (\ref{eq:slln_conv1}) so we omit the proof. We have
\begin{align*}
\bigl[ \mathscr{I}_n \bigl( \theta \bigr) \bigr]_{ij}  
 & = \tfrac{\sqrt{\Delta_n^3}}{n} 
 \sum_{k = 1}^n  \sum_{1 \le k_1, k_2 \le N} 
 R_{k_1 k_2} (1, \sample{k -1}, \theta) 
 m_{k}^{k_1} (\Delta_n, \trueparam) 
 m_{k}^{k_2} (\Delta_n, \trueparam)  \nonumber \\[0.1cm] 
 & \quad + \tfrac{\sqrt{\Delta_n^3}}{n} 
 \sum_{k = 1}^n  \sum_{1 \le k_1, k_2 \le N} 
 \biggl\{ \widetilde{R}_{k_1 k_2} (1, \sample{k -1}, \theta) 
 \bigl( m_{k}^{k_1} (\Delta_n, \theta) - m_{k}^{k_1} (\Delta_n, \trueparam) \bigr) 
 \bigl( m_{k}^{k_2} (\Delta_n, \theta) - m_{k}^{k_2} (\Delta_n, \trueparam) \bigr) \biggr\} \nonumber \\ 
 & \quad + \tfrac{1}{n} 
 \sum_{k = 1}^n  \sum_{1 \le k_1 \le N} 
  {R}_{k_1} (1, \sample{k -1}, \theta) 
 \bigl( m_{k}^{k_1} (\Delta_n, \theta) - m_{k}^{k_1} (\Delta_n, \trueparam) \bigr) 
 + \tfrac{1}{n} 
 \sum_{k = 1}^n {R} (\sqrt{\Delta_n}, \sample{k -1}, \theta),  
\end{align*}
for some $R_{k_1 k_2}, \widetilde{R}_{k_1k_2}, R_{k_1}, R \in \mathcal{S}$.  
Thus, we immediately obtain (\ref{eq:slln_conv2}) for $1 \le i \le N_{\beta_{S_1}}, \, N_\beta + 1 \le j \le N_\theta$ from Lemma \ref{lemma:ergodic_thm}, \ref{lemma:canonical_conv} and (\ref{eq:step1})-(\ref{eq:step2}).  
\subsubsection{Proof of (\ref{eq:slln_conv3})} 
It holds that for $N_\beta + 1 \le i, j \le N_\theta$, $\theta  = (\beta_{S_1}, \beta_{S_2}, \beta_R, \sigma) \in \Theta$,
\begin{align*}
\bigl[ \mathscr{I}_{n} \bigl( \theta \bigr) \bigr]_{ij}   
& = \tfrac{1}{n} \sum_{k = 1}^n  
m_{k} (\Delta_n, \trueparam)^\top \, 
\partial^{\sigma}_{(i, j)} \Lambda (\sample{k-1}, \theta) \, 
m_{k} (\Delta_n, \trueparam) 
+ \tfrac{1}{n} \sum_{k = 1}^n 
\partial^{\sigma}_{(i,j)} 
\log \bigl| \Sigma (\sample{k - 1}, \theta ) \bigr| \\[0.1cm]
& \quad + \tfrac{1}{n} \sum_{k = 1}^n \sum_{1 \le k_1, k_2 \le N}  
R_{k_1 k_2} (1 , \sample{k - 1}, \theta) m_{k}^{k_1} (\Delta_n, \trueparam) \bigl( m_{k}^{k_2} (\Delta_n, \trueparam)  - m_{k}^{k_2} (\Delta_n, \theta) \bigr)  \\[0.2cm] 
& \quad + \tfrac{1}{n} \sum_{k = 1}^n \sum_{1 \le k_1, k_2 \le N}  
\biggl\{ \widetilde{R}_{k_1 k_2} (1 , \sample{k - 1}, \theta) 
 \bigl( m_{k}^{k_1} (\Delta_n, \trueparam)  - m_{k}^{k_1} (\Delta_n, \theta)\bigr)
\bigl( m_{k}^{k_2} (\Delta_n, \trueparam)  - m_{k}^{k_2} (\Delta_n, \theta)\bigr) \biggr\}  \\[0.2cm]  
& \quad + \tfrac{1}{n} \sum_{k = 1}^n \sum_{1 \le k_1 \le N}  
R_{k_1} (\sqrt{\Delta_n}, \sample{k - 1}, \theta) 
\bigl( m_{k}^{k_1} (\Delta_n, \trueparam) 
   - m_{k}^{k_1} (\Delta_n, \theta)\bigr)  + \tfrac{1}{n} \sum_{k = 1}^n R (\Delta_n, \sample{k-1}, \theta),
\end{align*}
for some $R_{k_1 k_2}, \widetilde{R}_{k_1 k_2}, R_{k_1}, R \in \mathcal{S}$. Making use of Lemmas \ref{lemma:ergodic_thm}--\ref{lemma:canonical_conv}, (\ref{eq:step1}), (\ref{eq:step2}) and the consistency of the estimator, we obtain as \limit, 
\begin{align*} 
& \bigl[ \mathscr{I}_{n} \bigl( \trueparam + \lambda (\hat{\theta}_{n} - \trueparam) \bigr) \bigr]_{ij}  
\probconv \int 
    \Bigl\{ \mathrm{tr} \bigl(
    \partial^{\sigma}_{(i,j)}
    \Lambda (x, \trueparam ) 
    \, \Sigma  (x, \trueparam ) \bigr) +  \partial^{\sigma}_{(i,j)} \log \bigl| \Sigma (x, \trueparam) \bigr|  \Bigr\} \truedist (dx) \\[0.2cm] 
& \quad = \int 
     \mathrm{tr}  \bigl(
    \partial^{\sigma}_i \Sigma  (x, \trueparam) \, 
    \Lambda  (x, \trueparam)  
    \partial^{\sigma}_j \Sigma  (x, \trueparam ) \,
    \Lambda  (x, \trueparam) \bigr) \truedist (dx),
\end{align*}
uniformly in $\lambda \in [0,1]$, where we applied the following two formulae to the above equation: 
\begin{gather*}
\partial^{\sigma}_{(i,j)}  
\log \bigl| \Sigma (x, \trueparam) \bigr| 
= - \mathrm{tr} \bigl(
    \partial^{\sigma}_{(i,j)}  
    \Lambda (x, \trueparam) 
    \, \Sigma  (x, \trueparam) \bigr) 
  - \mathrm{tr} \bigl( \partial^{\sigma}_i \Sigma (x, \trueparam ) \partial^{\sigma}_j 
\Lambda (x, \trueparam )  \bigr); \\[0.2cm]
\mathrm{tr} \bigl( \partial^{\sigma}_i \Sigma (x, \trueparam) \partial^{\sigma}_j 
 \Lambda (x, \trueparam )  \bigr) 
= - \mathrm{tr} \bigl( \partial^{\sigma}_i  
\Sigma (x, \trueparam ) \, 
\Lambda  (x, \trueparam )
\partial^{\sigma}_j 
\Sigma  (x, \trueparam ) \,  
\Lambda  (x, \trueparam )   \bigr). 
\end{gather*} 
The proof is now complete. 
\subsection{Proof of Lemma \ref{lemma:clt}} \label{appendix:pf_clt}
We write
$ \mathscr{C}_{n}^k (\theta) = \textstyle{\sum_{i = 1}^n \zeta_{i}^k (\theta), \;  \theta \in \Theta, \;  1 \le k \le N_\theta}$,
where we have set: 
\begin{align*}
\zeta_{i}^k (\theta) 
\equiv \bigl[ M_{n}^{-1} \bigr]_{kk} \times
\partial^{\theta}_k  
\bigl\{   
m_{i} (\Delta_n, \theta)^\top
\Lambda (\sample{i-1}, \theta) 
m_{i} (\Delta_n, \theta) 
+ \log |\Sigma (\sample{i-1}, \theta) |  \bigr\}. 
\end{align*}
To prove the assertion, it suffices to show, from Theorem 3.2 and 3.4 in \cite{hall:14}, that: 
\begin{enumerate}
    \item[(i)] If \limit \; with $\Delta_n = o (n^{-1/2})$, then 
    \begin{align}
        \sum_{i = 1}^n \mathbb{E}_{\trueparam}
        \bigl[ \zeta_{i}^k (\trueparam) | \mathcal{F}_{t_{i-1}} \bigr] \probconv 0, \quad 1 \le k \le N_{\theta}.  \label{eq:clt_conv1}
    \end{align}
    \item[(ii)] If \limit, then 
    \begin{align}
    \sum_{i = 1}^n \mathbb{E}_{\trueparam}
    \bigl[  
     \zeta_{i}^{k_1} (\trueparam) 
     \zeta_{i}^{k_2} (\trueparam) 
    | \mathcal{F}_{t_{i-1}} \bigr] \probconv 4 [\Gamma (\trueparam)]_{k_1 k_2},  
    \quad  1 \le k_1, k_2 \le N_{\theta}. \label{eq:clt_conv2}
    \end{align}
    \item[(iii)] If \limit, then 
    \begin{align}
    \sum_{i = 1}^n 
    \mathbb{E}_{\trueparam}
    \bigl[  
    \zeta_{i}^{k_1} (\trueparam) 
    \zeta_{i}^{k_2} (\trueparam) 
    \zeta_{i}^{k_3} (\trueparam) 
    \zeta_{i}^{k_4} (\trueparam) 
    | \mathcal{F}_{t_{i-1}} \bigr] \probconv 0, \quad 1 \le k_1, k_2, k_3, k_4 \le N_{\theta}. \label{eq:clt_conv3} 
    \end{align}
\end{enumerate}
In what follows, we will check convergences (\ref{eq:clt_conv1}) and (\ref{eq:clt_conv2}). One can prove (\ref{eq:clt_conv3}) following similar arguments and by noticing that the left-hand-side of (\ref{eq:clt_conv3}) involves $1/n^2$.   
\subsubsection{Proof of (\ref{eq:clt_conv1})}
We recall (\ref{eq:mean_expectation}) and (\ref{eq:squared_mean_expectation}) that are immediately obtained from the definition of $m_i (\Delta_n, \trueparam)$ in (\ref{eq:m_simple}), that is, for $1 \le k_1, k_2 \le N$, 
\begin{gather} 
\mathbb{E}_{\trueparam} \bigl[ m_i^{k_1} (\Delta_n, \trueparam) | \mathcal{F}_{t_{i-1}} \bigr] = R_1 (\sqrt{\Delta_n^3}, \sample{i-1}, \trueparam); \label{eq:mean_expecatation_simple} \\[0.2cm] 
\mathbb{E}_{\trueparam} \bigl[ m_i^{k_1} (\Delta_n, \trueparam)
m_i^{k_2} (\Delta_n, \trueparam) | \mathcal{F}_{t_{i-1}} \bigr] = [\Sigma (\sample{i-1}, \trueparam )]_{k_1 k_2}
+ R_2 ({\Delta_n}, \sample{i-1}, \trueparam) 
\label{eq:squared_expecatation_simple}  
\end{gather}
for $R_1, R_2 \in \mathcal{S}$. We then write $\zeta_{i}^k (\theta), \, \theta \in \Theta, \, 1 \le i \le n, \, 1 \le k \le N_{\theta}$ as $\zeta_i^k (\theta) = \zeta_{i, 1}^{k} (\theta) + \zeta_{i, 2}^{k} (\theta)$, where we have set:
\begin{gather*}
\zeta_{i, 1}^{k} (\theta)
= 2 \bigl[ M_{n}^{-1} \bigr]_{kk} 
\times 
\bigl\{ 
\bigl( \partial^{\theta}_k  
m_{i} (\Delta_n, \theta) \bigr)^\top 
\Lambda (\sample{i-1}, \theta) 
\, m_{i} (\Delta_n, \theta) \bigr\};  \nonumber \\[0.2cm]
\zeta_{i, 2}^{k} (\theta) 
= \bigl[ M_{n}^{-1} \bigr]_{kk}  
\times
\bigl\{ 
\partial^{\theta}_k\log | \Sigma (\sample{i-1}, \theta) |
+ m_{i} (\Delta_n, \theta)^\top 
 \partial^{\theta}_k \Lambda (\sample{i-1}, \theta) 
 \, m_{i} (\Delta_n, \theta)  
\bigr\}.  
\end{gather*}
Exploiting (\ref{eq:mean_expecatation_simple}), we obtain
\begin{align*}
 \mathbb{E}_{\trueparam} \bigl[ \zeta_{i, 1 }^{k} (\trueparam)  | \mathcal{F}_{t_{i-1}} \bigr]
 & = \tfrac{1}{\sqrt{n}}  \sum_{j =1}^N R_k^j (1, \sample{i-1}, \trueparam) \,  
 \mathbb{E}_{\trueparam} \bigl[ m_{i}^j (\Delta_n, \trueparam) | \mathcal{F}_{t_{i-1}} \bigr] = \tfrac{1}{n}  \widetilde{R}_k^j (\sqrt{n \Delta_n^{3}}, \sample{i-1}, \trueparam),  
\end{align*}
for $R_k^j, \, \widetilde{R}_k^j \in \mathcal{S}$. 
Thus, from Lemma \ref{lemma:ergodic_thm}, we obtain 
\begin{align*}
\sum_{i = 1}^n \mathbb{E}_{\trueparam} \bigl[ \zeta_{i, 1 }^{k} (\trueparam)  | \mathcal{F}_{t_{i-1}} \bigr] 
= \tfrac{1}{n} \sum_{i=1}^n \widetilde{R}_k^j (\sqrt{n \Delta_n^{3}}, \sample{i-1}, \trueparam) \probconv 0,  
\end{align*}
if \limit \; and $\Delta_n = o (n^{-1/2})$. Next, we consider the second term $\zeta_{i, 2}^{k} (\theta)$. First, notice that for $N_{\beta_{S}} + 1 \le k \le N_{\beta}$, $\zeta_{i, 2}^k ( \trueparam) =0$, since $\Sigma (\sample{i-1}, \theta)$ and $\Lambda (\sample{i-1}, \theta)$ are independent of $\beta_R \in \Theta_{\beta_R}$. For $1 \le k \le N_{\beta_S}$ and $N_{\beta} + 1 \le k \le N_{\theta}$, we apply (\ref{eq:squared_expecatation_simple}) to obtain   
\begin{align*}
\mathbb{E}_{\trueparam} [\zeta^k_{i, 2} (\trueparam) | \mathcal{F}_{t_{i-1}} ]
& = [M_n^{-1}]_{kk} \times \Bigl\{
 \partial^{\theta}_k 
 \log |\Sigma (\sample{i-1}, \trueparam) |
 + \mrm{tr} \bigl(\partial^{\theta}_k  
 \Lambda (\sample{i-1}, \trueparam) 
 \Sigma (\sample{i-1}, \trueparam) \bigr) 
 \Bigr\}  + \tfrac{1}{n} R_k (\sqrt{n \Delta_n^2}, \sample{i-1}, \trueparam) \\[0.2cm]
& = \tfrac{1}{n} R_k (\sqrt{n \Delta_n^2}, \sample{i-1}, \trueparam)
\end{align*}
for $R_k \in \mathcal{S}$, where we used: 
\begin{align} \label{eq:deriv_logdet}
\partial^{\theta}_k \log |\Sigma (x, \theta)| 
= - \mathrm{tr} \bigl( \partial^{\theta}_k \Lambda (x, \theta) \Sigma (x, \theta) \bigr), \quad (x, \theta) \in \mathbb{R}^N \times \Theta. 
\end{align}
Thus, we have from Lemma \ref{lemma:ergodic_thm} that if \limit \;  with $\Delta_n = o (n^{-1/2})$, then
\begin{align*}
    \sum_{i = 1}^n \mathbb{E}_{\trueparam} 
    \bigl[  
    \zeta_{i, 2}^{k} (\trueparam)  | \mathcal{F}_{t_{i-1}} 
    \bigr]
    = \tfrac{1}{n} \sum_{i = 1}^n R_k ( \sqrt{n \Delta_n^{2}}, \sample{i-1}, \trueparam) \probconv 0, 
\end{align*}
and the proof of (\ref{eq:clt_conv1}) is now complete. 
\subsubsection{Proof of (\ref{eq:clt_conv2})} 
For simplicity, we write  
\begin{align*}
    \mathscr{Y}_{k_1 k_2} (\trueparam) \equiv \sum_{i = 1}^n 
    \mathbb{E}_{\trueparam}
    \bigl[  
     \zeta_{i}^{k_1} (\trueparam) 
     \zeta_{i}^{k_2} (\trueparam) 
    | \mathcal{F}_{t_{i-1}} 
    \bigr].  
\end{align*} 
We have that for $1 \le k_1 \le N_{\beta_S}$, $N_{\beta_S} + 1 \le k_2 \le N_{\beta}$, $N_{\beta} + 1 \le k_3 \le N_{\theta}$, %
\begin{align}
& 
\begin{aligned}
\sqrt{n} \zeta_{i}^{k_1} (\trueparam) 
 & = - 2 \, 
 \mu_{k_1} (\sample{i-1}, \trueparam)^\top  \Lambda (\sample{i-1},  \trueparam) 
\; m_{i} (\Delta_n, \trueparam)  \\[0.2cm]
& \quad 
+  \sum_{1 \le j_1, j_2 \le N} 
R_{k_1}^{j_1 j_2} (\sqrt{\Delta_n}, \sample{i-1}, \trueparam)
m_{i}^{j_1} (\Delta_n, \trueparam)
m_{i}^{j_2} (\Delta_n, \trueparam) 
+  \sum_{1 \le j_1 \le N} 
R_{k_1}^{j_1} (\sqrt{\Delta_n}, \sample{i-1}, \trueparam)
m_{i}^{j_1} (\Delta_n, \trueparam);  
\end{aligned}  \label{eq:partial_beta_S}  \\[0.3cm]
& 
\begin{aligned}
\sqrt{n} \zeta_{i}^{k_2} (\trueparam) 
&  = - 2 \, 
\mu_{k_2} (\sample{i-1}, \trueparam)^\top
\Lambda (\sample{i-1}, \trueparam)
m_{i} (\Delta_n, \trueparam); 
%
\end{aligned} \label{eq:partial_beta_R} \\[0.3cm] 
& 
\begin{aligned}
\sqrt{n} \zeta_{i}^{k_3} (\trueparam) 
& = m_{i} (\Delta_n, \trueparam)^\top
\, 
\bigl( 
\partial^{\theta}_{k_3}  
\Lambda (\sample{i-1}, \trueparam) \bigr)
\, m_{i} (\Delta_n, \trueparam) + \partial^{\theta}_{k_3} 
\log | \Sigma (\sample{i-1} \trueparam) | \\[0.2cm]
& \quad 
 + \sum_{1 \le j_1, j_2 \le N} R_{k_3}^{j_1 j_2} ({\Delta_n}, \sample{i-1}, \trueparam) \, m_{i}^{j_1} (\Delta_n, \trueparam) m_{i}^{j_2} (\Delta_n, \trueparam) + \sum_{1 \le j_1, j_2 \le N} R_{k_3}^{j_1} (\sqrt{\Delta_n}, \sample{i-1}, \trueparam) \, m_{i}^{j_1} (\Delta_n, \trueparam), 
\end{aligned} \label{eq:partial_sigma}
\end{align}
for $R_{k_1}^{j_1 j_2}, \, R_{k_1}^{j_1},  
\, R_{k_3}^{j_1 j_2}, \, R_{k_3}^{j_1} \in \mathcal{S}$, where we defined $\mu_k : \mathbb{R}^N \times \Theta \to \mathbb{R}^N, \, 1 \le k \le N_{\beta}$ as: 
\begin{align*} 
& \mu_k (x, \theta) 
= 
\begin{cases}
\Bigl[
\bigl( \partial^\beta_k V_{S_1, 0} (x_S, \beta_{S_1}) \bigr)^\top, 
\, \mathbf{0}_{N_{S_2}}^\top,  \, \mathbf{0}_{N_{R}}^\top 
 \Bigr]^\top,   
 &  1 \le k \le N_{\beta_{S_1}}; \\[0.2cm] 
\Bigl[ \tfrac{1}{2} \bigl( \partial^\beta_k \mathcal{L} V_{S_1, 0} (x, \theta) \bigr)^\top, 
\,  \bigl( \partial^\beta_k  V_{S_2, 0} (x, \beta_{S_2}) \bigr)^\top,  \, \mathbf{0}_{N_{R}}^\top 
 \Bigr]^\top,   &  N_{\beta_{S_1}} + 1 \le k \le N_{\beta_{S}};  \\[0.2cm] 
 \Bigl[ \tfrac{1}{6} \bigl( \partial^\beta_k \mathcal{L}^2 V_{S_1, 0} (x, \theta) \bigr)^\top, 
\,  \tfrac{1}{2} \bigl( \partial^\beta_k \mathcal{L} V_{S_2, 0} (x, \theta) \bigr)^\top,  \, \bigl( \partial^\beta_k  V_{R, 0} (x, \beta_{R} ) \bigr)^\top 
 \Bigr]^\top,   &  N_{\beta_{S}} + 1 \le k \le N_{\beta},
\end{cases}
\end{align*}
for $(x, \theta) \in \mathbb{R}^N \times \Theta$. From Lemma \ref{lemma:ergodic_thm}, \ref{lemma:canonical_conv} and (\ref{eq:partial_beta_S}), we have that for $1 \le k_1, k_2 \le N_{\beta}$,  
\begin{align*}
& \mathscr{Y}_{k_1 k_2} (\trueparam)  
\probconv 
4 \sum_{1 \le j_1, j_2, j_3, j_4 \le N} 
\int \mu_{k_1}^{j_1} (x, \trueparam)
\bigl[ \Lambda (x, \trueparam) \bigr]_{j_1 j_2}  
\mu_{k_2}^{j_3} (x, \trueparam) 
\bigl[ \Lambda(x, \trueparam ) \bigr]_{j_3 j_4}  
\bigl[ \Sigma (x, \trueparam) \bigr]_{j_2 j_4} 
\, \truedist (dx) \\[0.2cm]
& \qquad\qquad \qquad  = 4 \sum_{1 \le j_1, j_2 \le N} 
\int \mu_{k_1}^{j_1} (x, \trueparam)  
\bigl[ \Lambda (x, \trueparam) \bigr]_{j_1 j_2} 
\mu_{k_2}^{j_2} (x, \trueparam) 
\truedist (dx)  \equiv 4 \int  \widetilde{\mathscr{Y}}_{k_1 k_2} (x, \trueparam) \,  \truedist (dx). 
\end{align*}
%
%
%
We then have, for $1 \le k_1, k_2 \le N_{\beta_{S_1}}$, 
\begin{align}
\widetilde{\mathscr{Y}}_{k_1 k_2} (x, \trueparam) 
& = \bigl( \partial^\beta_{k_1} V_{S_1, 0} (x_S, \truebeta_{S_1}) \bigr)^\top
 \Lambda_{S_1 S_1} (x, \trueparam)
\, \partial^\beta_{k_2} V_{S_1, 0} (x_S, \truebeta_{S_1}) 
\nonumber \\[0.1cm]
 & = 720  \bigl( \partial^\beta_{k_1} V_{S_1, 0} (x_S, \truebeta_{S_1}) \bigr)^\top
 a_{S_1}^{-1} (x, \trueparam)  
\, \partial^\beta_{k_2} V_{S_1, 0} (x_S, \truebeta_{S_1}),
\label{eq:Y_1}
\end{align}
for $N_{\beta_{S_1}} + 1 \le k_1, k_2 \le N_{\beta_S}$, 
\begin{align}
\widetilde{\mathscr{Y}}_{k_1 k_2} (x, \trueparam) 
& = \tfrac{1}{4} \bigl( \partial^\beta_{k_1} \mathcal{L} V_{S_1,0 } (x, \trueparam) \bigr)^\top \Lambda_{S_1 S_1} (x, \trueparam) 
\, \partial^\beta_{k_2} \mathcal{L} V_{S_1,0 } (x, \trueparam) 
\nonumber \\[0.1cm] 
& \quad +  \tfrac{1}{2} \bigl( \partial^\beta_{k_1} \mathcal{L} V_{S_1,0 } (x, \trueparam) \bigr)^\top \Lambda_{S_1 S_2} (x, \trueparam) 
\, \partial^\beta_{k_2} V_{S_2,0 } (x, \truebeta_{S_2} ) 
\nonumber \\[0.1cm] 
& \quad +  \tfrac{1}{2} \bigl( 
\partial^\beta_{k_1} V_{S_2,0 } (x, \truebeta_{S_2} ) 
\bigr)^\top \Lambda_{S_2 S_1} (x, \trueparam) 
\, \partial^\beta_{k_2} \mathcal{L} V_{S_1,0 } (x, \trueparam)  
\nonumber \\[0.1cm]   
& \quad +  \bigl( 
\partial^\beta_{k_1} V_{S_2,0 } (x, \truebeta_{S_2} ) 
\bigr)^\top \Lambda_{S_2 S_2} (x, \trueparam) 
\, \partial^\beta_{k_2} V_{S_2,0 } (x, \truebeta_{S_2} ) 
\nonumber \\[0.2cm] 
& = 12 \bigl( 
\partial^\beta_{k_1} V_{S_2,0 } (x, \truebeta_{S_2} ) 
\bigr)^\top a_{S_2}^{-1} (x, \trueparam) 
\, \partial^\beta_{k_2} V_{S_2,0 } (x, \truebeta_{S_2} ),  
\label{eq:Y_2} 
\end{align}
and for $N_{\beta_S} + 1 \le k_1, k_2 \le N_{\beta}$, 
\begin{align*}
\widetilde{\mathscr{Y}}_{k_1 k_2} (x, \trueparam)  
= \tfrac{1}{2} \sum_{j = 1}^3 \bar{H}_{k_1 k_2, j} (x, \trueparam)
= \bigl( \partial^\beta_{k_1} V_{R, 0} (x, \truebeta_R) \bigr)^\top
\, a_R^{-1} (x, \truesigma) 
\, \partial^\beta_{k_2} V_{R, 0} (x, \truebeta_R),  
\end{align*}
where $\bar{H}_{k_1 k_2, j} \, 1 \le j \le 3$ are defined as (\ref{eq:bar_H1}), (\ref{eq:bar_H2}) and (\ref{eq:bar_H3}). Note that we made use of Lemma \ref{lemma:matrix}, \ref{lemma:lambda} to derive (\ref{eq:Y_1})--(\ref{eq:Y_2}). 
For $1 \le k_1 \le N_{\beta}, \, N_{\beta} + 1 \le k_2 \le N_\theta$, it follows from (\ref{eq:partial_beta_S}), (\ref{eq:partial_beta_R}) and (\ref{eq:partial_sigma}) that 
\begin{align*}
\mathscr{Y}_{k_1 k_2} (\trueparam) 
\probconv 0, 
\end{align*}
if \limit. For $N_{\beta_R} + 1 \le k_1, k_2 \le N_{\theta}$, it follows that if \limit, then 
\begin{align*} 
{\mathscr{Y}}_{k_1 k_2} (\trueparam) \probconv 
 \widetilde{\mathscr{Y}}_{k_1 k_2, 1} (\trueparam)
 + \widetilde{\mathscr{Y}}_{k_1 k_2, 2} (\trueparam) 
 + \widetilde{\mathscr{Y}}_{k_1 k_2, 3} (\trueparam),  
\end{align*}
where we have set:
\begin{align*}
\widetilde{\mathscr{Y}}_{k_1 k_2, 1} (\trueparam)   
&\equiv  \sum_{1 \le j_1, j_2, j_3, j_4 \le N} \int   
[ \partial^{\theta}_{k_1} \Lambda(x, \trueparam) ]_{j_1 j_2} 
[ \partial^{\theta}_{k_2} \Lambda (x, \trueparam) ]_{j_3 j_4} \mathscr{W}_{j_1 j_2 j_3 j_4} (x, \trueparam) \, \truedist (dx); \\[0.2cm]
\widetilde{\mathscr{Y}}_{k_1 k_2, 2}  (\trueparam) 
&\equiv \sum_{1 \le j_1, j_2 \le N} \int \Bigl\{ 
[\partial^{\theta}_{k_1} \Lambda(x, \trueparam)]_{j_1 j_2} 
[\Sigma (x, \trueparam)]_{j_1 j_2}
\partial^{\theta}_{k_2} \log |\Sigma (x, \trueparam)| \\[0.1cm]  
& \qquad \qquad \qquad \qquad \qquad \qquad  + 
[\partial^{\theta}_{k_2} \Lambda (x, \trueparam)]_{j_1 j_2} 
[\Sigma (x, \trueparam) ]_{j_1 j_2} 
\partial^{\theta}_{k_1} \log |\Sigma (x, \trueparam)| \Bigr\} \truedist (dx) \\[0.2cm]
&= - \sum_{1 \le j_1, j_2, j_3, j_4 \le N} 
\int \Bigl\{ 
[\partial^{\theta}_{k_1} \Lambda(x, \trueparam)]_{j_1 j_2} 
[ \Sigma (x, \trueparam) ]_{j_1 j_2}
[\partial^{\theta}_{k_2} \Lambda (x, \trueparam)]_{j_3 j_4}
[\Sigma (x, \trueparam)]_{j_3 j_4} \\[0.2cm] 
& \qquad \qquad \qquad 
+ [\partial^{\theta}_{k_2} \Lambda (x, \trueparam)]_{j_1 j_2} 
[ \Sigma (x, \trueparam) ]_{j_1 j_2} 
[ \partial^{\theta}_{k_1} \Lambda (x, \trueparam) ]_{j_3 j_4}
[\Sigma (x, \trueparam) ]_{j_3 j_4} \Bigr\} 
\truedist (dx);  
\\[0.3cm]  
 \widetilde{\mathscr{Y}}_{k_1 k_2, 3} (\trueparam) 
&\equiv  \int 
\partial^{\theta}_{k_1}  \log |\Sigma (x, \trueparam)|  
\partial^{\theta}_{k_2}  \log |\Sigma (x, \trueparam)| \truedist (dx) \\[0.2cm]
&= \sum_{1 \le j_1, j_2, j_3, j_4 \le N} \int 
[ \partial^{\theta}_{k_1} \Lambda (x, \trueparam) ]_{j_1 j_2} 
[ \Sigma (x, \trueparam) ]_{j_1 j_2} \, 
[\partial^{\theta}_{k_2} \Lambda (x, \trueparam)]_{j_3 j_4} \, 
[ \Sigma (x, \trueparam) ]_{j_3 j_4}  \truedist (dx), 
\end{align*}
with $\mathscr{W}_{j_1 j_2 j_3 j_4} : \mathbb{R}^N \times \Theta \to \mathbb{R}$ defined as follows, for 
$(x, \theta) \in \mathbb{R}^N \times \Theta$, 
\begin{align*} 
\mathscr{W}_{j_1 j_2 j_3 j_4} (x, \theta) 
& = [ \Sigma (x, \theta) ]_{j_1 j_2} 
[ \Sigma (x, \theta) ]_{j_3 j_4}  
+ 
[ \Sigma (x, \theta)]_{j_1 j_3} 
[ \Sigma (x, \theta)]_{j_2 j_4} 
+ 
[\Sigma (x, \theta)]_{j_1 j_4}
[\Sigma (x, \theta)]_{j_2 j_3}. 
\end{align*}
Notice that we have used (\ref{eq:deriv_logdet}) in the computation of $\widetilde{\mathscr{Y}}_{k_1 k_2, 2} (\trueparam)$ and  $\widetilde{\mathscr{Y}}_{k_1 k_2, 3} (\trueparam)$. 
Thus, we have 
\begin{align*} 
 \sum_{m = 1}^3 \widetilde{\mathscr{Y}}_{k_1 k_2, m} (\trueparam)  
&  = \sum_{1 \le j_1, j_2, j_3, j_4 \le N} \int 
\Bigl\{ [\partial^{\theta}_{k_1} \Lambda (x, \trueparam)]_{j_1 j_2}
[ \partial^{\theta}_{k_2}  \Lambda  (x, \trueparam)]_{j_3 j_4} \\[0.2cm] 
& \qquad \qquad \times \bigl( 
[\Sigma (x, \trueparam)]_{j_1 j_3} 
[ \Sigma (x, \trueparam) ]_{j_2 j_4} 
+ 
[\Sigma (x, \trueparam)]_{j_1 j_4} 
[ \Sigma (x, \trueparam)]_{j_2 j_3} 
\bigr) \Bigr\}   
\truedist (dx) \\[0.2cm]
& = 2 \int \mathrm{tr} \bigl(  
\partial^{\theta}_{k_1} \Sigma(x, \trueparam) \, 
\Lambda (x, \trueparam) \, 
\partial^{\theta}_{k_2} \Sigma(x, \trueparam)  \, 
\Lambda (x, \trueparam)  \bigr) \, \truedist (dx),  
\end{align*} 
where in the second equality, we have used the following formula:  
\begin{align*} 
 [\partial^{\theta}_k \Lambda (\trueparam)]_{j_1 j_2} 
 = - \sum_{1 \le j_3 j_4 \le N}  
 [\Lambda (x, \trueparam)]_{j_1 j_3}
 [\partial^{\theta}_k \Sigma (x, \trueparam)]_{j_3 j_4}
 [\Lambda (x, \trueparam)]_{j_4 j_2}, \ \  N_{\beta}  +  1 \le k \le N_{\theta}.    
\end{align*}
Furthermore, for other cases of $1 \le k_1, k_2 \le N_{\theta}$, it holds that: 
\begin{align*}
    \widetilde{\mathscr{Y}}_{k_1 k_2} (x, \trueparam) = 0, 
\end{align*}
where we have used Lemma \ref{lemma:matrix} and that for $N_{\beta_S} + 1 \le k \le N_{\beta}$, 
\begin{align*} 
\Lambda (x, \theta) \, \mu_k (x, \theta) 
= \Lambda (x, \theta)  \begin{bmatrix}
\Sigma_{S_1 R} (x, \theta) \\[0.1cm]
\Sigma_{S_2 R} (x, \theta) \\[0.1cm]  
\Sigma_{R R} (x, \theta)  
\end{bmatrix}
a_R^{-1} (x, \sigma) \partial^{\beta}_k V_{R, 0} (x, \beta_R)  
=  
\begin{bmatrix} 
\mathbf{0}_{N_{S_1}} \\[0.1cm]
\mathbf{0}_{N_{S_2}} \\[0.2cm] 
a_R^{-1} (x, \sigma) \partial^{\beta}_{k} V_{R, 0} (x, \beta_R) 
\end{bmatrix}. 
\end{align*}
for $x \in \mathbb{R}^N, \, \theta = (\beta_{S}, \beta_R, \sigma) \in \Theta$. The proof is now complete. 

\section{Proof for Case Study in Section \ref{sec:case_2}} \label{appendix:pf_case_study}
From (\ref{eq:toy2_cont}),  we have $(\hat{\sigma}_n)^2 = F_{1, n} + F_{2, n} + F_{3, n}$, where  
\begin{align*} 
    F_{1, n} 
    & \equiv \tfrac{6}{ \Delta_n^3} \times 
    \tfrac{1}{n} \sum_{i = 0}^{n-1} \bigl( \hat{p}_{i+1} - \hat{p}_i -  s_i \Delta_n  + s_i \tfrac{\Delta_n^2}{2}  \bigr)^2;  \nonumber \\ 
    F_{2, n} 
    & \equiv - \tfrac{6}{\Delta_n^2} \times 
    \tfrac{1}{n} \sum_{i = 0}^{n-1} 
    \bigl( \hat{p}_{i+1} - \hat{p}_i - s_i \Delta_n + s_i \tfrac{\Delta_n^2}{2} \bigr) \bigl( s_{i+1} - s_i - s_i \Delta_n  \bigr);   \nonumber \\ 
    F_{3, n}
    & \equiv \tfrac{2}{\Delta_n} \times \tfrac{1}{n} 
    \sum_{i = 0}^{n-1} 
    \bigl( s_{i+1} - s_i - s_i \Delta_n  \bigr)^2  
\end{align*}  
with the hidden components $\hat{p}_i$ estimated by numerical differentiation: 
\begin{align*}
   \hat{p}_i  = \frac{q_{i+1} - q_i}{\Delta_n},\quad  0 \le i \le n. 
\end{align*}
Since the rough component $s_t$ follows the linear SDE, the solution is explicitly given as: for $u \in [t_i, t_{i+1})$ 
\begin{align*}
  s_{u} = s_{t_i}  e^{ - (u - t_i)} 
  + \truesigma \int_{t_i}^{u} e^{- (u  - v)} dB_v  
\end{align*}
under the true parameter $\truesigma$. Thus, we have 
\begin{align*}
\hat{p}_{i}  
 & = p_i + \tfrac{s_{i}}{\Delta_n} \int_{t_i}^{t_{i+1}} \int_{t_i}^{u} e^{- (v-t_i)} dv du    
 + \tfrac{\truesigma}{\Delta_n} \int_{t_i}^{t_{i+1}} \int_{t_i}^{u} \int_{t_i}^{v} e^{- (v-w)} d B_w dv du, 
\end{align*}
and then  
\begin{align*}
 & \hat{p}_{i+1} - \hat{p}_{i}  = p_{i+1} - p_i +  \tfrac{1}{\Delta_n} 
 \bigl( \Delta_n - ( 1  - e^{- \Delta_n}) \bigr) (s_{i+1} -s_i)  \nonumber  \\[0.2cm] 
 & \quad + \tfrac{\truesigma}{\Delta_n} \int_{t_{i+1}}^{t_{i+2}} \int_{t_{i+1}}^{u} \int_{t_{i+1}}^{v} e^{- (v-w)} d B_w dv du 
 - \tfrac{\truesigma}{\Delta_n} \int_{t_i}^{t_{i+1}} \int_{t_i}^{u} \int_{t_i}^{v} e^{- (v-w)} d B_w dv du, 
\end{align*}
where we used: 
\begin{align*}
\int_{t_i}^{t_{i+1}} \int_{t_i}^{u} e^{- (v- t_i)} dv du 
 = \Delta_n - (1 - e^{- \Delta_n}). 
\end{align*} 
Since $\Delta_n \in [0, 1)$ is assumed to be small, we use the Taylor expansion for the terms $e^{-\Delta_n}, \, e^{-(v-w)}$ and the stochastic Taylor expansion of $(p_{i+1}, s_{i+1})$ around $(p_{i}, s_{i})$ under the true parameter $\truesigma$ to obtain  
\begin{align}
\begin{aligned}
& \hat{p}_{i+1} - \hat{p}_{i} 
= s_{i} \Delta_n + \truesigma B_{t_{i+1}-t_i} \tfrac{\Delta_n}{2} 
+ \truesigma \int_{t_i}^{t_{i+1}} \int_{t_i}^u d B_v du  \\[0.2cm]
& \qquad  + \tfrac{\truesigma}{\Delta_n} 
\int_{t_{i+1}}^{t_{i+2}} \int_{t_{i+1}}^u \int_{t_{i+1}}^v dB_w dv du 
 -  \tfrac{\truesigma}{\Delta_n} 
\int_{t_{i}}^{t_{i+1}} \int_{t_{i}}^u \int_{t_i}^v d B_w dv du + \Delta_n^2 \xi_i,
\end{aligned} \label{eq:diff_p}
\end{align}
where the random variable $\{ \xi_i \}_i $ 
appearing in the remainder term satisfies $ \mathbb{E} | \xi_i |^2 \leq C $ for some constant $C > 0$. 
We express the Gaussian random variables as: 
\begin{gather*}
B_{t_{i+1} - t_i} = \sqrt{\Delta_n} \times  z_i^{(1)}, \quad 
\int_{t_i}^{t_{i+1}} \int_{t_i}^u d B_v du = \sqrt{\Delta_n^3} \times \Bigl( \tfrac{z_{i}^{(1)}}{2} + \tfrac{z_{i}^{(2)}}{2 \sqrt{3}} \Bigr), \quad \\ \int_{t_{i}}^{t_{i+1}} \int_{t_{i}}^u \int_{t_i}^v d B_w dv du  = \sqrt{\Delta_n^5} \times 
\Bigl( \tfrac{z_{i}^{(1)}}{6}  
 + \tfrac{z_{i}^{(2)}}{4 \sqrt{3}} 
 + \tfrac{z_{i}^{(3)}}{12\sqrt{5}} \Bigr),  
\end{gather*}
where $\{ z_{i}^{(j)} \}_{0 \le i  \le n+1, \, j = 1,2,3}$ is an i.i.d. sequence of standard normal random variables so that it holds
\begin{align*}
& \mathbb{E} \bigl[ (B_{t_{i+1} - t_i})^2  \bigr] = \Delta_n, \quad
\mathbb{E} \Bigl[ B_{t_{i+1} - t_i} \times \Bigl(\int_{t_i}^{t_{i+1}} \int_{t_i}^u d B_v du \Bigr)  \Bigr] = \tfrac{\Delta_n^2}{2}; \\[0.2cm]  
& \mathbb{E} \Bigl[ B_{t_{i+1} - t_i} \times \Bigl( \int_{t_{i}}^{t_{i+1}} \int_{t_{i}}^u \int_{t_i}^v d B_w dv du \Bigr)  \Bigr] = \tfrac{\Delta_n^3}{6}, 
\quad 
\mathbb{E} \Bigl[ \Bigl(\int_{t_i}^{t_{i+1}} \int_{t_i}^u d B_v du \Bigr)^2 \Bigr] = \tfrac{\Delta_n^3}{3}; \\[0.2cm]
\quad 
& \mathbb{E} \Bigl[ \Bigl(\int_{t_i}^{t_{i+1}} \int_{t_i}^u d B_v du \Bigr) \times 
\Bigl( \int_{t_{i}}^{t_{i+1}} \int_{t_{i}}^u \int_{t_i}^v d B_w dv du \Bigr) \Bigr] = \tfrac{\Delta_n^4}{8},
\quad 
\mathbb{E} \Bigl[ 
\Bigl( \int_{t_{i}}^{t_{i+1}} \int_{t_{i}}^u \int_{t_i}^v d B_w dv du \Bigr)^2 \Bigr] = \tfrac{\Delta_n^5}{20}.  
\end{align*}
Then,  (\ref{eq:diff_p}) is written as: 
\begin{align}
\hat{p}_{i+1} - \hat{p}_{i}  
& =  s_{i} \Delta_n 
+  \truesigma \sqrt{\Delta_n^3}  
\tfrac{z_{i}^{(1)}}{2} 
+  \truesigma \sqrt{\Delta_n^3}  
\Bigl( \tfrac{z_{i}^{(1)}}{2} + \tfrac{z_{i}^{(2)}}{2 \sqrt{3}}  \Bigr)  \nonumber \\
& \qquad  + \truesigma \sqrt{\Delta_n^3}  
 \Bigl( \tfrac{z_{i+1}^{(1)}}{6}  
 + \tfrac{z_{i+1}^{(2)}}{4 \sqrt{3}} 
 + \tfrac{z_{i+1}^{(3)}}{12\sqrt{5}} \Bigr)   
 - \truesigma \sqrt{\Delta_n^3} 
 \Bigl( \tfrac{z_{i}^{(1)}}{6}  
 + \tfrac{z_{i}^{(2)}}{4 \sqrt{3}} 
 + \tfrac{z_{i}^{(3)}}{12\sqrt{5}} \Bigr) 
 + \Delta_n^2 \xi_i  \nonumber  \\[0.2cm] 
& = s_{i} \Delta_n + \truesigma \sqrt{\Delta_n^3} 
\Bigl\{ \Bigl( \tfrac{z_{i+1}^{(1)}}{6}  
 + \tfrac{z_{i+1}^{(2)}}{4 \sqrt{3}} 
 + \tfrac{z_{i+1}^{(3)}}{12\sqrt{5}} \Bigr) 
 + \Bigl( \tfrac{5 z_{i}^{(1)}}{6}  
 + \tfrac{z_{i}^{(2)}}{4 \sqrt{3}} 
 - \tfrac{z_{i}^{(3)}}{12\sqrt{5}} \Bigr)  
\Bigr\} + \Delta_n^2 \xi_i.    \label{eq:diff_p_2}
\end{align} 
From the ergodicity of the process $\{s_t\}_t$ and (\ref{eq:diff_p_2}), we have that, 
as $n \to \infty$, $\Delta_n \to 0$ and $n \Delta_n \to \infty$, 
\begin{align} 
F_{1,n}  & = 
 \tfrac{6 (\truesigma)^2 }{n} \sum_{i = 0}^{n-1} \Bigl( \tfrac{z_{i+1}^{(1)}}{6} 
 + \tfrac{z_{i+1}^{(2)}}{4 \sqrt{3}} 
 + \tfrac{z_{i+1}^{(3)}}{12\sqrt{5}}  
 + \tfrac{5}{6}z_{i}^{(1)}
 + \tfrac{z_{i}^{(2)}}{4 \sqrt{3}} 
 - \tfrac{z_{i}^{(3)}}{12\sqrt{5}}    
\Bigr)^2  + \tfrac{1}{n} \sum_{i = 0}^n R_i^{(1)} (\Delta_n)  
\probconv \tfrac{23}{5} (\truesigma)^2;  \label{eq:conv_F1} \\[0.2cm] 
F_{2, n} 
& = - \tfrac{6 (\truesigma)^2}{n} \sum_{i = 0}^{n-1}
 \Bigl( \tfrac{z_{i+1}^{(1)}}{6} 
 + \tfrac{z_{i+1}^{(2)}}{4 \sqrt{3}} 
 + \tfrac{z_{i+1}^{(3)}}{12\sqrt{5}}  
 + \tfrac{5}{6}z_{i}^{(1)}
 + \tfrac{z_{i}^{(2)}}{4 \sqrt{3}} 
 - \tfrac{z_{i}^{(3)}}{12\sqrt{5}}    
\Bigr) z_i^{(1)} 
+ \tfrac{1}{n} \sum_{i = 0}^n R_i^{(2)} (\Delta_n) 
\probconv - 5 (\truesigma)^2, \label{eq:conv_F2}
\end{align}
where each $\{ R_i^{(1)} (\Delta_n) \}_i$ and $\{R_i^{(2)} (\Delta_n) \}_i$ is sequence of random variables such that for $0 \le i \le n$, $j = 1,2$, 
\begin{align*} 
 \{ \mathbb{E} [ |R_i^{(j)} (\Delta_n) |^2 ]  \} \leq C \Delta_n
\end{align*} 
for some constant $C > 0$. Similarly, we have that 
\begin{align} \label{eq:conv_F3}
 F_{3, n}  \probconv  2 (\truesigma)^2.   
\end{align}
From (\ref{eq:conv_F1}), (\ref{eq:conv_F2}) and (\ref{eq:conv_F3}), we immediately obtain the convergence (\ref{eq:conv_toy_2}). 
\section{Kalman Filter for Sub-Class of (\ref{eq:hypo-II})} \label{appendix:kalman} 
For simplicity, we write $x_i = (x_{S_1,i}, x_{S_2, i}, x_{R,i}) \in \mathbb{R}^N = \mathbb{R}^{N_{S_1}} \times  \mathbb{R}^{N_{S_2}} \times   \mathbb{R}^{N_{R}}$ for the state of scheme (\ref{eq:scheme_linear}) at time $t_i$. Component $x_{S_1,i}$ is observable and $h_{i} = (x_{S_2, i}, x_{R, i}) \in \mathbb{R}^{N_H}$, $N_H = N_{S_1} + N_R$, is the hidden component,  in agreement with applications. 
Thus, scheme (\ref{eq:scheme_linear}) is now expressed as
\begin{align} \label{eq:scheme_linear_appendix}
x_{i+1} 
=  b (\Delta_n, x_{S_1, i}, \theta)  
+ A (\Delta_n, {x}_{S_1, i}, \theta) h_i
+ w (\Delta_n, \theta).  
\end{align}
We set $\Sigma (\Delta_n, \theta) = \mathbb{E}\,[\,w (\Delta_n, \theta) w (\Delta_n, \theta)^\top] $ and assume that $h_0 | x_{S_1, 0} \sim  \mathscr{N} ({m}_0, {Q}_0)$ for some 
${m}_0 \in \mathbb{R}^{N_H}$ and ${Q}_0 \in \mathbb{R}^{N_H \times N_H}$. Then, the filtering formula and the marginal likelihood are obtained as follows.
\begin{itemize}
\item \emph{Filtering Recursion}:  We have that
\begin{align} 
\label{eq:filter}
    h_{k} | x_{S_1, 0:k} \sim \mathscr{N} ({m}_k, {Q}_k), 
    \quad 0 \le k \le n, 
\end{align}
with the filter mean ${m}_k$ and covariance ${Q}_k$ given as: 
\begin{align*} 
{m}_k & =  \mu_{H, k-1} + \Lambda_{H S_1, k-1} \,  \Lambda^{-1}_{S_1 S_1, k-1} \,  \bigl( x_{S_1, k} -  \mu_{S_1, k-1} \bigr); \\[0.2cm]
Q_k & = \Lambda_{HH, k-1} - \Lambda_{H S_1, k-1} \, 
\Lambda_{S_1 S_1, k-1}^{-1} \Lambda_{S_1 H, k-1},
\end{align*}
where $\mu_{H, k-1} \in \mathbb{R}^{N_{H}}$, $\mu_{S_1, k-1} \in \mathbb{R}^{N_{S_1}}$, $\Lambda_{S_1 S_1, k-1} \in \mathbb{R}^{N_{S_1} \times N_{S_1} }$, $\Lambda_{S_1 H, k-1} \in \mathbb{R}^{N_{S_1} \times N_{H}}$, $\Lambda_{H S_1, k-1} \in \mathbb{R}^{N_{H} \times N_{S_1} }$, $\Lambda_{H H, k-1} \in \mathbb{R}^{N_{H} \times N_{H}}$ are found via the following equations:
\begin{align*}
\mu_{k-1} & =  
\begin{bmatrix}
    \mu_{S_1, k-1} \\[0.1cm]
    \mu_{H, k-1} 
\end{bmatrix}
= 
b (\Delta_n, x_{S_1, k-1}, \theta) + A (\Delta_n, x_{S_1, k-1}, \theta) 
m_{k-1}; \\[0.2cm] 
\Lambda_{k-1} 
& = 
\begin{bmatrix}
    \Lambda_{S_1 S_1, k-1}  & \Lambda_{S_1 H, k-1}  \\ 
    \Lambda_{H S_1, k-1}  & \Lambda_{H H, k-1}
\end{bmatrix} 
= \Sigma (\Delta_n, \theta) + A (\Delta_n, x_{S_1, k-1}, \theta) \, Q_{k-1} \, A (\Delta_n, x_{S_1, k-1}, \theta)^\top. 
\end{align*}
\item \emph{Marginal likelihood}: For a given initial distribution $p_{\theta} (x_{S_1, 0})$, we have that
\begin{align} \label{eq:marginal_likelihood}
 p_{\theta} (x_{S_1, 0:n}) = p_{\theta} (x_{S_1, 0}) \times 
\prod_{k = 1}^n p_{\theta} (x_{S_1, k} | x_{S_1, 0:k-1}), 
\end{align}
where $p_{\theta} (x_{S_1, k} | x_{S_1, 0:k-1})$ is the density of $x_{S_1, k}$ given $x_{S_1, 0:k-1}$ whose conditional distribution is given by:  
\begin{align*}
x_{S_1, k} \, | \, x_{S_1, 0:k-1} \sim \mathscr{N} (\mu_{S_1, k-1}, \, \Lambda_{S_1 S_1, k-1}). 
\end{align*}
\end{itemize}
\subsection{Derivation of Filter (\ref{eq:filter})} 
We assume that the filter in the previous time step is obtained as: 
\begin{align} \label{eq:previous_filter}
h_{k-1} | x_{S_1, 0:k-1} \sim \mathscr{N} (m_{k-1}, Q_{k-1}). 
\end{align}
It follows that 
\begin{align*}  
p_\theta (h_k | x_{S_1, 0:k} )
= \frac{p_\theta (x_k | x_{S_1, 0:k-1} ) }{p_\theta (x_{S_1, k} | x_{S_1, 0:k-1})}, 
\end{align*}
and 
\begin{align} \label{eq:conv_filter}
p_\theta (x_k | x_{S_1, 0:k-1} )  
& = \int p_\theta (x_k, h_{k-1} | x_{S_1, 0:k-1} )  d h_{k-1}
= \int p_\theta (x_k | x_{k-1} )  p_\theta (h_{k-1} | x_{S_1, 0:k-1} )  d h_{k-1}. 
\end{align}
From the definition of scheme (\ref{eq:scheme_linear_appendix}), we have that 
\begin{align} \label{eq:scheme_dist}
x_k | x_{k-1} \sim \mathscr{N} \bigl( b (\Delta_n, x_{S_1, k-1}, \theta)  
+ A (\Delta_n, {x}_{S_1, k-1}, \theta) h_{k-1}, \Sigma (\Delta_n, \theta) \bigr). 
\end{align}
From (\ref{eq:previous_filter}), (\ref{eq:conv_filter}) and (\ref{eq:scheme_dist}), we obtain 
\begin{align} \label{eq:state_given_obs}
x_k | x_{S_1, 0:k-1} 
\sim \mathscr{N} \bigl(\mu_{k-1}, \Lambda_{k-1}\bigr), 
\end{align} 
where
\begin{gather*}
\mu_{k-1} 
= b (\Delta_n, x_{S_1, k-1}, \theta)  
+ A (\Delta_n, {x}_{S_1, k-1}, \theta) m_{k-1}, 
\qquad  
\Lambda_{k-1} = \Sigma (\Delta_n, \theta) + 
A (\Delta_n, {x}_{S_1, k-1}, \theta) 
\, Q_{k-1}  \, 
A (\Delta_n, {x}_{S_1, k-1}, \theta)^\top. 
\end{gather*}
Finally, applying the conditional Gaussian distribution formula, we obtain (\ref{eq:filter}).
\subsection{Derivation of the Marginal Likelihood (\ref{eq:marginal_likelihood})}
From the marginal of the Gaussian distribution (\ref{eq:state_given_obs}) for $x_{k} | x_{S_1, 0:k-1}$, we immediately obtain: 
\begin{align*} 
x_{S_1, k} | x_{S_1, 0:k-1} \sim \mathscr{N} 
\bigl( \mu_{S_1, k-1}, \,  \Lambda_{S_1 S_1, k-1} \bigr), 
\end{align*}
where $\mu_{S_1, k-1} = \mrm{proj}_{1, N_{S_1}} ( \mu_{k-1}) $, and $\Lambda_{S_1 S_1, k-1} = \bigl[\Lambda_{k-1}^{ij} \bigr]_{1 \le i,j \le N_{S_1}}$. 

\end{appendices}






\bibliographystyle{elsarticle-harv}  
\bibliography{Hypo}  

\begin{thebibliography}{37}
\expandafter\ifx\csname natexlab\endcsname\relax\def\natexlab#1{#1}\fi
\providecommand{\url}[1]{\texttt{#1}}
\providecommand{\href}[2]{#2}
\providecommand{\path}[1]{#1}
\providecommand{\DOIprefix}{doi:}
\providecommand{\ArXivprefix}{arXiv:}
\providecommand{\URLprefix}{URL: }
\providecommand{\Pubmedprefix}{pmid:}
\providecommand{\doi}[1]{\href{http://dx.doi.org/#1}{\path{#1}}}
\providecommand{\Pubmed}[1]{\href{pmid:#1}{\path{#1}}}
\providecommand{\bibinfo}[2]{#2}
\ifx\xfnm\relax \def\xfnm[#1]{\unskip,\space#1}\fi
\bibitem[{Ayaz et~al.(2021)Ayaz, Tepper, Br{\"u}nig, Kappler, Daldrop and
  Netz}]{ay:21}
\bibinfo{author}{Ayaz, C.}, \bibinfo{author}{Tepper, L.},
  \bibinfo{author}{Br{\"u}nig, F.N.}, \bibinfo{author}{Kappler, J.},
  \bibinfo{author}{Daldrop, J.O.}, \bibinfo{author}{Netz, R.R.},
  \bibinfo{year}{2021}.
\newblock \bibinfo{title}{Non-{M}arkovian modeling of protein folding}.
\newblock \bibinfo{journal}{Proc. Natl. Acad. Sci. U. S. A.}
  \bibinfo{volume}{118}, \bibinfo{pages}{e2023856118}.
\bibitem[{Bally and Talay(1996)}]{bally:96}
\bibinfo{author}{Bally, V.}, \bibinfo{author}{Talay, D.}, \bibinfo{year}{1996}.
\newblock \bibinfo{title}{The law of the {E}uler scheme for stochastic
  differential equations: I. {C}onvergence rate of the distribution function}.
\newblock \bibinfo{journal}{Probab. Theory Relat. Fields}
  \bibinfo{volume}{104}, \bibinfo{pages}{43--60}.
\bibitem[{Buckwar et~al.(2020)Buckwar, Tamborrino and Tubikanec}]{buc:20}
\bibinfo{author}{Buckwar, E.}, \bibinfo{author}{Tamborrino, M.},
  \bibinfo{author}{Tubikanec, I.}, \bibinfo{year}{2020}.
\newblock \bibinfo{title}{Spectral density-based and measure-preserving {ABC}
  for partially observed diffusion processes. an illustration on {H}amiltonian
  {SDE}s}.
\newblock \bibinfo{journal}{Stat. Comput.} \bibinfo{volume}{30},
  \bibinfo{pages}{627--648}.
\bibitem[{Cass(2009)}]{cass:09}
\bibinfo{author}{Cass, T.}, \bibinfo{year}{2009}.
\newblock \bibinfo{title}{Smooth densities for solutions to stochastic
  differential equations with jumps}.
\newblock \bibinfo{journal}{Stoch. Process. their Appl.} \bibinfo{volume}{119},
  \bibinfo{pages}{1416--1435}.
\bibitem[{Ceriotti et~al.(2010)Ceriotti, Bussi and Parrinello}]{ce:10}
\bibinfo{author}{Ceriotti, M.}, \bibinfo{author}{Bussi, G.},
  \bibinfo{author}{Parrinello, M.}, \bibinfo{year}{2010}.
\newblock \bibinfo{title}{Colored-noise thermostats {\`a} la carte}.
\newblock \bibinfo{journal}{J. Chem. Theory Comput.} \bibinfo{volume}{6},
  \bibinfo{pages}{1170--1180}.
\bibitem[{Chen(2023)}]{ch:23}
\bibinfo{author}{Chen, N.}, \bibinfo{year}{2023}.
\newblock \bibinfo{title}{Stochastic methods for modeling and predicting
  complex dynamical systems: uncertainty quantification, state estimation, and
  reduced-order models}.
\newblock \bibinfo{publisher}{Springer Nature}.
\bibitem[{Coti~Zelati and Hairer(2021)}]{co:21}
\bibinfo{author}{Coti~Zelati, M.}, \bibinfo{author}{Hairer, M.},
  \bibinfo{year}{2021}.
\newblock \bibinfo{title}{A {N}oise-induced transition in the {L}orenz system}.
\newblock \bibinfo{journal}{Commun. Math. Phys.} \bibinfo{volume}{383},
  \bibinfo{pages}{2243--2274}.
\bibitem[{Ditlevsen and Samson(2019)}]{dit:19}
\bibinfo{author}{Ditlevsen, S.}, \bibinfo{author}{Samson, A.},
  \bibinfo{year}{2019}.
\newblock \bibinfo{title}{Hypoelliptic diffusions: {F}iltering and inference
  from complete and partial observations}.
\newblock \bibinfo{journal}{J. R. Stat. Soc., B: Stat. Methodol.}
  \bibinfo{volume}{81}, \bibinfo{pages}{361--384}.
\bibitem[{Ditlevsen et~al.(2023)Ditlevsen, Tamborrino and Tubikanec}]{dit:23}
\bibinfo{author}{Ditlevsen, S.}, \bibinfo{author}{Tamborrino, M.},
  \bibinfo{author}{Tubikanec, I.}, \bibinfo{year}{2023}.
\newblock \bibinfo{title}{Network inference in a stochastic multi-population
  neural mass model via approximate bayesian computation}.
\newblock \bibinfo{journal}{arXiv preprint arXiv:2306.15787} .
\bibitem[{Douc et~al.(2014)Douc, Moulines and Stoffer}]{douc:14}
\bibinfo{author}{Douc, R.}, \bibinfo{author}{Moulines, E.},
  \bibinfo{author}{Stoffer, D.}, \bibinfo{year}{2014}.
\newblock \bibinfo{title}{Nonlinear time series: {T}heory, methods and
  applications with {R} examples}.
\newblock \bibinfo{publisher}{CRC press}.
\bibitem[{Dureau et~al.(2013)Dureau, Kalogeropoulos and Baguelin}]{dur:13}
\bibinfo{author}{Dureau, J.}, \bibinfo{author}{Kalogeropoulos, K.},
  \bibinfo{author}{Baguelin, M.}, \bibinfo{year}{2013}.
\newblock \bibinfo{title}{Capturing the time-varying drivers of an epidemic
  using stochastic dynamical systems}.
\newblock \bibinfo{journal}{Biostatistics} \bibinfo{volume}{14},
  \bibinfo{pages}{541--555}.
\bibitem[{Ferretti et~al.(2020)Ferretti, Chardes, Mora, Walczak and
  Giardina}]{fe:20}
\bibinfo{author}{Ferretti, F.}, \bibinfo{author}{Chardes, V.},
  \bibinfo{author}{Mora, T.}, \bibinfo{author}{Walczak, A.M.},
  \bibinfo{author}{Giardina, I.}, \bibinfo{year}{2020}.
\newblock \bibinfo{title}{Building general {L}angevin models from discrete
  datasets}.
\newblock \bibinfo{journal}{Phys. Rev. X} \bibinfo{volume}{10},
  \bibinfo{pages}{031018}.
\bibitem[{Genon-Catalot and Jacod(1993)}]{genon:93}
\bibinfo{author}{Genon-Catalot, V.}, \bibinfo{author}{Jacod, J.},
  \bibinfo{year}{1993}.
\newblock \bibinfo{title}{On the estimation of the diffusion coefficient for
  multi-dimensional diffusion processes}, in: \bibinfo{booktitle}{Annales de
  l'IHP Probabilit{\'e}s et statistiques}, pp. \bibinfo{pages}{119--151}.
\bibitem[{Gloter and Yoshida(2020)}]{glot:20}
\bibinfo{author}{Gloter, A.}, \bibinfo{author}{Yoshida, N.},
  \bibinfo{year}{2020}.
\newblock \bibinfo{title}{Adaptive and non-adaptive estimation for degenerate
  diffusion processes}.
\newblock \bibinfo{journal}{arXiv preprint} ,
  \bibinfo{pages}{arXiv:2002.10164}.
\bibitem[{Gloter and Yoshida(2021)}]{glot:21}
\bibinfo{author}{Gloter, A.}, \bibinfo{author}{Yoshida, N.},
  \bibinfo{year}{2021}.
\newblock \bibinfo{title}{Adaptive estimation for degenerate diffusion
  processes}.
\newblock \bibinfo{journal}{Electron. J. Stat.} \bibinfo{volume}{15},
  \bibinfo{pages}{1424--1472}.
\bibitem[{Gobet and Labart(2008)}]{go:08}
\bibinfo{author}{Gobet, E.}, \bibinfo{author}{Labart, C.},
  \bibinfo{year}{2008}.
\newblock \bibinfo{title}{Sharp estimates for the convergence of the density of
  the {E}uler scheme in small time}.
\newblock \bibinfo{journal}{Electron. Commun. Probab.} \bibinfo{volume}{13},
  \bibinfo{pages}{352--363}.
\bibitem[{Hall and Heyde(1980)}]{hall:14}
\bibinfo{author}{Hall, P.}, \bibinfo{author}{Heyde, C.C.},
  \bibinfo{year}{1980}.
\newblock \bibinfo{title}{Martingale limit theory and its application}.
\newblock \bibinfo{publisher}{Academic press}.
\bibitem[{Iguchi et~al.(2024)Iguchi, Beskos and Graham}]{iguchi:22}
\bibinfo{author}{Iguchi, Y.}, \bibinfo{author}{Beskos, A.},
  \bibinfo{author}{Graham, M.M.}, \bibinfo{year}{2024}.
\newblock \bibinfo{title}{Parameter estimation with increased precision for
  elliptic and hypo-elliptic diffusions}.
\newblock \bibinfo{journal}{arxiv preprint} , \bibinfo{pages}{arXiv:2211.16384.
  To appear in Bernoulli}.
\bibitem[{Iguchi and Yamada(2021)}]{iguchi:21-2}
\bibinfo{author}{Iguchi, Y.}, \bibinfo{author}{Yamada, T.},
  \bibinfo{year}{2021}.
\newblock \bibinfo{title}{Operator splitting around {E}uler-{M}aruyama scheme
  and high order discretization of heat kernels}.
\newblock \bibinfo{journal}{ESAIM: Math. Model. Numer. Anal.}
  \bibinfo{volume}{55}, \bibinfo{pages}{S323--S367}.
\bibitem[{Kalliadasis et~al.(2015)Kalliadasis, Krumscheid and
  Pavliotis}]{kal:15}
\bibinfo{author}{Kalliadasis, S.}, \bibinfo{author}{Krumscheid, S.},
  \bibinfo{author}{Pavliotis, G.A.}, \bibinfo{year}{2015}.
\newblock \bibinfo{title}{A new framework for extracting coarse-grained models
  from time series with multiscale structure}.
\newblock \bibinfo{journal}{J. Comput. Phys.} \bibinfo{volume}{296},
  \bibinfo{pages}{314--328}.
\bibitem[{Kessler(1997)}]{kess:97}
\bibinfo{author}{Kessler, M.}, \bibinfo{year}{1997}.
\newblock \bibinfo{title}{Estimation of an ergodic diffusion from discrete
  observations}.
\newblock \bibinfo{journal}{Scand. J. Stat.} \bibinfo{volume}{24},
  \bibinfo{pages}{211--229}.
\bibitem[{Kloeden and Platen(1992)}]{kloe:92}
\bibinfo{author}{Kloeden, P.E.}, \bibinfo{author}{Platen, E.},
  \bibinfo{year}{1992}.
\newblock \bibinfo{title}{Numerical Solution of Stochastic Differential
  Equations}.
\newblock \bibinfo{publisher}{Springer}.
\bibitem[{Leimkuhler and Matthews(2015)}]{lei:15}
\bibinfo{author}{Leimkuhler, B.}, \bibinfo{author}{Matthews, C.},
  \bibinfo{year}{2015}.
\newblock \bibinfo{title}{{M}olecular {D}ynamics: With {D}eterministic and
  {S}tochastic {N}umerical {M}ethods}.
\newblock \bibinfo{journal}{Interdisciplinary applied mathematics}
  \bibinfo{volume}{39}, \bibinfo{pages}{443}.
\bibitem[{Leimkuhler and Sachs(2022)}]{lei:22}
\bibinfo{author}{Leimkuhler, B.}, \bibinfo{author}{Sachs, M.},
  \bibinfo{year}{2022}.
\newblock \bibinfo{title}{Efficient numerical algorithms for the generalized
  {L}angevin equation}.
\newblock \bibinfo{journal}{SIAM J. Sci. Comput.} \bibinfo{volume}{44},
  \bibinfo{pages}{A364--A388}.
\bibitem[{Li et~al.(2017)Li, Lee, Darve and Karniadakis}]{li:17}
\bibinfo{author}{Li, Z.}, \bibinfo{author}{Lee, H.S.}, \bibinfo{author}{Darve,
  E.}, \bibinfo{author}{Karniadakis, G.E.}, \bibinfo{year}{2017}.
\newblock \bibinfo{title}{Computing the non-{M}arkovian coarse-grained
  interactions derived from the {M}ori--{Z}wanzig formalism in molecular
  systems: application to polymer melts}.
\newblock \bibinfo{journal}{J. Chem. Phys.} \bibinfo{volume}{146},
  \bibinfo{pages}{014104}.
\bibitem[{Mitterwallner et~al.(2020)Mitterwallner, Schreiber, Daldrop,
  R{\"a}dler and Netz}]{mi:20}
\bibinfo{author}{Mitterwallner, B.G.}, \bibinfo{author}{Schreiber, C.},
  \bibinfo{author}{Daldrop, J.O.}, \bibinfo{author}{R{\"a}dler, J.O.},
  \bibinfo{author}{Netz, R.R.}, \bibinfo{year}{2020}.
\newblock \bibinfo{title}{Non-{M}arkovian data-driven modeling of single-cell
  motility}.
\newblock \bibinfo{journal}{Phys. Rev. E} \bibinfo{volume}{101},
  \bibinfo{pages}{032408}.
\bibitem[{Ness et~al.(2015)Ness, Stella, Lorenz and Kantorovich}]{ne:15}
\bibinfo{author}{Ness, H.}, \bibinfo{author}{Stella, L.},
  \bibinfo{author}{Lorenz, C.}, \bibinfo{author}{Kantorovich, L.},
  \bibinfo{year}{2015}.
\newblock \bibinfo{title}{Applications of the generalized {L}angevin equation:
  {T}owards a realistic description of the baths}.
\newblock \bibinfo{journal}{Phys. Rev. B} \bibinfo{volume}{91},
  \bibinfo{pages}{014301}.
\bibitem[{Nualart(2006)}]{nua:06}
\bibinfo{author}{Nualart, D.}, \bibinfo{year}{2006}.
\newblock \bibinfo{title}{The Malliavin {C}alculus and {R}elated {T}opics}.
  volume \bibinfo{volume}{1995}.
\newblock \bibinfo{publisher}{Springer}.
\bibitem[{Pavliotis(2014)}]{pav:14}
\bibinfo{author}{Pavliotis, G.A.}, \bibinfo{year}{2014}.
\newblock \bibinfo{title}{Stochastic {P}rocesses and {A}pplications:
  {D}iffusion {P}rocesses, the {F}okker-Planck and {L}angevin {E}quations}.
  volume~\bibinfo{volume}{60}.
\newblock \bibinfo{publisher}{Springer}.
\bibitem[{Pigato(2018)}]{piga:18}
\bibinfo{author}{Pigato, P.}, \bibinfo{year}{2018}.
\newblock \bibinfo{title}{Tube estimates for diffusion processes under a weak
  h{\"o}rmander condition}.
\newblock \bibinfo{journal}{Ann. Inst. H. Poincaré Probab. Statist.}
  \bibinfo{volume}{54}, \bibinfo{pages}{299--342}.
\bibitem[{Pilipovic et~al.(2024)Pilipovic, Samson and Ditlevsen}]{pil:24}
\bibinfo{author}{Pilipovic, P.}, \bibinfo{author}{Samson, A.},
  \bibinfo{author}{Ditlevsen, S.}, \bibinfo{year}{2024}.
\newblock \bibinfo{title}{Parameter estimation in nonlinear multivariate
  stochastic differential equations based on splitting schemes.}
\newblock \bibinfo{journal}{A preprint. hal-04457892} .
\bibitem[{Pokern et~al.(2009)Pokern, Stuart and Wiberg}]{poke:09}
\bibinfo{author}{Pokern, Y.}, \bibinfo{author}{Stuart, A.M.},
  \bibinfo{author}{Wiberg, P.}, \bibinfo{year}{2009}.
\newblock \bibinfo{title}{Parameter estimation for partially observed
  hypoelliptic diffusions}.
\newblock \bibinfo{journal}{J. R. Stat. Soc., B: Stat. Methodol.}
  \bibinfo{volume}{71}, \bibinfo{pages}{49--73}.
\bibitem[{Samson and Thieullen(2012)}]{sam:12}
\bibinfo{author}{Samson, A.}, \bibinfo{author}{Thieullen, M.},
  \bibinfo{year}{2012}.
\newblock \bibinfo{title}{A contrast estimator for completely or partially
  observed hypoelliptic diffusion}.
\newblock \bibinfo{journal}{Stoch. Process. their Appl.} \bibinfo{volume}{122},
  \bibinfo{pages}{2521--2552}.
\bibitem[{Spannaus et~al.(2022)Spannaus, Papamarkou, Erwin and
  Christian}]{sir:22}
\bibinfo{author}{Spannaus, A.}, \bibinfo{author}{Papamarkou, T.},
  \bibinfo{author}{Erwin, S.}, \bibinfo{author}{Christian, J.B.},
  \bibinfo{year}{2022}.
\newblock \bibinfo{title}{Inferring the spread of {COVID}-19: the role of
  time-varying reporting rate in epidemiological modelling}.
\newblock \bibinfo{journal}{Sci. Rep.} \bibinfo{volume}{12},
  \bibinfo{pages}{1--12}.
\bibitem[{Tsai and Chan(2000)}]{ts:00}
\bibinfo{author}{Tsai, H.}, \bibinfo{author}{Chan, K.}, \bibinfo{year}{2000}.
\newblock \bibinfo{title}{Testing for nonlinearity with partially observed time
  series}.
\newblock \bibinfo{journal}{Biometrika} \bibinfo{volume}{87},
  \bibinfo{pages}{805--821}.
\bibitem[{Uchida and Yoshida(2012)}]{uchi:12}
\bibinfo{author}{Uchida, M.}, \bibinfo{author}{Yoshida, N.},
  \bibinfo{year}{2012}.
\newblock \bibinfo{title}{Adaptive estimation of an ergodic diffusion process
  based on sampled data}.
\newblock \bibinfo{journal}{Stoch. Process. their Appl.} \bibinfo{volume}{122},
  \bibinfo{pages}{2885--2924}.
\bibitem[{Vroylandt et~al.(2022)Vroylandt, Gouden{\`e}ge, Monmarch{\'e},
  Pietrucci and Rotenberg}]{vr:22}
\bibinfo{author}{Vroylandt, H.}, \bibinfo{author}{Gouden{\`e}ge, L.},
  \bibinfo{author}{Monmarch{\'e}, P.}, \bibinfo{author}{Pietrucci, F.},
  \bibinfo{author}{Rotenberg, B.}, \bibinfo{year}{2022}.
\newblock \bibinfo{title}{Likelihood-based non-{M}arkovian models from
  molecular dynamics}.
\newblock \bibinfo{journal}{Proc. Natl. Acad. Sci. U. S. A.}
  \bibinfo{volume}{119}, \bibinfo{pages}{e2117586119}.

\end{thebibliography}
 

\end{document}